\documentclass[12pt]{article} 
\usepackage[sectionbib]{natbib}
\usepackage{array,epsfig,fancyheadings,rotating}
\usepackage[]{hyperref}  
\usepackage{sectsty, secdot}
\usepackage{subfigure}
\usepackage{booktabs}
\usepackage{tablefootnote}
\sectionfont{\fontsize{12}{14pt plus.8pt minus .6pt}\selectfont}
\renewcommand{\theequation}{\thesection\arabic{equation}}
\subsectionfont{\fontsize{12}{14pt plus.8pt minus .6pt}\selectfont}

\usepackage{amsmath}
\usepackage{amssymb}
\usepackage{amsfonts}
\usepackage{multirow}
\usepackage{amsthm}

\newtheorem{theorem}{Theorem}
\newtheorem{lemma}{Lemma}

\theoremstyle{definition}

\newtheorem{remark}{Remark}
\usepackage[usenames,dvipsnames]{color}
\usepackage{colortbl}
\definecolor{mygray}{gray}{.8}

\begin{document}


\renewcommand{\baselinestretch}{2}

\markright{ \hbox{\footnotesize\rm Statistica Sinica
}\hfill\\[-13pt]
\hbox{\footnotesize\rm
}\hfill }

\markboth{\hfill{\footnotesize\rm FIRSTNAME1 LASTNAME1 AND FIRSTNAME2 LASTNAME2} \hfill}
{\hfill {\footnotesize\rm FILL IN A SHORT RUNNING TITLE} \hfill}

\renewcommand{\thefootnote}{}
$\ $\par


\fontsize{12}{14pt plus.8pt minus .6pt}\selectfont \vspace{0.8pc}
\centerline{\large\bf DOUBLY ROBUST ESTIMATION FOR}
\vspace{2pt}
\centerline{\large\bf CONDITIONAL TREATMENT EFFECT:}
\vspace{2pt}
\centerline{\large\bf A STUDY ON ASYMPTOTICS}
\vspace{.4cm}
\centerline{Chuyun Ye$^{1}$, Keli Guo$^{2}$ and Lixing Zhu$^{1,2}$}
\vspace{.4cm}
\centerline{\it $^1$ Beijing Normal University, Beijing, China}
\centerline{\it $^2$ Hong Kong Baptist University, Hong Kong}
 \vspace{.55cm} \fontsize{9}{11.5pt plus.8pt minus.6pt}\selectfont


\begin{quotation}
\noindent {\it Abstract:}
In this paper, we apply doubly robust approach  to estimate, when some covariates are given, the conditional average treatment effect under parametric, semiparametric and nonparametric structure of the nuisance propensity score and outcome regression models. We then conduct a systematic study on the asymptotic distributions of nine estimators with different combinations of  estimated propensity score and outcome regressions. The study covers the asymptotic properties with all models correctly specified; with either propensity score  or outcome regressions locally / globally misspecified; and with all models locally / globally misspecified. The asymptotic variances are compared and  the asymptotic bias correction under model-misspecification is discussed. The phenomenon that the asymptotic variance, with model-misspecification, could sometimes be even smaller than that with all models correctly specified is explored. We also conduct a numerical study to examine the theoretical results.

\vspace{9pt}
\noindent {\it Key words and phrases:}
Asymptotic variance, Conditional average treatment effect, Doubly robust estimation.
\par
\end{quotation}\par

\def\thefigure{\arabic{figure}}
\def\thetable{\arabic{table}}

\renewcommand{\theequation}{\thesection.\arabic{equation}}

\fontsize{12}{14pt plus.8pt minus .6pt}\selectfont

\section{Introduction} \label{introduction}


 To explore the heterogeneity of treatment effect under Rubin's protential outcome framework (\cite{rosenbaum1983central}) to reveal the casuality of a treatment, conditional average treatment effect (CATE)  is useful, which is conditional on some covariates of interest. See \cite{abrevaya2015estimating}  as an example. 
\cite{shi2019sparse} showed that the existence of optimal individualized treatment regime (OITR) has a close connection with CATE.

To estimate CATE, there are some standard approaches available in the literature. When either propensity score function or outcome regression functions or both are unknown, we need to estimate them first such that we can then estimate the CATE function. Regard these functions as nuisance models. \cite{abrevaya2015estimating} used the propensity score-based (PS-based) estimation under parametric (P-IPW) and nonparametric structure (N-IPW), and showed that N-IPW is asymptotically more efficient than P-IPW. \cite{zhou2020on} suggested  the PS-based estimation under a semiparametric dimension reduction structure (S-IPW) to show the advantage of semiparametric estimation and \cite{li2020outcome} considered outcome regression-based (OR-based) estimation under parametric (P-OR), semiparametric (S-OR) and nonparametric structure (N-OR) to derive their asymptotic properties and suggested also the use of semiparametric method. Both of the works together give an estimation efficiency comparison  between PS-based and OR-based estimators. A clear asymptotic efficiency ranking was shown by \cite{li2020outcome} when the propensity score and outcome regression models are all correctly specified and  the underlying nonparametric models is sufficiently smooth such that, with delicately selecting bandwidths and kernel functions, the nonparametric estimation can achieve sufficiently fast rates of convergence:
\begin{align}
	\overbrace{{\tiny\text{O-OR}} \cong {\tiny\text{P-OR}} \preceq {\tiny\text{S-OR}} \preceq {\tiny\text{N-OR}}}^{\text{OR-based estimators}} \cong \overbrace{{\tiny\text{N-IPW}} \preceq {\tiny\text{S-IPW}} \preceq {\tiny\text{P-IPW}} \cong {\tiny\text{O-IPW}}}^{\text{PS-based CATE estimators}}
\end{align}\label{1.1}
where $ A \preceq B$ denotes the asymptotic efficiency advantage, with smaller variance, of $A$ over $B$, $ A \cong B$ the efficiency equivalence and O-OR and O-IPW stand for OR-based and PS-based estimator respectively assuming the nuisance models are known with no need to estimate.

As well known, the doubly robust (DR) method   that was first suggested as  the augmented inverse probability weighting (AIPW) estimation proposed by  \cite{robins1994estimation}. Later developments provide the estimation consistency (\cite{scharfstein1999adjusting}) for more general doubly robust  estimation, not restricted to AIPW, that  even has  one  misspecified in the two involved models. For further discussion and introduction on  DR estimation, readers can refer to, as an example, \cite{seaman2018introduction}. Like \cite{abrevaya2015estimating}, \cite{lee2017doubly} brought up a two-step AIPW estimator of CATE also under parametric structure. For the cases with  high-dimensional covariate, \cite{fan2019estimation} and \cite{zimmert2019nonparametric} combined such an estimator with statistical learning.

In the current paper, we focus on investigating the asymptotic efficiency comparisons among nine doubly robust estimators  under parametric, semiparametric dimension reduction and nonparametric structure.   To this end, we will give a systematic study to provide insight into which combinations may have merit in an asymptotic sense and in practice, which ones would be worth of recommendation for use. We also further consider the asymptotic efficiency when nuisance models are globally or locally misspecified, which will be defined later. Roughly  speaking, local misspecification means that  misspecified model can converge, at a certain rate, to the corresponding correctly specified model as the sample size $n$ goes to infinity, while globally misspecified model cannot. Denote $c_n$, $d_{1n}$ and $d_{0n}$ respectively  the departure degrees of used models to the corresponding correctly specified models, and $V_i(x_1)$ for $i=1,2,3,4$, which will be clarified in Theorems 1, 2, 3 and 5 respectively, of the asymptotic variance functions of $x_1$ for all nine estimators in difference scenarios. Here $V_1(x_1)$ is the asymptotic variance when all models are correctly specified, which is regarded as a benchmark for comparisons. We have that $V_1(x_1)\le V_3(x_1)$, but $V_2(x_1)$ and $V_4(x_1)$ are not necessarily larger than $V_1(x_1)$.   Here we display main findings in this paper.

\begin{itemize}
	\item When all nuisance models are correctly specified, and the tuning parameters including the bandwidths in nonparametric estimations are delicately selected, the asymptotic variances are  all equal to $V_1(x_1)$.   Write all DR estimators as $DRCATE$. Together with (1.1), the asymptotic efficiency ranking is as:
\begin{align*}
	\overbrace{{\tiny\text{O-OR}} \cong {\tiny\text{P-OR}} \preceq {\tiny\text{S-OR}} \preceq {\tiny\text{N-OR}}}^{\text{OR-based estimators}} \cong {\tiny \text{DRCATE}} \cong \overbrace{{\tiny\text{N-IPW}} \preceq {\tiny\text{S-IPW}} \preceq {\tiny\text{P-IPW}} \cong  {\tiny\text{O-IPW}}}^{\text{PS-based CATE estimators}}
\end{align*}
	\item If only one  of the nuisance models, either propensity score or outcome regressions, is (are) misspecified, the estimators remain unbiased as expectably. But globally misspecified outcome regressions or propensity score lead to asymptotic variance changes. We can  give examples of propensity score to show that the variance can be even smaller than that with correctly specified models.
 Further, when the nuisance models are locally misspecified, the asymptotic efficiency remains the same as that with no misspecification.
	\item Further, when all nuisance models are globally misspecified, we need to take care of estimation bias.  When the misspecifications are all local, but the convergence rates $c_n d_{1n}$ and $c_n d_{0n}$ are all faster than the  convergence rate of nonparametric estimation that will be specified later, the asymptotic distributions  remain unchanged. 
\end{itemize}

To give a quick access to the results about the asymptotic variances, we present a summary in Table~\ref{summary table}. 
Denote $PS(P)$, $PS(N)$ and $PS(S)$ as estimators with parametrically, nonparametrically and semiparametrically estimated PS function respectively, $OR(P)$, $OR(N)$ and $OR(S)$ as estimators with parametrically,  nonparametrically and semiparametrically estimated OR functions respectively. Dark cells mean no such  combinations.

\begin{table}[t!]
	\label{summary table}
	\newcommand{\tabincell}[2]{\begin{tabular}{@{}#1@{}}#2\end{tabular}}
	\caption{Asymptotic variance result summary}
	\resizebox{\linewidth}{!}{
		\begin{tabular}{|c|c|c|c|c|c|}
			\hline
			Combination  & \tabincell{c}{All \\ Correctly \\ specified} & \tabincell{c}{Globally \\ Misspecified \\ PS} & \tabincell{c}{Locally \\ Misspecified \\ PS} & \tabincell{c}{Globally \\ Misspecified \\ OR} & \tabincell{c}{Locally \\ Misspecified \\ OR} \\
			\hline
			$PS(P) + OR(P)$ & $V_1(x_1)$ & \tabincell{c}{$V_2(x_1)$ \\ (Not necessarily \\ enlarged \footnotemark[1])} & $V_1(x_1)$ & \tabincell{c}{$V_3(x_1)$ \\ (Enlarged)} & $V_1(x_1)$ \\
			\hline
			$PS(P) + OR(N)$ & $V_1(x_1)$ & $V_1(x_1)$ & $V_1(x_1)$ & \cellcolor{mygray} & \cellcolor{mygray}\\
			\hline
			$PS(N) + OR(P)$ & $V_1(x_1)$ & \cellcolor{mygray} & \cellcolor{mygray} & $V_1(x_1)$ & $V_1(x_1)$\\
			\hline
			$PS(N) + OR(N)$ & $V_1(x_1)$ & \cellcolor{mygray} & \cellcolor{mygray} & \cellcolor{mygray} & \cellcolor{mygray}\\
			\hline
			$PS(P) + OR(S)$ & $V_1(x_1)$ &  \tabincell{c}{$V_2(x_1)$ \\ (Not necessarily \\ enlarged)} & $V_1(x_1)$ & \cellcolor{mygray} & \cellcolor{mygray} \\
			\hline
			$PS(S) + OR(P)$ & $V_1(x_1)$ & \cellcolor{mygray} & \cellcolor{mygray} & \tabincell{c}{$V_3(x_1)$ \\ (Enlarged)} & $V_1(x_1)$\\
			\hline
			$PS(S) + OR(N)$ & $V_1(x_1)$ & \cellcolor{mygray} & \cellcolor{mygray} & \cellcolor{mygray} & \cellcolor{mygray} \\
			\hline
			$PS(N) + OR(S)$ & $V_1(x_1)$ & \cellcolor{mygray} & \cellcolor{mygray} & \cellcolor{mygray} & \cellcolor{mygray} \\
			\hline
			$PS(S) + OR(S)$ & $V_1(x_1)$ & \cellcolor{mygray} & \cellcolor{mygray} & \cellcolor{mygray} & \cellcolor{mygray} \\
			\hline \hline
			Combination & \cellcolor{mygray} & \tabincell{c}{All \\ Globally \\ Misspecified} & \tabincell{c}{All \\ Locally \\ Misspecified} & \tabincell{c}{Globally Misspecified PS \\ + \\ Locally Misspecified OR } & \tabincell{c}{Locally  Misspecified  PS \\ + \\ Globally  Misspecified OR} \\
			\hline
			$PS(P) + OR(P)$ & \cellcolor{mygray}
& \tabincell{c}{ Biased + $V_4(x_1)$ \\ (Not necessarily \\ enlarged variance) } & $V_1(x_1)$ & \tabincell{c}{$V_2(x_1)$ \\ (Not necessarily \\ enlarged)} & \tabincell{c}{$V_3(x_1)$ \\ (Enlarged)} \\
			\hline
		\end{tabular}}
\end{table}

The remaining parts of this article are organized as follows. We first describe the Rubin's potential outcome framework and the relevant notations in Section~2. Section~3 contains a general two-step estimation of CATE, while Section~4 describes the corresponding asymptotic properties under different situations. Section~5 presents the results of Monte Carlo simulations and  Section~6 includes some concluding remarks. We would like to point out that such comparisons do  not mean the estimations that are of asymptotic efficiency advantage are always worthwhile to recommend because, particularly, the nonparametric-based estimations may  have severe difficulties to handle high- even moderate-dimensional models in practice. But the comparisons can provide a good insight into the nature of various estimations such that the practitioners can have a relatively complete picture about them and have idea for when and how to use these estimations.

\section{Framework and Notation} \label{framework and notation}

For any individual, datum $W=(X^{\top},Y,D)^{\top}$ is observable, including the observed effect $Y$, the treatment status $D$, and the $p$-dimensional covariates $X$. $D=1$ implies that the individual is treated, and $D=0$ means untreated.  Denote $Y(1)$ and $Y(0)$ as the potential outcomes with and without treatment, respectively. The observed effect $Y$ can be expressed as $Y=DY(1)+(1-D)Y(0)$.
Denote that $p(X) = P(D=1 | X), m_1(X) = E(Y(1) | X), m_0(X) = E(Y(0) | X)$  as propensity score function and outcome regression functions. The following conditions are commonly used when we discuss the potential outcome framework.
\begin{description}
	\item (C1) (Sampling distribution) $\{ W_i \}_{i=1}^n$ is a set of identically distributed samples.
	\item (C2) (Ignorability condition)
	\begin{description}
		\item (\romannumeral1) (Unconfoundedness) $(Y(1), Y(0)) \perp D | X$
		\item (\romannumeral2) Denote $\mathcal{X}$ as the support of $X$, where $\mathcal{X}$ is a Cartersian product of compact intervals. For any $x \in \mathcal{X}$, $p(x)$ is bounded away from 0 and 1.
	\end{description}
\end{description}


Denote $\tau(x_1)$ as CATE:
\begin{align*}
	\tau(x_1) = E [ Y(1) - Y(0) | X_1 = x_1 ]
\end{align*}
where $X_1$ is a strict subset of $X$. That is, $X_1$ is a $k$-dimension covariate, and $k < p$. Also denote $f(x_1)$ as the density function of $X_1$.

\section{Doubly robust estimation} \label{proposed}

Rewrite  $\tau(x_1)$ as
\begin{align}\label{3.3}
	\tau(x_1) & = E \left\{ \left. \left[ m_1(X) - m_0(X) \right] \right| X_1 = x_1 \right\} \nonumber\\
	& = E \left\{ \left. \left[ \frac{DY}{p(X)} - \frac{(1-D)Y}{1-p(X)} \right] \right| X_1 = x_1 \right\} \nonumber\\
	& = E \left\{ \left. \left[ \frac{D}{p(X)}[Y - m_1(X)] - \frac{1 - D}{1 - p(X)}[Y - m_0(X)] + m_1(X) - m_0(X) \right] \right| X_1 = x_1 \right\}
\end{align}

 The first two equations in (\ref{3.3}) show how OR and PS method work for estimating CATE.  The third equation in (\ref{3.3}) is an essential expression to construct a doubly robust estimator of $\tau(x_1)$. Under which, we propose a two-step estimation.
In the first step, we estimate the  function in (\ref{3.3}):
\begin{align*}
	\frac{D}{p(X)}[Y - m_1(X)] - \frac{1 - D}{1 - p(X)}[Y - m_0(X)] + m_1(X) - m_0(X).
\end{align*}
       To study the influence from  estimating the nuisance functions,  $p(X)$ and  $m_1(X)$, $m_0(X)$ under parametric, nonparametric, and semiparametric dimension reduction framework, we will construct the corresponding estimations below.

After this, we can then estimate the conditional expectation given $x_1$.  This is a standard nonparametric estimation. We utilize the Nadaraya-Watson type estimator to define  
the  resulting estimator:
\begin{align*}
	\widehat{\tau}(x_1) = \frac{\frac{1}{nh_1^k} \sum_{i=1}^n \left[ \frac{D_i}{\widehat{p}_i} \left( Y_i - \widehat{m}_{1i} \right) - \frac{1-D_i}{1-\widehat{p}_i} \left( Y_i - \widehat{m}_{0i} \right) + \widehat{m}_{1i} - \widehat{m}_{0i} \right] K_1 \left( \frac{X_{1i} - x_1}{h_1} \right)}{\frac{1}{nh_1^k} \sum_{i=1}^n K_1 \left( \frac{X_{1i} - x_1}{h_1} \right)},
\end{align*}
where $K_1(u)$ is a kernel function of order $s_1$, which is $s*$ times continuously differentiable, and $h_1$ is the corresponding bandwidth and $\widehat{p}_i, \widehat{m}_{1i}, \widehat{m}_{0i}$ denote the estimators of $p(X_i),m_1(X_i),m_0(X_i)$ respectively, which are general notations and have different formulas under different model structures.

We now consider the estimations of the nuisance functions. Under the parametric structures with $\widetilde{p}(x; \beta)$, $\widetilde{m}_1(x; \gamma_1)$ and $\widetilde{m}_0(x; \gamma_0)$ as the specified parametric models of $p(x)$, $m_1(x)$ and $m_0(x)$ respectively, where $\beta$, $\gamma_1$ and $\gamma_0$ are unknown parameters. By maximum likelihood estimation, we can obtain $\widehat{\beta}$, $\widehat{\gamma_1}$ and $\widehat{\gamma_0}$ so as to have  $\widetilde{p}(X_i; \widehat{\beta})$, $\widetilde{m}_1(X_i; \widehat{\gamma_1})$ and $\widetilde{m}_0(X_i; \widehat{\gamma_0})$ as the parametric estimators. Note that the specified models are not necessarily equal to true data generate mechanism. Now we further distinguish correctly specified, globally misspecified and locally misspecified case. For all $x \in \mathcal{X}$, there exist $\beta_0,\gamma_{10},\gamma_{00}$, such that the true models have the relationship with the specified models:
\begin{align}
	p(x) & = \widetilde{p}(x;\beta_0)[1+c_n a(x)], \nonumber\\
	m_1(x) & = \widetilde{m}_1(x;\gamma_{10}) + d_{1n} b_1(x), \\
	m_0(x) & = \widetilde{m}_0(x;\gamma_{00}) + d_{0n} b_0(x). \nonumber
\end{align}\label{true}
Take propensity score function as an example. If $c_n = 0$, then the parametric propensity score model $\widetilde{p}(x;\beta_0)$ is correctly specified, otherwise, it is not. If $c_n$ converges to 0 as $n$ goes to infinity,  the parametric model is locally misspecified. If $c_n$ remains a nonzero constant,  it is a globally misspecified case. Similarly for the models with $d_{1n}$ and $d_{0n}$. Recall that $\widehat{\beta}$, $\widehat{\gamma}_1$ and $\widehat{\gamma}_0$ are the maximum likelihood estimators of the corresponding unknown parameters. Denote $\beta^*$, $\gamma_1^*$ and $\gamma_0^*$ as the limits of $\widehat{\beta}$, $\widehat{\gamma}_1$ and $\widehat{\gamma}_0$ as $n$ goes to infinity.

Under the nonparametric structure, we utilize the kernel-based nonparametric estimators as
\begin{align*}
	\widehat{p}(X_i) & = \frac{\sum_{j=1}^n D_j K_2 \left( \frac{X_j - X_i}{h_2} \right)}{\sum_{t=1}^n K_2 \left( \frac{X_t - X_i}{h_2} \right)}, \\
	\widehat{m}_1(X_i) & = \frac{\sum_{j=1}^n D_j Y_j K_3 \left( \frac{X_j - X_i}{h_3} \right)}{\sum_{t=1}^n D_t K_3 \left( \frac{X_t - X_i}{h_3} \right)}, \\
	\widehat{m}_0(X_i) & = \frac{\sum_{j=1}^n (1 - D_j) Y_j K_4 \left( \frac{X_j - X_i}{h_4} \right)}{\sum_{t=1}^n (1 - D_t) K_4 \left( \frac{X_t - X_i}{h_4} \right)}
\end{align*}
where $K_2(u)$, $K_3(u)$ and $K_4(u)$ are kernels of order $s_2 \geq d$, $s_3 \geq d$ and $s_4 \geq d$, with the corresponding bandwidths $h_2$, $h_3$ and $h_4$. The conditions on  the kernel functions and bandwidths will be listed in the supplement.

Under the semiparametric  structure on the baseline covariate $X$ for propensity score and outcome regressions, we have the following dimension reduction framework. Denote the matrix $A \in \mathbb{R}^{d \times d_2}$ such that
\begin{align} \label{SDR structure 1}
	p(X) \perp X | A^{\top} X ,
\end{align}
where $d_2 \le d$. The $A$ spanned space $\mathcal{S}_{E(D|X)}$ is called the central mean subspace if it is the intersection of all subspaced spanned by all $A$ satisfy the above conditional independence. The dimension of $\mathcal{S}_{E(D|X)}$ is called the structural dimension that is often  smaller than  or equal to $d_2$. Without confusion, still write it  as $d_2$.   Formula (\ref{SDR structure 1}) implies that $p(X) = E(D | X) = E(D | A^{\top} X) := g(A^{\top} X)$. Note that a nonparametric estimation of $p(X)$ may have very slow rate of convergence when $p$ is large. However, under (\ref{SDR structure 1}) we can estimate the matrix $A$ first to reduce the dimension $d$ to $d_2$, the nonparametric estimation of  $E(D | A^{\top} X)$ can achieve a faster rate of  convergence. The semiparametric estimator $p(X_i)$ is then defined as, when $A$ is root-$n$ consistently by an estimator $\widehat{A}$,
\begin{align*}
	\widehat{g}(\widehat{A}^{\top} X_i) = \frac{\sum_{j=1}^n D_j K_5 \left( \frac{\widehat{A}^{\top} X_j - \widehat{A}^{\top} X_i}{h_5} \right)}{\sum_{t=1}^n K_5 \left( \frac{\widehat{A}^{\top} X_t - \widehat{A}^{\top} X_i}{h_5} \right)}.
\end{align*}
Similarly, for regression models, denote  matrixes $B_1 \in \mathbb{R}^{d \times d_1}$ and $B_0 \in \mathbb{R}^{d \times d_0}$, such that
\begin{align}
	E(Y(1) | X) \perp X | B_1^{\top} X , \nonumber\\ 
	E(Y(0) | X) \perp X | B_0^{\top} X . \label{SDR structure 3}
\end{align}
The corresponding dimension reduction subspaces are called the central mean subspaces (see Cook and Li 2002). Thus, $m_1(X) = E(Y(1)|X) = E(Y(1)|B_1^{\top} X) := r_1(B_1^{\top} X)$  and $m_0(X) = E(Y(0)|X) = E(Y(0)|B_0^{\top} X) := r_0(B_0^{\top} X)$. The semiparametric estimators $m_1(X_i)$ and $m_0(X_i)$ are defined as, with $\widehat{B}_i$ being the estimators of $B_i$, $i=0,1$,
\begin{align*}
	\widehat{r}_1(\widehat{B}_1^{\top} X_i) & = \frac{\sum_{j=1}^n D_j Y_j K_6 \left( \frac{\widehat{B}_1^{\top} X_j - \widehat{B}_1^{\top} X_i}{h_6} \right)}{\sum_{t=1}^n D_t K_6 \left( \frac{\widehat{B}_1^{\top} X_t - \widehat{B}_1^{\top} X_i}{h_6} \right)}, \\
	\widehat{r}_0(\widehat{B}_0^{\top} X_i) & = \frac{\sum_{j=1}^n D_j Y_j K_7 \left( \frac{\widehat{B}_0^{\top} X_j - \widehat{B}_0^{\top} X_i}{h_7} \right)}{\sum_{t=1}^n D_t K_7 \left( \frac{\widehat{B}_1^{\top} X_t - \widehat{B}_1^{\top} X_i}{h_7} \right)}
\end{align*}
where $K_5(u)$, $K_6(u)$ and $K_7(u)$ are kernels of order $s_5 \geq d$, $s_6 \geq d$ and $s_7 \geq d$, with the corresponding bandwidths $h_5$, $h_6$ and $h_7$.

\section{Asymptotic Properties} \label{asymptotic}

 Define the following functions
\begin{align*}
	\Psi_1(X,Y,D) &:= \frac{D[Y-m_1(X)]}{p(X)} - \frac{(1-D)[Y-m_0(X)]}{1-p(X)} + m_1(X) - m_0(X), \\
	\Psi_2(X,Y,D) & :=\frac{D\{Y-m_{1}(X)\}}{\widetilde{p}(X;\beta^*)}-\frac{(1-D)\{Y-m_{0}(X)\}}{1-\widetilde{p}(X;\beta^*)}+m_{1}(X)-m_{0}(X), \\
	\Psi_3(X,Y,D) & :=\frac{D\{Y-\widetilde{m}_{1}(X;\gamma_1^*)\}}{p(X)}-\frac{(1-D)\{Y-\widetilde{m}_{0}(X;\gamma_0^*)\}}{1-p(X)} + \widetilde{m}_1(X;\gamma_1^*) - \widetilde{m}_0(X;\gamma_0^*), \\
	\Psi_4(X,Y,D) & :=\frac{D\{Y-\widetilde{m}_{1}(X;\gamma_1^*)\}}{\widetilde{p}(X;\beta^*)}-\frac{(1-D)\{Y-\widetilde{m}_{0}(X;\gamma_0^*)\}}{1-\widetilde{p}(X;\beta^*)} + \widetilde{m}_1(X;\gamma_1^*) - \widetilde{m}_0(X;\gamma_0^*).
\end{align*}

\subsection{The Cases With No Model Misspecification} \label{no model is misspecified}

 The following theorem shows all asymptotic distributions of the estimators are identical.
\begin{theorem}
\label{thm all correct} Suppose Conditions (C1) -- (C6), (A1), (A2) and (B1) are satisfied for $s^* \geq s_2 \geq d$, $s^* \geq s_3 \geq d$,  $s^* \geq s_4 \geq d$, $s^* \geq s_5 \geq d_2$, $s^* \geq s_6 \geq d_1$, $s^* \geq s_7 \geq d_0$, and formulas \ref{SDR structure 1} and  
\ref{SDR structure 3}
hold. Then, for each point $x_1$, we have
\begin{align*}
	\sqrt{nh_1^k} \left[ \widehat{\tau}(x_1)-\tau(x_1) \right] = & \frac{1}{\sqrt{nh_1^k}}\frac{1}{f(x_1)} \sum_{i=1}^n [\Psi_1(X_i,Y_i,D_i)-\tau(x_1)] K_1 \left( \frac{X_{1i}-x_1}{h_1} \right)\\
& + o_p(1),
\end{align*}and
\begin{align*}
	\sqrt{nh_1^k} \left[ \widehat{\tau}(x_{1})-\tau(x_{1}) \right] \xrightarrow{d} N \left( 0, V_1(x_1) \right),
\end{align*}
where
\begin{align*}
	V_1(x_1) = \frac{\sigma_1^2(x_1) \int K_1^2(u) du }{f(x_1)}, \quad 	\sigma_1^2(x_1) = E \left\{ \left. [\Psi_1(X,Y,D) - \tau(x_1)]^2 \right| X_1 = x_1 \right\}.
\end{align*}
\end{theorem}



\subsection{The Cases With Misspecified Models}

 Now we discuss the asymptotic behaviours of the proposed estimators if either  outcome regression models or propensity score model is (are) misspecified. The following results show how  global misspecification affects the asymptotic properties.

\begin{theorem}
\label{globally misspecified propensity score} Assume that the propensity score is globally misspecified in which $c_n=C$ is a nonzero constant. Suppose conditions (C1) -- (C6), (A1), (A2) and (B1) are satisfied for $s^* \geq s_3 \geq d$, $s^* \geq s_4 \geq d$, $s^* \geq s_6 \geq d_1$, $s^* \geq s_7 \geq d_0$, $s_6 < (2 s_6 + k) (d - d_1)$, $s_7 < (2 s_7 + k) (d - d_0)$.\\
1). When the outcome regression functions are estimated nonparametrically, then, for each value $x_1$, we have
\begin{align*}
	\sqrt{nh_1^k} \left[ \widehat{\tau}(x_{1})-\tau(x_{1}) \right] \xrightarrow{d} N \left( 0,V_1(x_1)\right).
\end{align*}
2). When the outcome regression functions have  a dimension reduction structure specified in (\ref{SDR structure 3})  or are correctly specified with $d_{1n} = d_{0n} = 0$ with parametric estimation, for each value $x_1$, the asymptotic distributions are identical:
\begin{align*}
	\sqrt{nh_1^k} \left[ \widehat{\tau}(x_{1})-\tau(x_{1}) \right] \xrightarrow{d} N \left( 0,V_2(x_1) \right),
\end{align*}
where
\begin{align*}
	V_2(x_1) = \frac{\sigma_2^2(x_1) \int K_1^2(u) du}{f(x_1)}, \quad \sigma_2^2(x_1) = E \left\{ \left. [\Psi_2(X,Y,D) - \tau(x_1)]^2 \right| X_1 = x_1 \right\}.
\end{align*}
\end{theorem}

Now we consider the cases with global misspecification of the outcome regression models.
\begin{theorem}
\label{globally misspecified outcome regression} Assume that the outcome regression models are globally misspecified with  fixed nonzero constants $d_{1n}=d_1$ and $d_{0n}=d_2$.  Suppose conditions (C1) -- (C6), (A1), (A2) and (B1) are satisfied for $s^* \geq s_2 \geq d$, $s^* \geq s_5 \geq d_2$, $s_5 < (2 s_5 + k) (d - d_2)$.\\
1). When the propensity score is estimated nonparametrically, then, for each  $x_1$,
\begin{align*}
	\sqrt{nh_1^k} \left[ \widehat{\tau}(x_{1})-\tau(x_{1}) \right] \xrightarrow{d} N \left( 0, V_1(x_1) \right).
\end{align*}
2). When the propensity score has  a dimension reduction structure in (\ref{SDR structure 1}) or is correctly specified with $c_n = 0$ and parametric estimation, for each value $x_1$,  the asymptotic distributions are identical:
\begin{align*}
	\sqrt{nh_1^k} \left[ \widehat{\tau}(x_{1})-\tau(x_{1}) \right] \xrightarrow{d} N \left( 0, V_3(x_1) \right),
\end{align*}
where
\begin{align*}
	V_3(x_1) = \frac{\sigma_3^2(x_1) \int K_1^2(u) du }{f(x_1)}, \quad \sigma_3^2(x_1) = E \left\{ \left. [\Psi_3(X,Y,D) - \tau(x_1)]^2 \right| X_1 = x_1 \right\}.
\end{align*}
\end{theorem}

\begin{remark}
\label{global misspecification - variance}
By some  calculations, we can obtain in Proposition~\ref{proposition: sigma} below in Section~4.4 that $ \sigma_1^2(x_1)\le \sigma_3^2(x_1) $, while the analogy does not hold between $\sigma_2^2(x_1) $ and $\sigma_1^2(x_1)$. That is, the asymptotic variance of the proposed estimator inflates when the outcome regression models are misspecified, and the propensity score model is parametrically estimated (correctly specified) or semiparametrically estimated. However, whether the asymptotic variance gets larger with a misspecified propensity score model is model-dependent. We show the following example. Suppose that the outcome regression models are correctly specified, while the propensity score model is globally misspecified. Consider a situation that $p(x) = p_1, \widetilde{p}(x;\beta^*) = p_2$, where $p_1, p_2$ are free of $x$, and $p_1 \neq p_2$. We have
\begin{align*}
	\sigma_2^2(x_1) - \sigma_1^2(x_1) = & E \left\{ \left. \frac{p^2(X) - \widetilde{p}^2(X;\beta^*)}{\widetilde{p}^2(X;\beta^*)p(X)} Var (Y|X, D=1) \right| X_1 = x_1\right\} \\
	& + E \left\{ \left. \frac{[1 - p(X)]^2 - [1 - \widetilde{p}^2(X;\beta^*)]^2}{[1 - \widetilde{p}(X;\beta^*)]^2[1 - p(X)]} Var (Y|X, D=0) \right| X_1 = x_1 \right\} \\
	= & \frac{p_1^2 - p_2^2}{p_1 p_2^2} E \left[ \left. Var (Y|X, D=1) \right| X_1 = x_1 \right] \\
	& + \frac{(1-p_1)^2 - (1 - p_2)^2}{(1 - p_1) (1 - p_2)^2} E \left[ \left. Var (Y|X, D=0) \right| X_1 = x_1 \right].
\end{align*}

To give a clear picture, we further assume that the outcome regression models are homoscedastic that $Var (Y|X, D=1) = Var (Y|X, D=0) = \xi^2$, which is free of $X$. Then we have, $\sigma_2^2(x_1) - \sigma_1^2(x_1) = \xi^2 \left( \frac{p_1^2 - p_2^2}{p_1 p_2^2} + \frac{(1-p_1)^2 - (1 - p_2)^2}{(1 - p_1) (1 - p_2)^2} \right)$. Define the function $vd(p_1, p_2) = \left( \frac{p_1^2 - p_2^2}{p_1 p_2^2} + \frac{(1-p_1)^2 - (1 - p_2)^2}{(1 - p_1) (1 - p_2)^2} \right)$. A negative $vd(p_1,p_2)$ implies the variance shrinkage. Consider three true propensity score values $p(x) = p_1 = 0.3, 0.5, 0.7$. The following three curves of $vd(p_1,p_2)$ show how the variance inflation or shrinkage occurs.

\begin{figure}[h]
	\centering
	\subfigure[$p_1 = 0.3$]{
		\includegraphics[width=0.3 \textwidth]{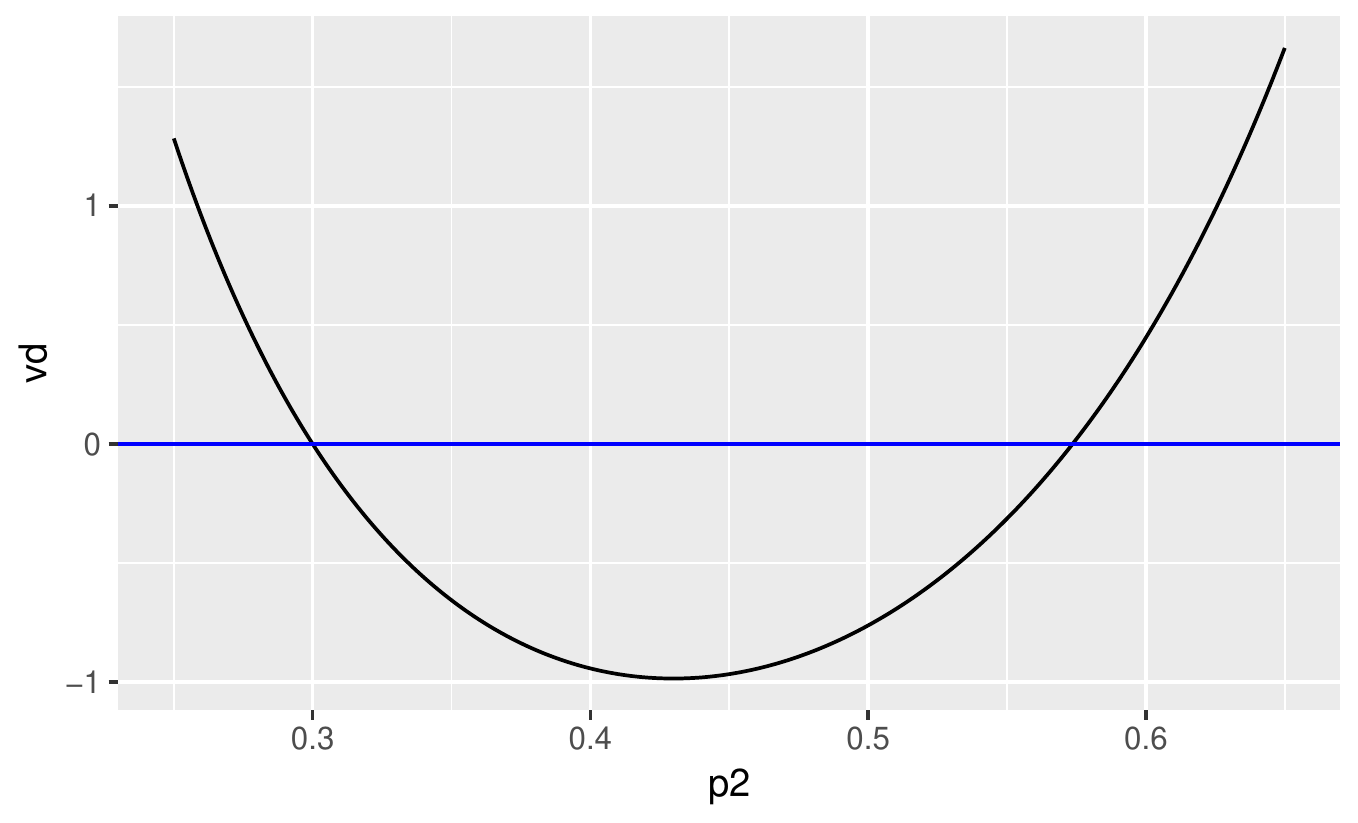}
	}
	\subfigure[$p_1 = 0.5$]{
		\includegraphics[width=0.3 \textwidth]{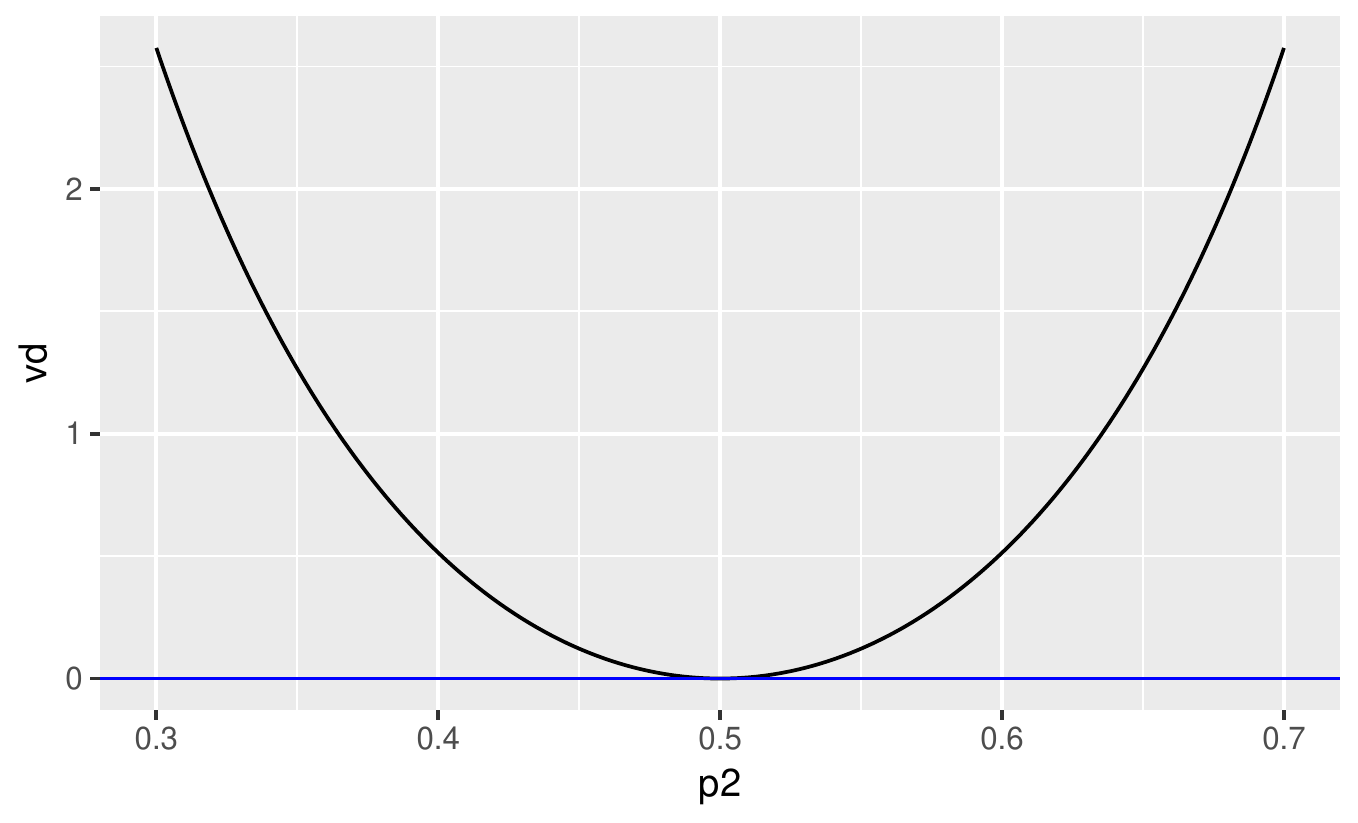}
	}
	\subfigure[$p_1 = 0.7$]{
		\includegraphics[width=0.3 \textwidth]{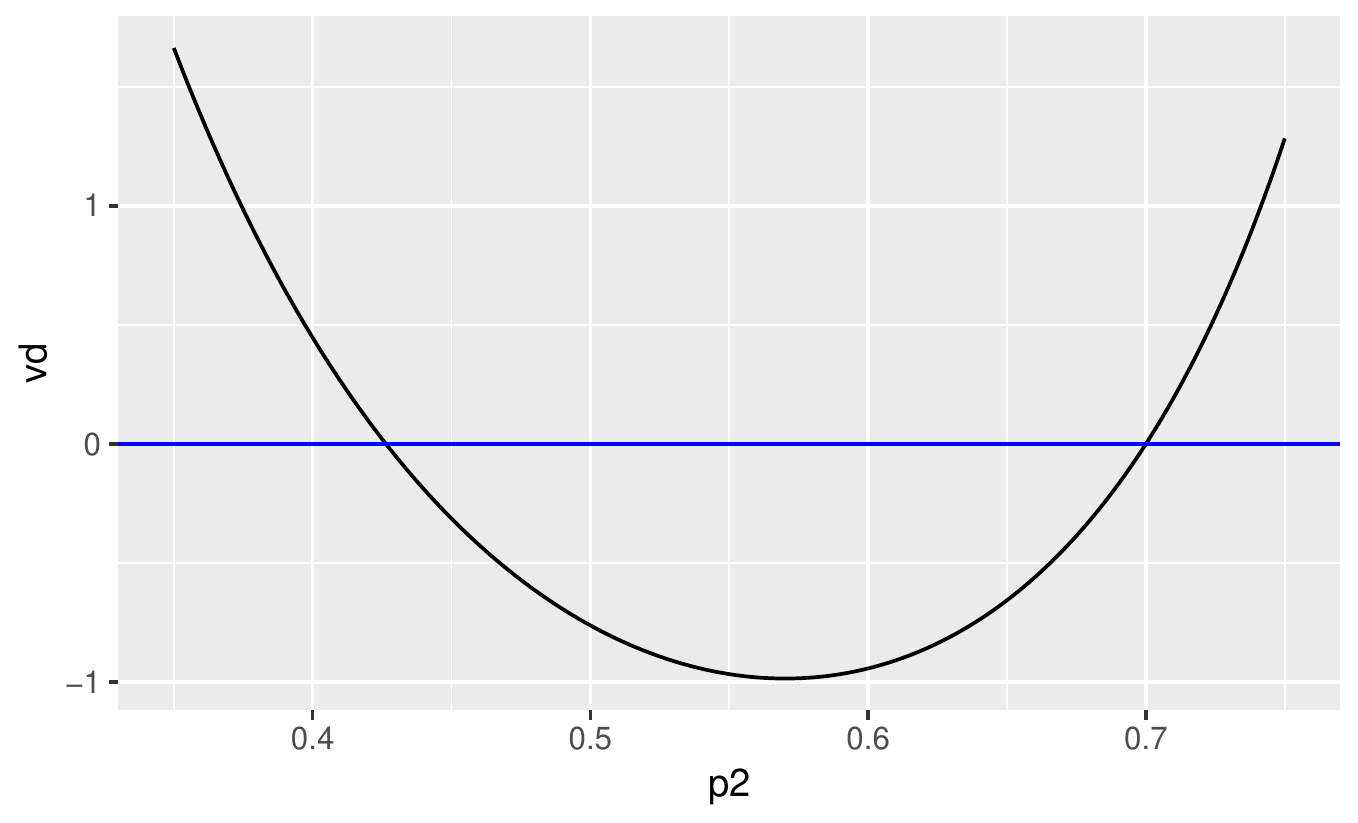}
	}
\caption{Curves of $vd(p_1,p_2)$ with different $p_1$}
\end{figure}

When $p_1 = 0.3$ or $0.7$,  appropriately overestimated propensity score may result in an asymptotic variance shrinking in some cases. When $p_1 = 0.5$, which means that every individual have an $0.5$ probability to be treated regardless of any covariates, misspecification leads to the asymptotic variance augmentation.

We can in effect obtain some more examples since $Var (Y|X, D=1)$ and $Var (Y|X, D=0)$ are not necessarily equal. Such simple examples show that when only propensity score is misspecified, augmenting or shrinking asymptotic variances are all possible. 
\end{remark}

\begin{remark}
\label{global misspecification - absorb}
Another interesting phenomenon is that once propensity score model is misspecified and outcome regressions are nonparametrically estimated, or vice versa, the asymptotic performance of the proposed estimator is identical to  that when all models are correctly specified. As nonparametric estimation takes no risk of misspecification, such an estimation procedure ``absorbs'' the influence brought by model misspecification due to the doubly robust property. But it is clear that in high-dimensional scenarios, a purely nonparametric estimation is not worthwhile to recommend. Thus, this property  mainly serves as an investigation with theoretical interest unless the dimension of the covariates is small.
\end{remark}

The results with local misspecification are stated in the following.
\begin{theorem}
\label{local misspecification} Assume that  the propensity score is locally specified with $c_n \to 0$. Suppose  conditions (C1) -- (C6), (A1), (A2) and (B1) are satisfied for $s^* \geq s_3 \geq d$, $s^* \geq s_4 \geq d$, $s^* \geq s_6 \geq d_1$, $s^* \geq s_7 \geq d_0$, $s_6 < (2 s_6 + k) (d - d_1)$, $s_7 < (2 s_7 + k) (d - d_0)$. Then, for each value $x_1$, we have
\begin{align*}
	\sqrt{nh_1^k} \left[ \widehat{\tau}(x_1)-\tau(x_1) \right] \xrightarrow{d}N \left( 0, V_1(x_1) \right).
\end{align*}
Similarly, assume that the outcome regression functions are locally misspecified with $d_{1n} \to 0$ and $d_{0n} \to 0$.  Under the same conditions as those in Theorem~\ref{local misspecification} for $s^* \geq s_2 \geq d$, $s^* \geq s_5 \geq d(2)$, $s_5 < (2 s_5 + k) (p - p(2))$. For each value $x_1$, the asymptotic distribution of $\widehat{\tau}(x_1)$ is identical to the above.
\end{theorem}

\subsection{A Further Study: All Models are Misspecified}

We study this case as it then has a non-ignorable bias in general and goes to zero unless the rate of convergence of local misspecification is sufficiently fast. Recall the definitions of $\gamma_0^*$ $\gamma_1^*$ and $\beta^*$ below (\ref{true}).

\begin{theorem}
\label{all globally misspecified}
Suppose that all models are globally misspecified with nonzero constants $c_n$, $d_{1n}$ and $d_{0n}$. Assume that conditions (C1) -- (C6) are satisfied. Then, for each value $x_1$, we have
\begin{align*}
	\sqrt{nh_1^k} \left[ \widehat{\tau}(x_{1})-\tau(x_{1}) -bias(x_1)\right] \xrightarrow{d} N \left( 0, V_4(x_1) \right),
\end{align*}
where
\begin{align*}
	bias(x_1) = &  E \left\{ \frac{\left[ m_1(X) - \widetilde{m}_1(X;\gamma_1^*) \right] \left[ p(X) - \widetilde{p}(X;\beta^*) \right]}{\widetilde{p}(X;\beta^*)} \right. \\
	 & \left. \left. - \frac{\left[ m_0(X) - \widetilde{m}_0(X;\gamma_0^*) \right] \left[ \widetilde{p}(X;\beta^*) - p(X) \right]}{1 - \widetilde{p}(X;\beta^*)} \right| X_1 = x_1 \right\},
\end{align*}
\begin{align*}
	V_4(x_1) = \frac{\sigma_4^2(x_1) \int K_1^2(u) du}{f(x_1)}, \quad \sigma_4^2(x_1) = E \left\{ \left. [\Psi_4(X,Y,D) - \widetilde{\tau}(x_1)]^2 \right| X_1 = x_1 \right\},
\end{align*}
and
\begin{align*}
	\widetilde{\tau}(x_1) = & E \left\{ \left[ \frac{D}{\widetilde{p}(X;\beta^*)}[Y - \widetilde{m}_1(X;\gamma_1^*)]  \right. \right.  \\
	& -\left. \left. \left. \frac{1 - D}{1 - \widetilde{p}(X;\beta^*)}[Y - \widetilde{m}_0(X;\gamma_0^*)] + \widetilde{m}_1(X;\gamma_1^*) - \widetilde{m}_0(X;\gamma_0^*) \right] \right| X_1 = x_1 \right\} .
\end{align*}
\end{theorem}

The following results show the importance of the convergence rates of $c_n$, $d_{1n}$ and $d_{0n}$ to zero for  bias reduction and variance change.

\begin{theorem}
\label{all locally misspecified} Under the conditions in Theorem~\ref{all globally misspecified}, when
\begin{align*}
	c_n d_{1n} & = o \left( \frac{1}{\sqrt{nh_1^k}} \right), \, \, c_n d_{0n} = o\left( \frac{1}{\sqrt{nh_1^k}} \right),
\end{align*}
then, for each  $x_1$, we have
\begin{align*}
	\sqrt{nh_1^k} \left[ \widehat{\tau}(x_{1})-\tau(x_{1}) \right] \xrightarrow{d} N \left( 0,V_1(x_1)\right).
\end{align*}
\end{theorem}

\begin{remark}
\label{local misspecification - rate}
     This theorem show that to make the bias vanished, $c_n d_{1n}$ and $c_n d_{0n}$ need to tend to zero at the rates faster than the nonparametric convergence rate, $O ( 1/\sqrt{nh_1^k} )$. Recall that  Theorems \ref{globally misspecified propensity score} and   \ref{globally misspecified outcome regression} show that when $c_n = o(1)$, then the variance  is $V_3(x_1)$; when $d_{1n}=o(1)$ and $d_{0n}=o(1)$  the variance is $V_2(x_1)$. Altogether,  when all misspecifications are local, the asymptotic variances reduces to $V_1(x_1)$. We can then further discuss four cases:\\
 1) All nuisance models are globally misspecified;
 2) All nuisance models are locally misspecified;
 3) The propensity score function is globally misspecified, and the outcome regression functions are locally misspecified;
 4) The propensity score function is locally misspecified, and the outcome regression functions are globally misspecified.

The first  is the  case exactly described in Theorem \ref{all globally misspecified}, the second shows that  if  $c_n d_{1n} = o ( 1/\sqrt{nh_1^k} )$ and $c_n d_{0n} = o ( 1/\sqrt{nh_1^k} )$, the bias term is negligible, which is the situation  in Theorem \ref{all locally misspecified}. Otherwise, the estimator is biased.
Cases 3 and 4 can be regarded as a combination of those in Theorems \ref{all globally misspecified} and \ref{all locally misspecified}. In case~3, once $d_{1n} = o ( 1/\sqrt{nh_1^k} ) $ and $d_{0n} = o ( 1/\sqrt{nh_1^k} ) $, the bias  goes to $0$, and the variance  goes to $||K_1||_2^2 \sigma_2^2(x_1)/f(x_1)$. In other words, if $d_{1n}$ and $d_{0n}$ go to $0$ at a rate faster than $O ( 1/\sqrt{nh_1^k} )$, Case 3 turns to the case in Theorem \ref{globally misspecified propensity score}. We can then also derive that if $c_n = o ( 1/\sqrt{nh_1^k} )$, Case 4 is similar to that in Theorem \ref{globally misspecified outcome regression}.
\end{remark}

\subsection{A summary on the comparison among the asymptotic Variances }

 We summarize the comparison among  the  $4$  variances $V_j(x_1)$ for $j=1,2,3,4$ as listed in Section~1. Note that the variances are $V_j(x_1) = ||K_1||_2^2 \sigma_j^2(x_1)/f(x_1)$ for $j=1,2,3,4$  and thus the comparison among them is equivalent to the comparison among  $\sigma_j^2(x_1)$ for $j=1,2,3,4$.
\begin{remark}\label{proposition: sigma}
For any $x_1$,\\
1). $\sigma_1^2(x_1)$ is not necessarily smaller than $\sigma_2^2(x_1)$ and  as shown in the example in Remark \ref{global misspecification - variance}, $\sigma_1^2(x_1)$ can be larger than $\sigma_2^2(x_1)$ for some $x_1$;\\
2). $\sigma_1^2(x_1) \le \sigma_3^2(x_1)$;\\
3). We have no definitive answer to say whether $\sigma_1^2(x_1)$ is necessarily smaller than $\sigma_4^2(x_1)$.
\end{remark}

\section{Numerical Study}

In this section, we present some Monte Carlo simulations to examine the finite sample performances of the estimators.

\subsection{Data-Generating Process}

Consider two data-generating processes (DGPs) similarly as those in \cite{abrevaya2015estimating}, the case of $d=2$ and $d=4$. Here we only consider that the conditioning covariate $X_1$ is  univariate, i.e. $k=1$. So in the simulations, $\tau(x_1) = E[Y(1) - Y(0) | X_1 = x_1]$.

{\bf Model 1}. It  is featured by a 2-dimensional unconfounded covariate, $X = (X_1, X_2)^{\top}$. In other words, $d=2$. For further information,
\begin{align*}
	X_1 = \rho_1, \qquad X_2 = ( 1 + 2 X_1 )^2 ( -1 + X_1 )^2 + \rho_2,
\end{align*}
where $\rho_1, \rho_2$ are independently identically $U(-0.5,0.5)$ distributed. The  potential outcomes and the propensity score function are given as:
\begin{align*}
	Y(1) = X_1 X_2 + \epsilon, \qquad Y(0) = 0, \\
	p(X) = \frac{exp(X_1 + X_2)}{1 + exp(X_1 + X_2)},
\end{align*}
where $\epsilon \sim N\left( 0, 0.25^2 \right)$. The true CATE conditioning on $X_1$ can be derived as $\tau(x_1) = x_1 (1 + 2 x_1)^2 ( -1 + x_1)^2$. Since  the misspecification effect is a concern,  we use the misspecified parametric model respectively:
\begin{align*}
	\widetilde{m}_1(X;\gamma_1) = (1, X^{\top}) \gamma_1, \qquad \widetilde{p}(X;\beta) = \frac{exp \left( (1, X_1)\beta \right)}{1 + exp \left( (1, X_1)\beta \right)}.
\end{align*}
where $\gamma_1 \in \mathbb{R}^3$, $\beta \in \mathbb{R}^2$.

{\bf Model 2}. Another DGP is featured by a 4-dimensional unconfounded covariate for the purpose of a further investigation on higher dimension cases. Write $X = (X_1, X_2, X_3, X_4)^{\top}$ and
\begin{align*}
	X_1 = \rho_1, & \qquad X_2 = 1 + 2 X_1 + \rho_2, \\
	X_3 = 1 + 2 X_1 + \rho_3,  & \qquad X_4 = (-1 + X_1)^2 + \rho_4,
\end{align*}
where $\rho_1, \rho_2, \rho_3, \rho_4$ are independently identically $U(-0.5,0.5)$ distributed. The  potential outcomes and the propensity score function are defined as:
\begin{align*}
	Y(1) = X_1 X_2 X_3 X_4 + \epsilon, \qquad Y(0) = 0, \\
	p(X) = \frac{exp\left[ \frac{1}{2}(X_1 + X_2 + X_3 + X_4) \right]}{1 + exp\left[ \frac{1}{2}(X_1 + X_2 + X_3 + X_4) \right]},
\end{align*}
where $\epsilon \sim \mathcal{N}\left( 0, 0.25^2 \right)$. The true CATE conditioning on $X_1$ remains as $\tau(x_1) = x_1 (1 + 2 x_1)^2 ( -1 + x_1)^2$. Still we use the misspecified parametric model respectively:
\begin{align*}
	\widetilde{m}_1(X;\gamma_1) = (1, X^{\top}) \gamma_1, \qquad \widetilde{p}(X;\beta) = \frac{exp \left( (1, X_1)\beta \right)}{1 + exp \left( (1, X_1)\beta \right)}.
\end{align*}
where $\gamma_1 \in \mathbb{R}^5$, $\beta \in \mathbb{R}^2$.

\subsection{Kernel Functions and Bandwidths}

As  the selections of kernel functions and bandwidths (listed in the supplementary material) have great influence on the asymptotic property when the nuisance models are nonparametrically or semiparametrically estimated, we  first discuss this issue.

Let $h = a n^{- \eta}$ for $\eta > 0$. Together with condition (A2), how to determine the value $\eta$ goes to a linear programming problem.

For {\bf model 1} ($d=2$), we consider a kernel function of order 4 ($s_1 = 4$) as the kernel in the second step of N-W estimation, $K_1$. Write $h_1 = a_1 h^{- \eta_1}$. For the other bandwidths, take $h_2$ as an example. The results in Section~4 requires that $s^* \geq s_2 \geq d$, we then  choose $s_2 = 2$. Also let $h_2 = a_2 n^{- \eta_2}$. Then let $(\eta_1, \eta_2) = \left( \frac{1}{9}, \frac{1}{4} \right)$. The  other bandwidths can also be determined similarly as $h_j = a_j n^{- \frac{1}{4}}, (j = 2, 3, 5, 6)$, when $s_j = 2, (j = 1, 2, 3, 5, 6)$. Also, these convergence rates  of $h_i$ to meet condition (A16). To choose $a_j, (j = 2, 3, 5, 6)$, we, by the rule of thumb, choose $a_1 = 0.1$, $a_2= 0.7$, $a_3 = 1.5$, $a_5 = 0.5$ and $a_6 = 1$. For {\bf model 2} ($d=4$),  consider $s_1 = 6$ and $s_j = 4, (j = 2, 3, 5, 6)$ and $h_1 = a_1 n^{- \frac{1}{13}}$ and $h_j = a_j n^{- \frac{1}{8}}, (j = 2, 3, 5, 6)$. Further, let $a_1 = 0.1$, $a_2 = 2$, $a_3 = 2.5$, $a_5 = 2.8$ and $a_6 = 1$. In simulations, we chose many other values and found that the above values are recommendable as the values around them can make the estimators relatively stable.

 Consider the Gaussian kernel $K_1$ of order $s_1$ under condition (A1)(\romannumeral1). For other kernel functions, use Epanechnikov kernels of the corresponding orders under  conditions (A1)(\romannumeral2) and (\romannumeral3).

\subsection{Simulation Results}

As there are many estimators $\widehat{\tau}(x_1)$ with different estimated nuisance models, we then,  in Table \ref{simulation estimator}, list them and the corresponding notations for convenience.

\begin{table}[t!]
		\newcommand{\tabincell}[2]{\begin{tabular}{@{}#1@{}}#2\end{tabular}}
		\caption{Estimators involved in simulation}
		\label{simulation estimator}
		\resizebox{\linewidth}{!}{\begin{tabular}{|c|c|c|}
			\hline
			\textbf{ DRCATE} & $p(x)$ & $m_1(x)$ \\
			\hline
			\textbf{(O, O)} & oracle & oracle \\
			\hline
			\textbf{(cP, cP)} & \tabincell{c}{parametrically estimated \\ (correctly specified)} & \tabincell{c}{parametrically estimated \\ (correctly specified)} \\
			\hline
			\textbf{(N, N)} & nonparametrically estimated & nonparametrically estimated \\
			\hline
			\textbf{(S, S)} & semiparametrically estimated & semiparametrically estimated \\
			\hline
			\textbf{(mP, cP)} & \tabincell{c}{parametrically estimated \\ (misspecified)} & \tabincell{c}{parametrically estimated \\ (correctly specified)} \\
			\hline
			\textbf{(mP, N)} & \tabincell{c}{parametrically estimated \\ (misspecified)} & nonparametrically estimated \\
			\hline
			\textbf{(mP, S)} & \tabincell{c}{parametrically estimated \\ (misspecified)} & semiparametrically estimated \\
			\hline
			\textbf{(cP, mP)} & \tabincell{c}{parametrically estimated \\ (correctly specified)} & \tabincell{c}{parametrically estimated \\ (misspecified)} \\
			\hline
			\textbf{(N, mP)} & nonparametrically estimated & \tabincell{c}{parametrically estimated \\ (misspecified)} \\
			\hline
			\textbf{(S, mP)} & semiparametrically estimated & \tabincell{c}{parametrically estimated \\ (misspecified)} \\
			\hline
		\end{tabular}}
\end{table}

\noindent To guarantee the regularity conditions and the  estimation stability,  all estimated propensity scores are trimmed within $[0.005, 0.995]$ as many authors did.

\begin{table}[t!]
\caption{The simulation results under {\bf model 1} (part 1)}
\label{result 1}
\resizebox{\linewidth}{!}{\begin{tabular}{|c|r|rrrrr|rrrrr|}
\hline
 & \multicolumn{1}{c}{}   & \multicolumn{5}{c}{n=500} & \multicolumn{5}{c}{n=5000} \\
\cline{3-12}
DRCATE & \multicolumn{1}{c}{$x_1$} & \multicolumn{1}{c}{bias} & \multicolumn{1}{c}{sam-SD} & \multicolumn{1}{c}{MSE} & \multicolumn{1}{c}{$P_{0.05}$} & \multicolumn{1}{c}{$P_{0.95}$} & \multicolumn{1}{c}{bias} & \multicolumn{1}{c}{sam-SD} & \multicolumn{1}{c}{MSE} & \multicolumn{1}{c}{$P_{0.05}$} & \multicolumn{1}{c}{$P_{0.95}$} \\
\hline
\multirow{5}{*}{\begin{tabular}[c]{@{}c@{}} (O,O)\end{tabular}} & -0.4 & 0.0001 & 0.2776 & 0.0770 & 0.052 & 0.046 & 0.0004 & 0.2724 & 0.0742 & 0.044 & 0.052 \\
                         			& -0.2 & -0.0023 & 0.2378 & 0.0567 & 0.056 & 0.044 & -0.0005 & 0.2333 & 0.0544 & 0.049 & 0.050 \\
                         			& 0 & -0.0002 & 0.2088 & 0.0436 & 0.049 & 0.050 & 0.0003 & 0.2014 & 0.0405 & 0.047 & 0.048 \\
                        			& 0.2 & 0.0003 & 0.1997 & 0.0399 & 0.052 & 0.047 & 0.0002 & 0.1999 & 0.0400 & 0.050 & 0.054 \\
                                    & 0.4 & 0.0027 & 0.2003 & 0.0403 & 0.045 & 0.058 & 0.0004 & 0.2006 & 0.0403 & 0.048 & 0.054 \\
\hline
\multirow{5}{*}{\begin{tabular}[c]{@{}c@{}} (cP,cP)\end{tabular}} & -0.4 & 0.0000 & 0.2797 & 0.0782  & 0.053 & 0.048 & 0.0004 & 0.2725 & 0.0743 & 0.044 & 0.052 \\
                         			  & -0.2 & -0.0023 & 0.2378 & 0.0567 & 0.056 & 0.042 & -0.0005 & 0.2333 & 0.0544 & 0.051 & 0.048 \\
                         			  & 0 & -0.0002 & 0.2089 & 0.0436 & 0.048 & 0.050 & 0.0003 & 0.2014 & 0.0405 & 0.047 & 0.047 \\
                         			  & 0.2 & 0.0003 & 0.1994 & 0.0397 & 0.051 & 0.048 & 0.0002 & 0.2001 & 0.0400 & 0.051 & 0.054 \\
                        			  & 0.4 & 0.0027 & 0.2003 & 0.0403 & 0.044 & 0.058 & 0.0004 & 0.2007 & 0.0403 & 0.047 & 0.054 \\
\hline
\multirow{5}{*}{\begin{tabular}[c]{@{}c@{}} (N,N)\end{tabular}} & -0.4 & 0.0008 & 0.2716 & 0.0738 & 0.050 & 0.053 &  0.0001 & 0.2845 & 0.0809 & 0.050 & 0.049 \\
                         			 & -0.2 & 0.0015 & 0.2366 & 0.0560 & 0.042 & 0.058 & -0.0001 & 0.2344 & 0.0549 & 0.050 & 0.050 \\
                         			 & 0 & 0.0002 & 0.2046 & 0.0419 & 0.043 & 0.052 & -0.0005 & 0.1996 & 0.0399 & 0.057 & 0.041 \\
                         			 & 0.2 & 0.0010 & 0.2000 & 0.0400 & 0.044 & 0.051 & -0.0001 & 0.1941 & 0.0377 & 0.052 & 0.056 \\
                         			 & 0.4 & 0.0014 & 0.2081 & 0.0433 & 0.045 & 0.054 & 0.0009 & 0.2012 & 0.0406 & 0.045 & 0.056 \\
\hline
\multirow{5}{*}{\begin{tabular}[c]{@{}c@{}} (S,S)\end{tabular}} & -0.4 & -0.0022 & 0.2815 & 0.0794 & 0.051 & 0.044 & 0.0002 & 0.2862 & 0.0819 & 0.045 & 0.050 \\
                         			& -0.2 & 0.0004 & 0.2365 & 0.0559 & 0.046 & 0.052 & -0.0004 & 0.2302 & 0.0530 & 0.046 & 0.048 \\
                         			& 0 & 0.0005 & 0.2082 & 0.0433 & 0.053 & 0.052 & 0.0003 & 0.2059 & 0.0424 & 0.052 & 0.052 \\
                         			& 0.2 & -0.0015 & 0.1992 & 0.0397 & 0.061 & 0.041 & -0.0002 & 0.2011 & 0.0404 & 0.053 & 0.051 \\
                         			& 0.4 & 0.0002 & 0.2021 & 0.0408 & 0.050 & 0.046 & 0.0012 & 0.2048 & 0.0422 & 0.043 & 0.059 \\
\hline
\end{tabular}}
\end{table}

\begin{table}[t!]
\caption{The simulation results under {\bf model 1} (part 2)}
\label{result 2}
\resizebox{\linewidth}{!}{\begin{tabular}{|c|r|rrrrr|rrrrr|}
\hline
 & \multicolumn{1}{c}{}   & \multicolumn{5}{c}{n=500} & \multicolumn{5}{c}{n=5000} \\
\cline{3-12}
DRCATE & \multicolumn{1}{c}{$x_1$} & \multicolumn{1}{c}{bias} & \multicolumn{1}{c}{sam-SD} & \multicolumn{1}{c}{MSE} & \multicolumn{1}{c}{$P_{0.05}$} & \multicolumn{1}{c}{$P_{0.95}$} & \multicolumn{1}{c}{bias} & \multicolumn{1}{c}{sam-SD} & \multicolumn{1}{c}{MSE} & \multicolumn{1}{c}{$P_{0.05}$} & \multicolumn{1}{c}{$P_{0.95}$} \\
\hline
\multirow{5}{*}{\begin{tabular}[c]{@{}c@{}} (O,O)\end{tabular}} & -0.4 & 0.0001 & 0.2776 & 0.0770 & 0.052 & 0.046 & 0.0004 & 0.2724 & 0.0742 & 0.044 & 0.052 \\
                         			& -0.2 & -0.0023 & 0.2378 & 0.0567 & 0.056 & 0.044 & -0.0005 & 0.2333 & 0.0544 & 0.049 & 0.050 \\
                         			& 0 & -0.0002 & 0.2088 & 0.0436 & 0.049 & 0.050 & 0.0003 & 0.2014 & 0.0405 & 0.047 & 0.048 \\
                        			& 0.2 & 0.0003 & 0.1997 & 0.0399 & 0.052 & 0.047 & 0.0002 & 0.1999 & 0.0400 & 0.050 & 0.054 \\
                                    & 0.4 & 0.0027 & 0.2003 & 0.0403 & 0.045 & 0.058 & 0.0004 & 0.2006 & 0.0403 & 0.048 & 0.054 \\
\hline
\multirow{5}{*}{\begin{tabular}[c]{@{}c@{}} (mP,cP)\end{tabular}} & -0.4 & 0.0000 & 0.2599 & 0.0675 & 0.052 & 0.049 & 0.0004 & 0.2530 & 0.0640 & 0.044 & 0.052 \\
                         			   & -0.2 & -0.0022 & 0.2363 & 0.0559 & 0.056 & 0.041 & -0.0005 & 0.2323 & 0.0540 & 0.050 & 0.050 \\
                         			   & 0 & -0.0002 & 0.2203 & 0.0485 & 0.049 & 0.048 & 0.0003 & 0.2116 & 0.0448 & 0.047 & 0.052 \\
                         			   & 0.2 & 0.0003 & 0.2041 & 0.0417 & 0.051 & 0.046 & 0.0002 & 0.2048 & 0.0419 & 0.050 & 0.053 \\
                         			   & 0.4 & 0.0027 & 0.1953 & 0.0383 & 0.044 & 0.058 & 0.0004 & 0.1955 & 0.0382 & 0.046 & 0.054 \\
\hline
\multirow{5}{*}{\begin{tabular}[c]{@{}c@{}} (mP,N)\end{tabular}} & -0.4 & -0.0046 & 0.2666 & 0.0716 & 0.064 & 0.040 & -0.0011 & 0.2629 & 0.0693 & 0.054 & 0.044 \\
                         			  & -0.2 & -0.0035 & 0.2373 & 0.0566 & 0.059 & 0.044 & -0.0029 & 0.2383 & 0.0584 & 0.074 & 0.037 \\
                         			  & 0 & -0.0068 & 0.2152 & 0.0474 & 0.072 & 0.032 & -0.0027 & 0.2107 & 0.0458 & 0.072 & 0.034 \\
                         			  & 0.2 & -0.0011 & 0.2041 & 0.0417 & 0.052 & 0.047 & -0.0004 & 0.1952 & 0.0381 & 0.050 & 0.045 \\
                         			  & 0.4 & -0.0008 & 0.2003 & 0.0401 & 0.049 & 0.049 & 0.0007 & 0.2002 & 0.0402 & 0.043 & 0.056 \\
\hline
\multirow{5}{*}{\begin{tabular}[c]{@{}c@{}} (mP,S)\end{tabular}} & -0.4 & -0.0143 & 0.2701 & 0.0781 & 0.082 & 0.029 & -0.0115 & 0.2722 & 0.0996 & 0.146 & 0.010 \\
                         			  & -0.2 & -0.0094 & 0.2453 & 0.0624 & 0.070 & 0.032 & -0.0073 & 0.2302 & 0.0634 & 0.114 & 0.016 \\
                         			  & 0 & -0.0046 & 0.2116 & 0.0453 & 0.064 & 0.043 & -0.0038 & 0.2099 & 0.0469 & 0.083 & 0.032 \\
                         			  & 0.2 & -0.0019 & 0.2041 & 0.0417 & 0.050 & 0.046 & -0.0006 & 0.1970 & 0.0388 & 0.054 & 0.047 \\
                         			  & 0.4 & 0.0022 & 0.2002 & 0.0402 & 0.046 & 0.058 & 0.0017 & 0.1968 & 0.0393 & 0.037 & 0.062 \\
\hline
\end{tabular}}
\end{table}

\begin{table}[t!]
\caption{The simulation results under {\bf model 1} (part 3)}
\label{result 3}
\resizebox{\linewidth}{!}{\begin{tabular}{|c|r|rrrrr|rrrrr|}
\hline
 & \multicolumn{1}{c}{}   & \multicolumn{5}{c}{n=500} & \multicolumn{5}{c}{n=5000} \\
\cline{3-12}
DRCATE & \multicolumn{1}{c}{$x_1$} & \multicolumn{1}{c}{bias} & \multicolumn{1}{c}{sam-SD} & \multicolumn{1}{c}{MSE} & \multicolumn{1}{c}{$P_{0.05}$} & \multicolumn{1}{c}{$P_{0.95}$} & \multicolumn{1}{c}{bias} & \multicolumn{1}{c}{sam-SD} & \multicolumn{1}{c}{MSE} & \multicolumn{1}{c}{$P_{0.05}$} & \multicolumn{1}{c}{$P_{0.95}$} \\
\hline
\multirow{5}{*}{\begin{tabular}[c]{@{}c@{}} (O,O)\end{tabular}} & -0.4 & 0.0001 & 0.2776 & 0.0770 & 0.052 & 0.046 & 0.0004 & 0.2724 & 0.0742 & 0.044 & 0.052 \\
                         			& -0.2 & -0.0023 & 0.2378 & 0.0567 & 0.056 & 0.044 & -0.0005 & 0.2333 & 0.0544 & 0.049 & 0.050 \\
                         			& 0 & -0.0002 & 0.2088 & 0.0436 & 0.049 & 0.050 & 0.0003 & 0.2014 & 0.0405 & 0.047 & 0.048 \\
                        			& 0.2 & 0.0003 & 0.1997 & 0.0399 & 0.052 & 0.047 & 0.0002 & 0.1999 & 0.0400 & 0.050 & 0.054 \\
                                    & 0.4 & 0.0027 & 0.2003 & 0.0403 & 0.045 & 0.058 & 0.0004 & 0.2006 & 0.0403 & 0.048 & 0.054 \\
\hline
\multirow{5}{*}{\begin{tabular}[c]{@{}c@{}} (cP,mP)\end{tabular}} & -0.4 & -0.0012 & 0.3230 & 0.1044 & 0.051 & 0.048 & 0.0001 & 0.3201 & 0.1024 & 0.050 & 0.049 \\
                         			   & -0.2 & -0.0021 & 0.2400 & 0.0577 & 0.052 & 0.042 & -0.0005 & 0.2362 & 0.0558 & 0.054 & 0.044 \\
                         			   & 0 & 0.0004 & 0.2147 & 0.0461 & 0.052 & 0.049 & 0.0003 & 0.2050 & 0.0420 & 0.049 & 0.049 \\
                         			   & 0.2 & 0.0004 & 0.2012 & 0.0405 & 0.054 & 0.046 & 0.0001 & 0.2016 & 0.0406 & 0.048 & 0.049 \\
                         			   & 0.4 & 0.0028 & 0.2059 & 0.0426 & 0.043 & 0.061 & 0.0004 & 0.2039 & 0.0416 & 0.045 & 0.053 \\
\hline
\multirow{5}{*}{\begin{tabular}[c]{@{}c@{}} (N,mP)\end{tabular}} & -0.4 & -0.0105 & 0.2840 & 0.0834 & 0.075 & 0.040 & -0.0013 & 0.2970 & 0.0885 & 0.060 & 0.045 \\
                         			  & -0.2 & 0.0014 & 0.2353 & 0.0554 & 0.047 & 0.050 & 0.0007 & 0.2288 & 0.0525 & 0.040 & 0.053 \\
                         			  & 0 & 0.0013 & 0.2104 & 0.0443 & 0.048 & 0.054 & 0.0002 & 0.2065 & 0.0426 & 0.047 & 0.044 \\
                         			  & 0.2 & -0.0014 & 0.1995 & 0.0398 & 0.056 & 0.048 & -0.0004 & 0.2022 & 0.0409 & 0.052 & 0.044 \\
                         			  & 0.4 & 0.0008 & 0.2034 & 0.0414 & 0.046 & 0.046 & 0.0000 & 0.2077 & 0.0431 & 0.048 & 0.050 \\
\hline
\multirow{5}{*}{\begin{tabular}[c]{@{}c@{}} (S,mP)\end{tabular}} & -0.4 & -0.0051 & 0.2964 & 0.0884 & 0.055 & 0.046 & -0.0005 & 0.3089 & 0.0955 & 0.050 & 0.045 \\
                         			  & -0.2 & -0.0002 & 0.2421 & 0.0586 & 0.049 & 0.050 & 0.0001 & 0.2394 & 0.0573 & 0.048 & 0.051 \\
                         			  & 0 & 0.0005 & 0.2076 & 0.0431 & 0.050 & 0.050 & -0.0001 & 0.2051 & 0.0421 & 0.048 & 0.049 \\
                         			  & 0.2 & -0.0008 & 0.2082 & 0.0433 & 0.049 & 0.049 & -0.0001 & 0.1966 & 0.0386 & 0.054 & 0.048 \\
                         			  & 0.4 & 0.0005 & 0.2104 & 0.0443 & 0.044 & 0.052 & 0.0006 & 0.2085 & 0.0435 & 0.048 & 0.054  \\
\hline
\end{tabular}}
\end{table}

\begin{figure}[h]
	\centering
	\subfigure[$n=500$]{
		\includegraphics[width=0.33\textwidth]{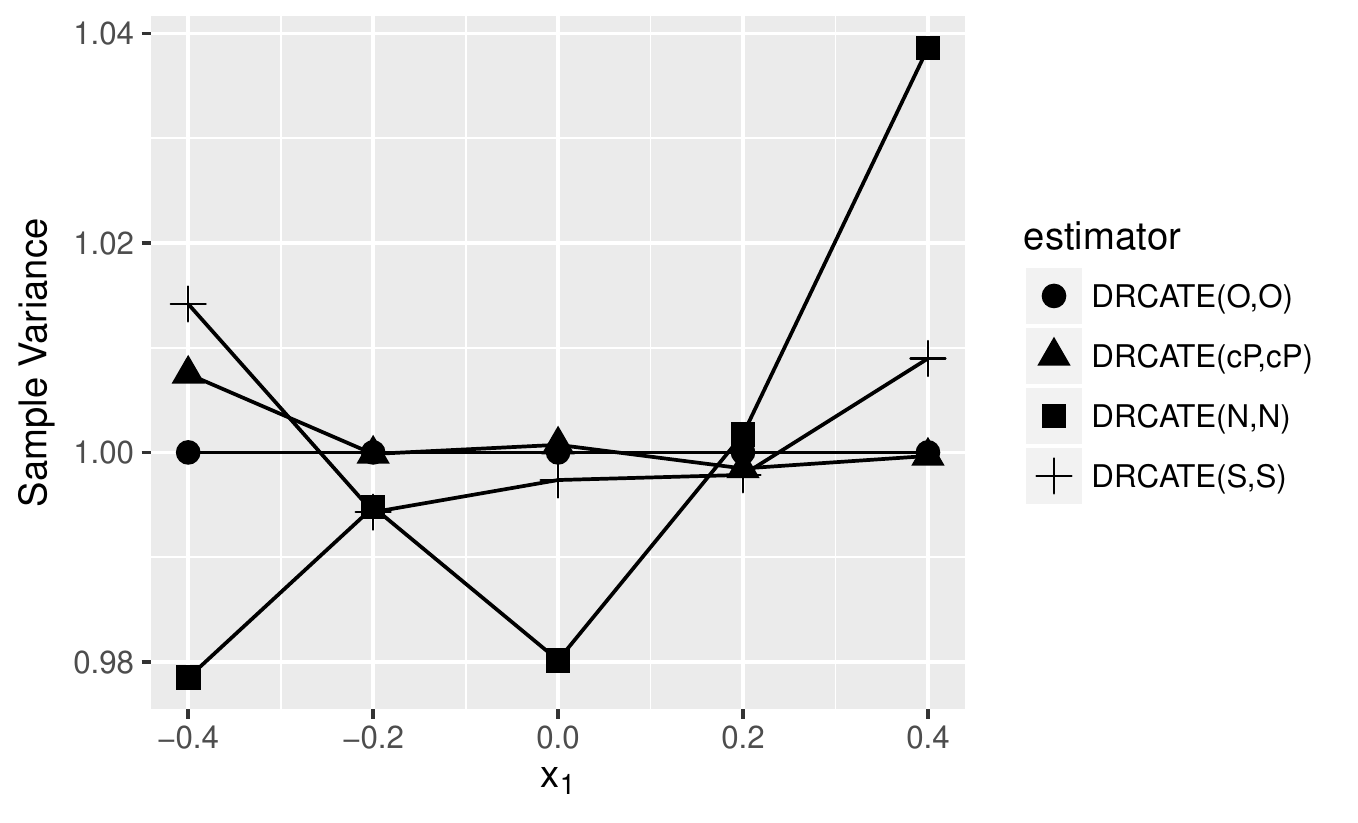}
		\includegraphics[width=0.33\textwidth]{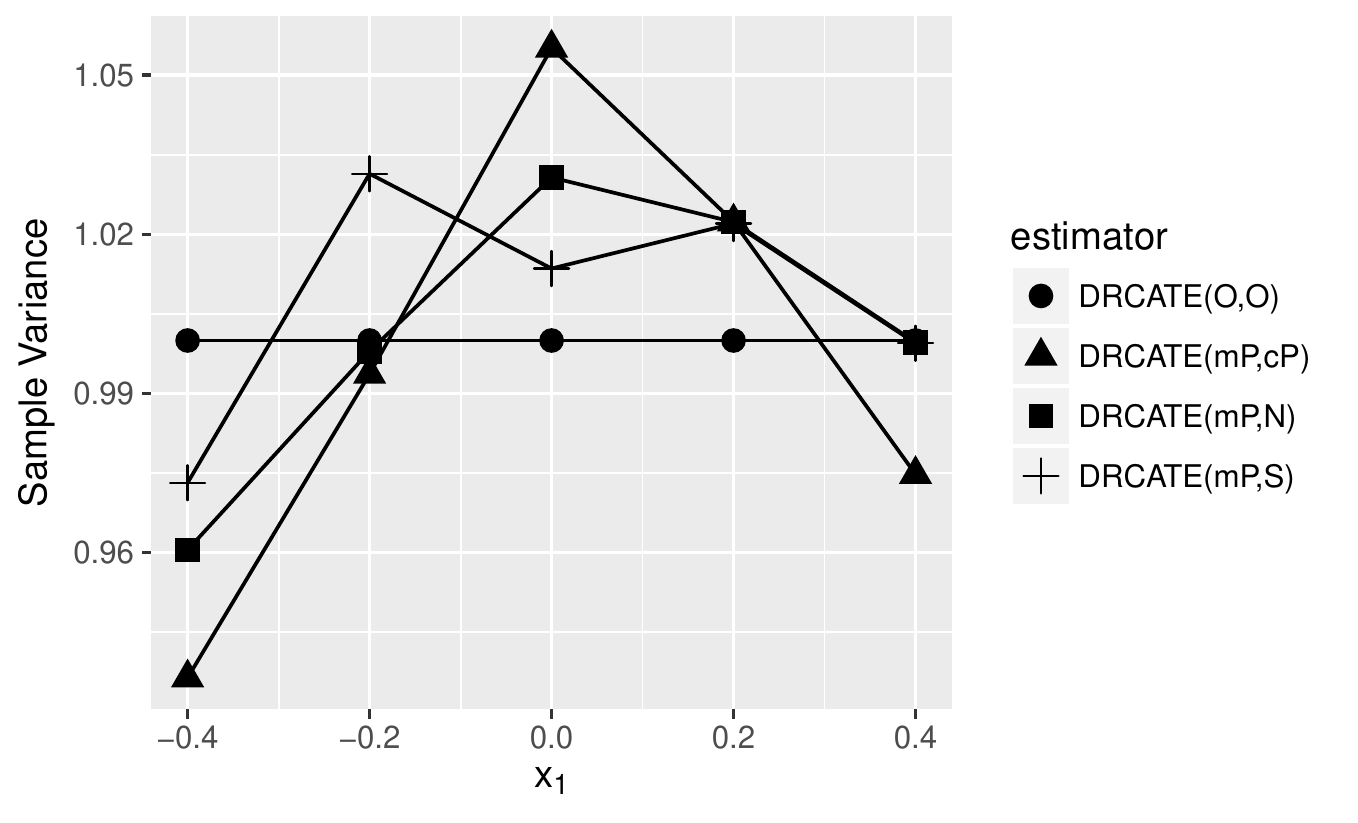}
		\includegraphics[width=0.33\textwidth]{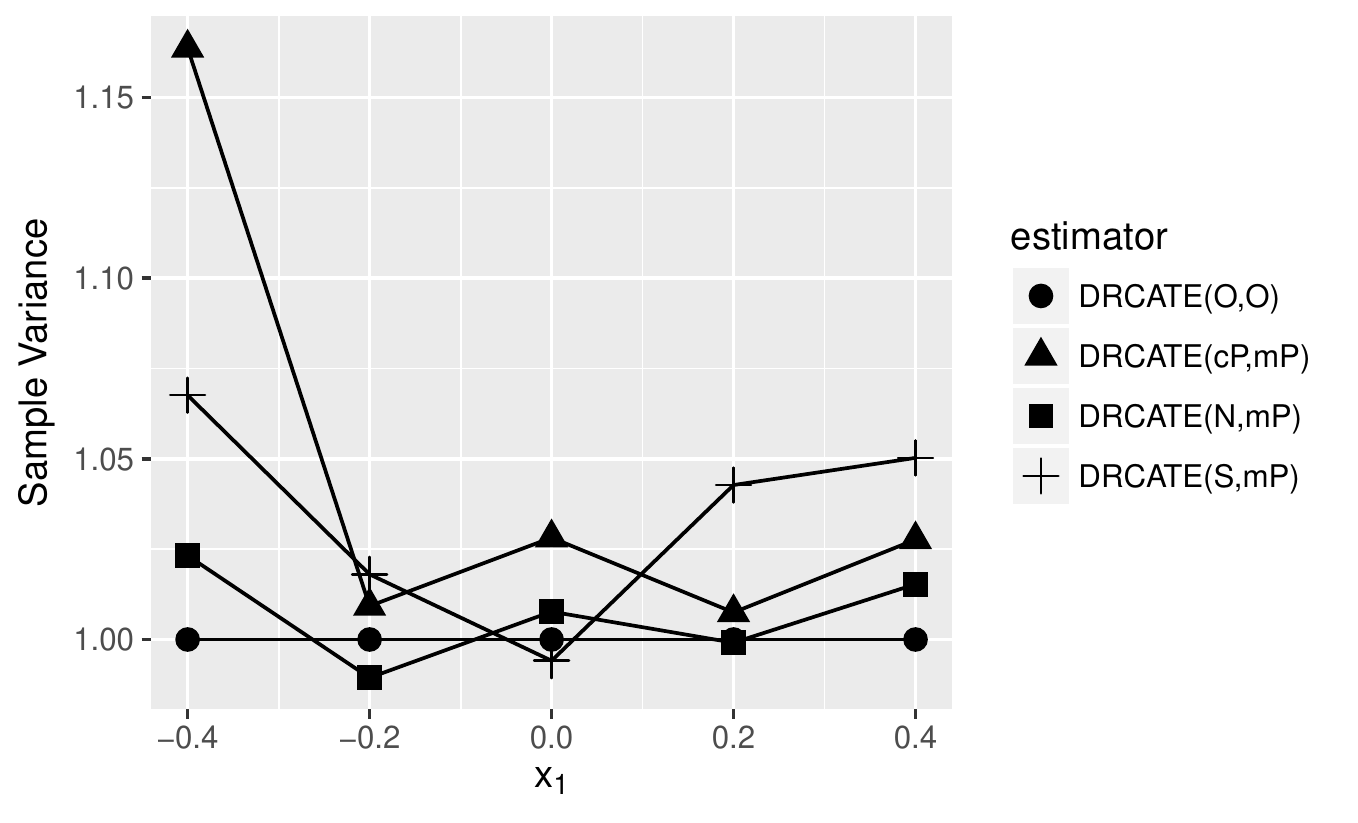}
	}
	\\
	\subfigure[$n=5000$]{
		\includegraphics[width=0.33\textwidth]{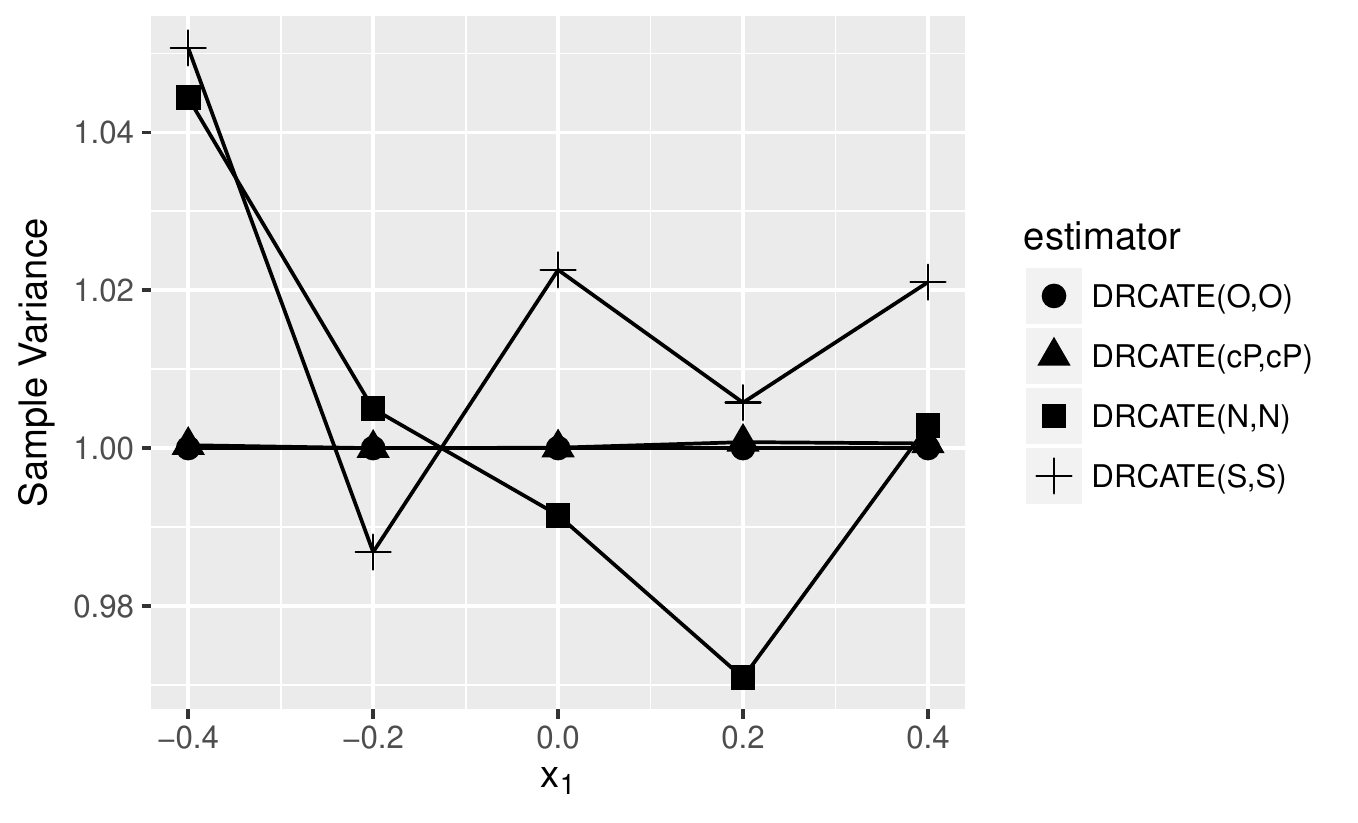}
		\includegraphics[width=0.33\textwidth]{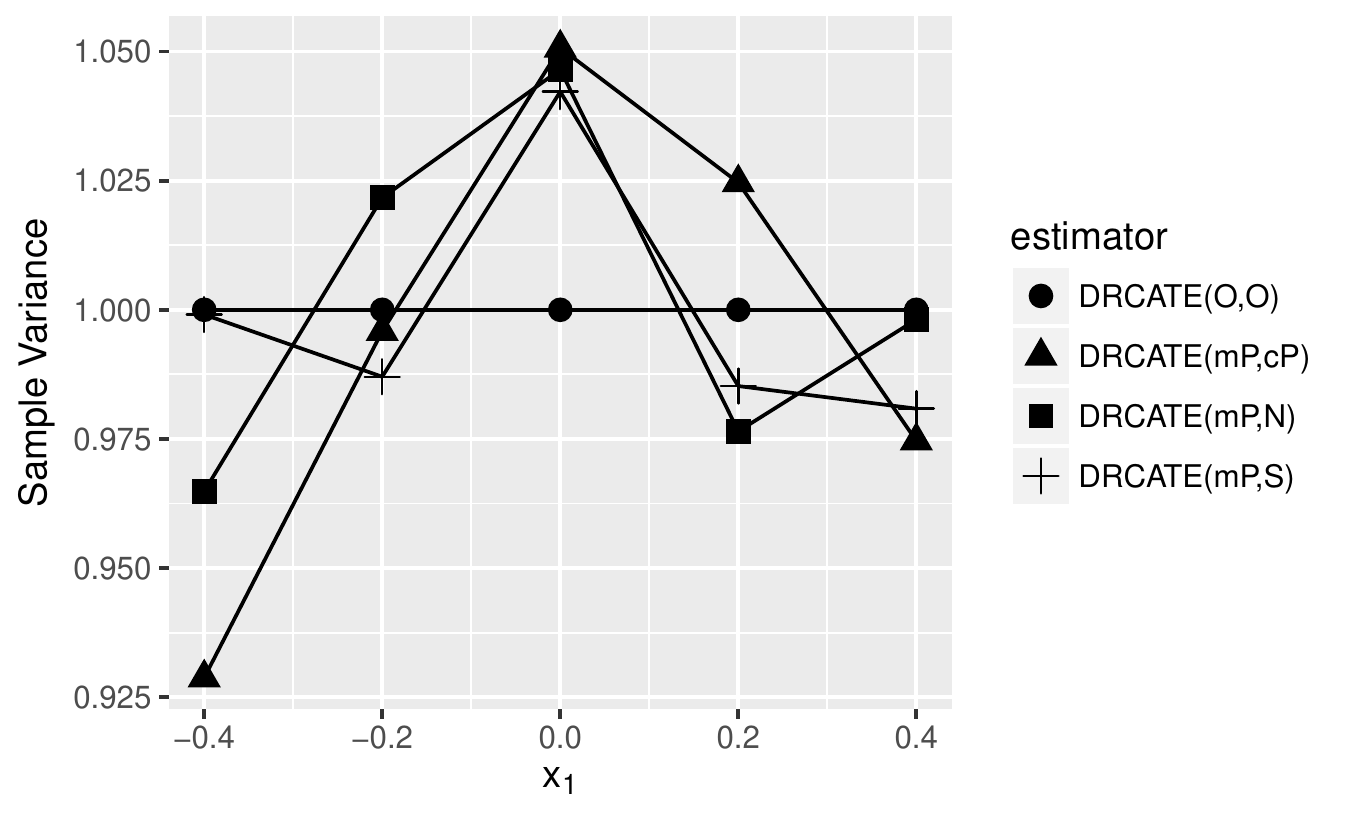}
		\includegraphics[width=0.33\textwidth]{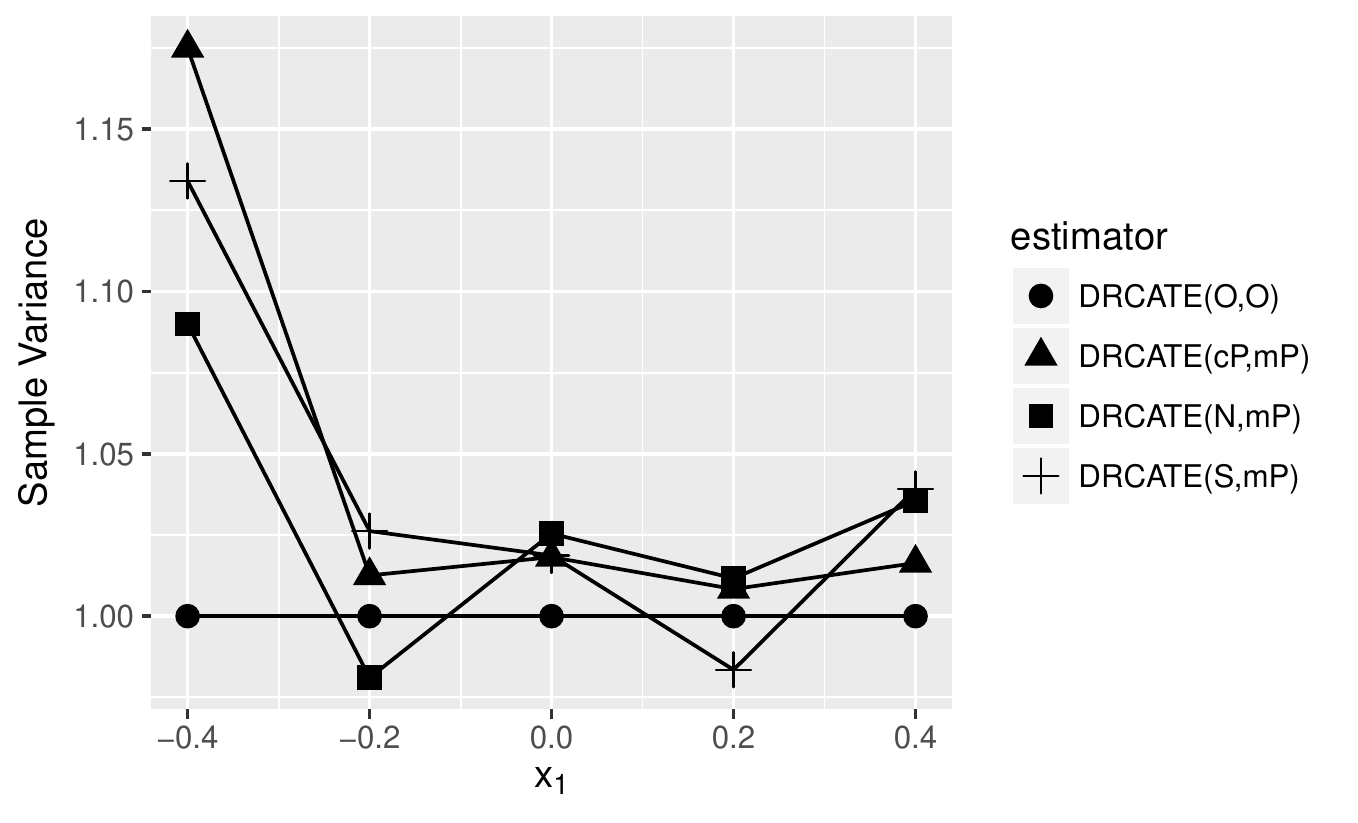}
	}
	\caption{Relative variance against DRCATE(O,O) in {\bf model 1}}
	\label{relative variance plot 1}
\end{figure}

\begin{figure}[h]
	\centering
	\subfigure[$n=500$]{
		\includegraphics[width=0.33\textwidth]{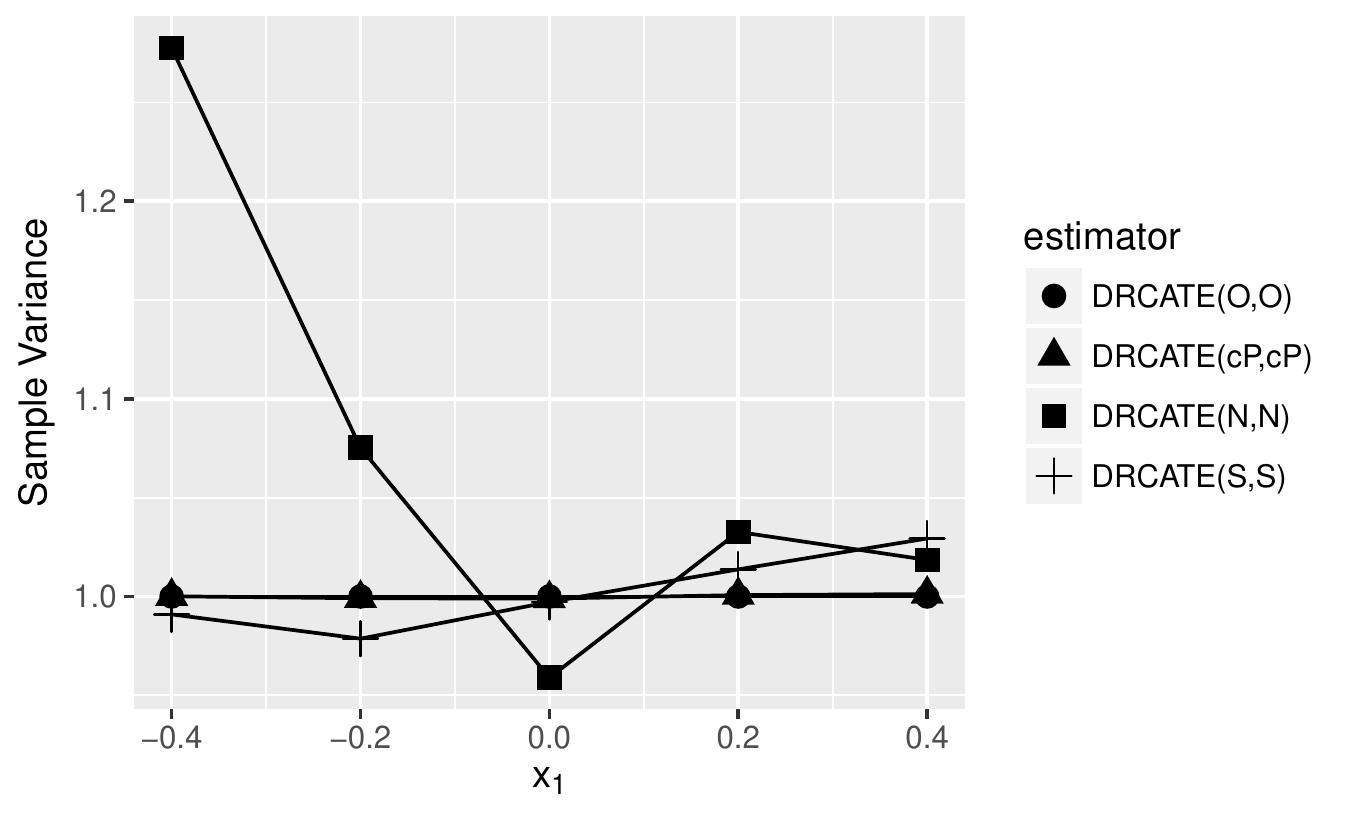}
		\includegraphics[width=0.33\textwidth]{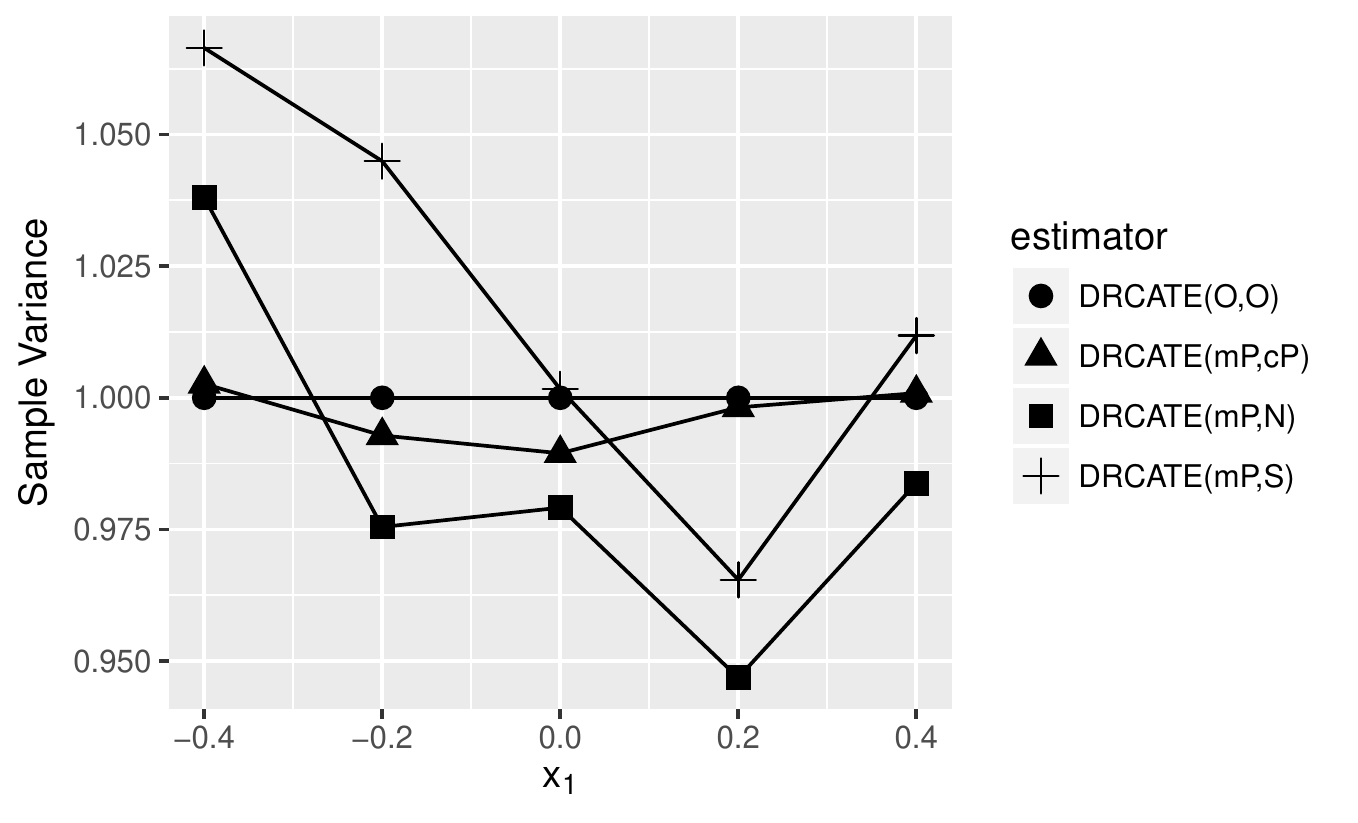}
		\includegraphics[width=0.33\textwidth]{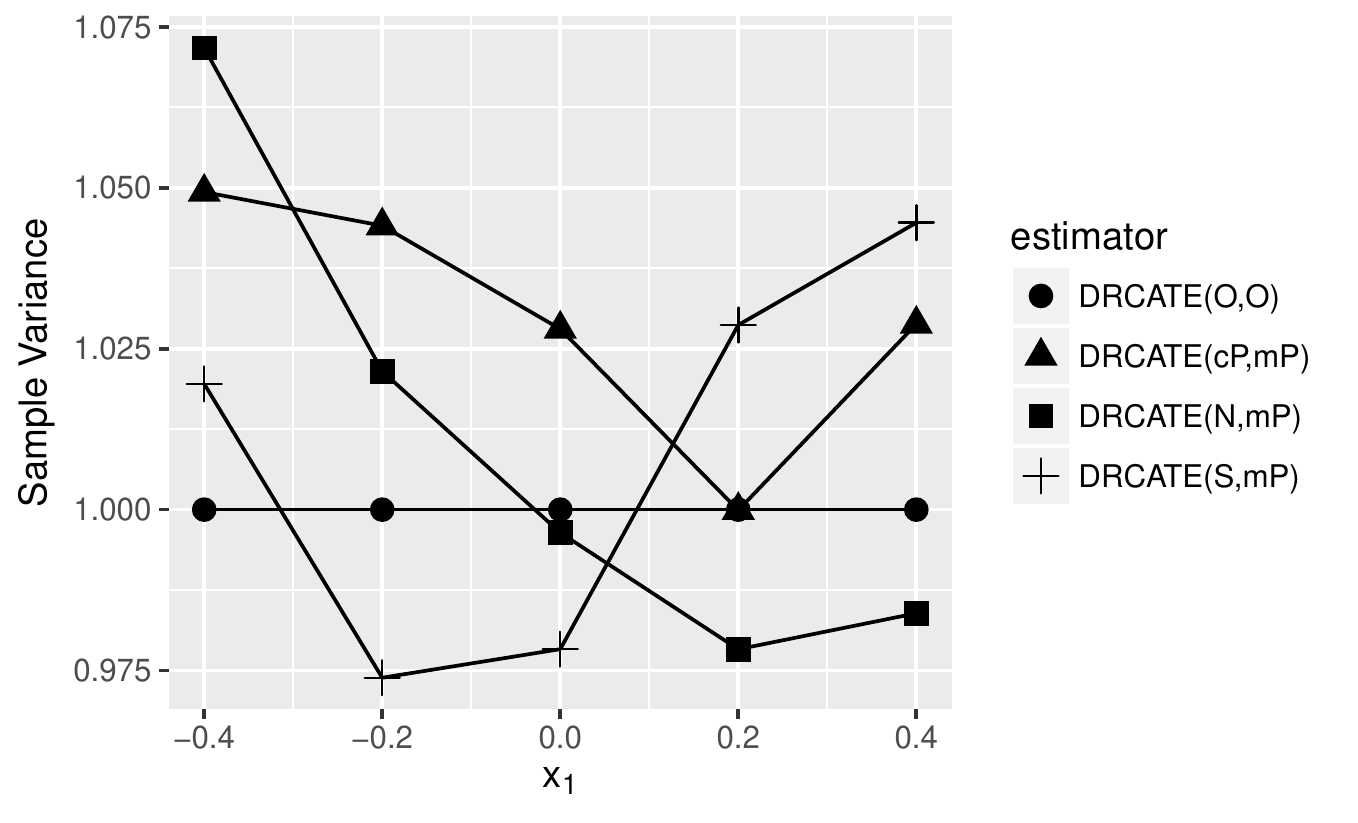}
	}
	\\
	\subfigure[$n=5000$]{
		\includegraphics[width=0.33\textwidth]{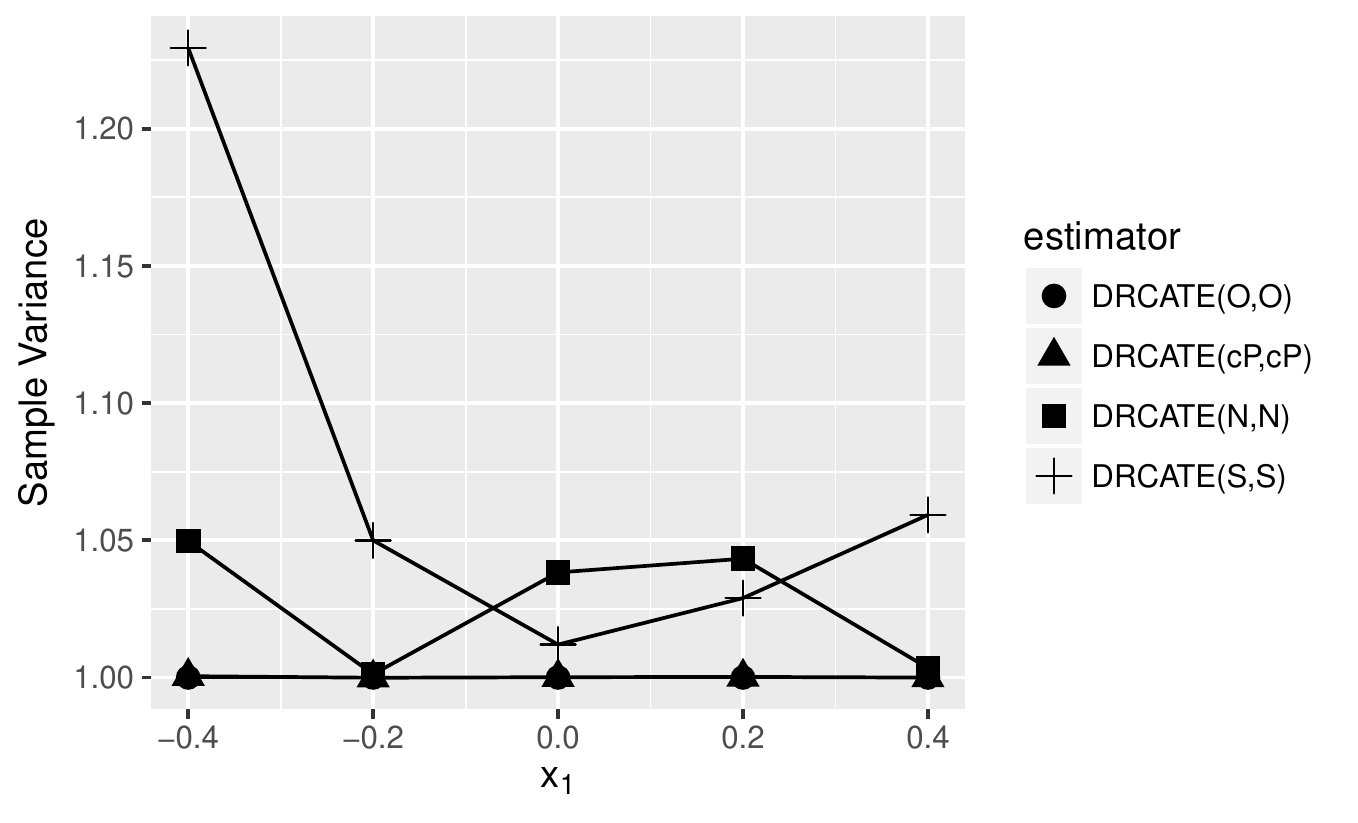}
		\includegraphics[width=0.33\textwidth]{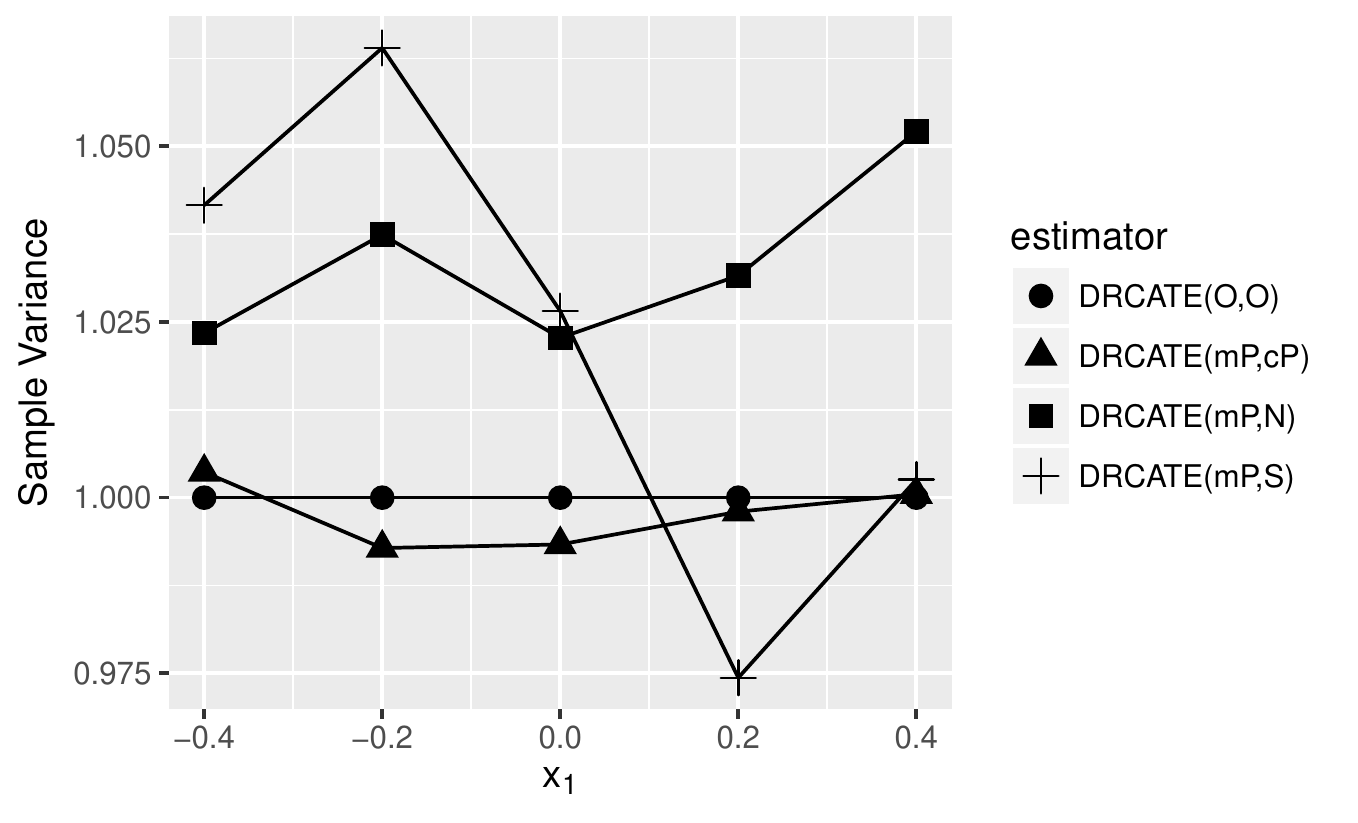}
		\includegraphics[width=0.33\textwidth]{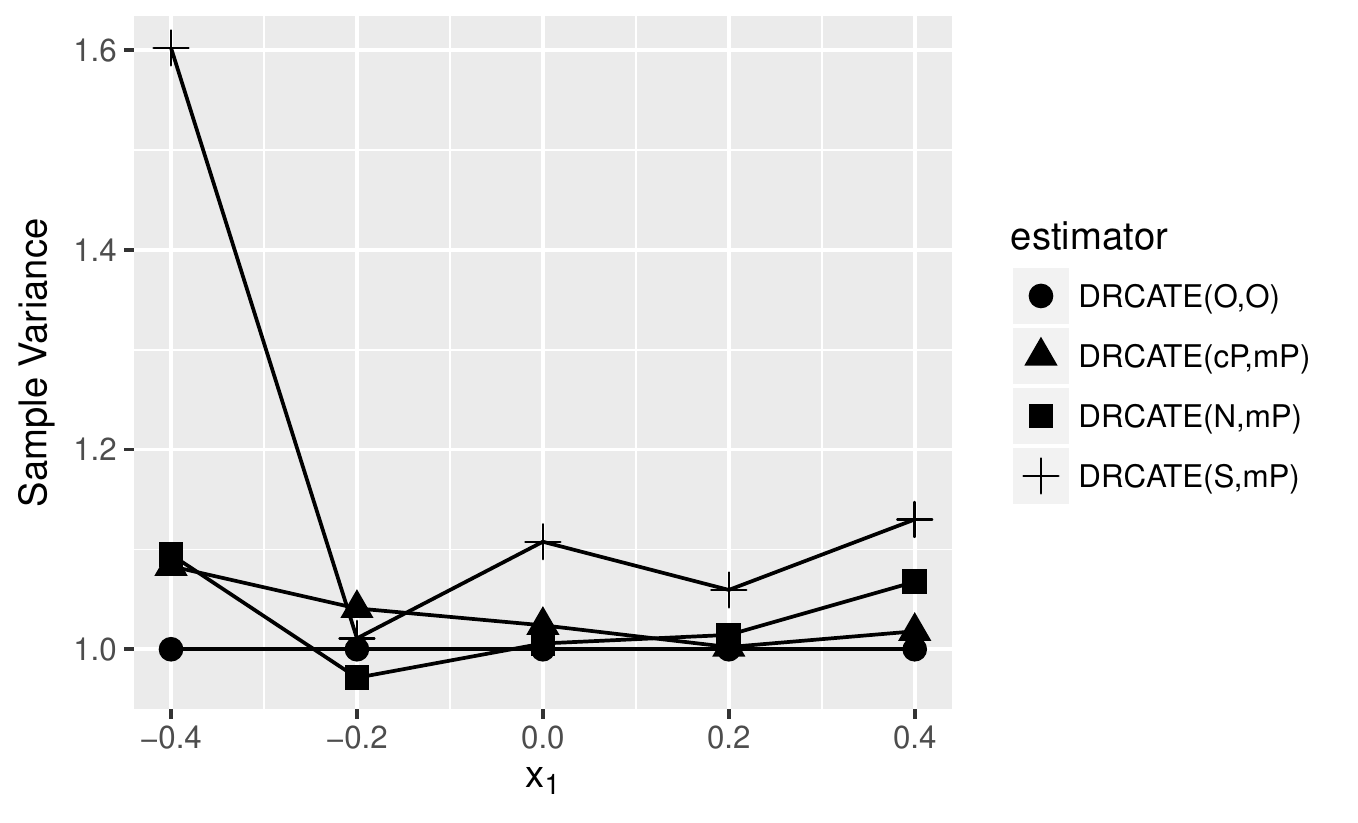}
	}
	\caption{Relative variance against DRCATE(O,O) in {\bf model 2}}
	\label{relative variance plot 2}
\end{figure}

In the simulations, we estimate $\tau(x_1)$ for $x_1 \in \{ -0.4, -0.2, 0, 0.2, 0.4 \}$. The sample sizes  are $n = 500$ and $n = 5,000$ respectively to see their asymptotic behaviours. The experiments are repeated $2,500$. Denote $T(x_1) = \sqrt{n h_1}\left( \widehat{\tau}(x_1) - \tau(x_1) \right)$. we evaluate the estimators based on following criteria: bias of $\widehat{\tau}(x_1)$; sample standard deviation (sam-SD) of $T(x_1)$; mean square error (MSE) of $T(x_1)$. We also report the proportions ($P_{0.05}$, $P_{0.95}$) of the standardized $T(x_1)$ below  the $5\%$ quantile and  above the $95\%$ quantile of $\mathcal{N} (0,1)$ to verify the asymptotic Normality. We display the efficiency comparisons among different estimators under models~1 and 2 in Figures \ref{relative variance plot 1} and \ref{relative variance plot 2} and  the detailed results under {\bf model 1} are displayed in Tables~\ref{result 1}, \ref{result 2} and \ref{result 3}. To save space, the other  simulation results about {\bf model 2} are reported in the supplementary material.

Here we present some observations from the simulation results.

First, with the  sample size growth, the bias and the standard deviation of $\widehat{\tau}(x_1)$ reasonably tend to be smaller due to the estimation consistency. The reported proportions $P_{0.05}$ and $P_{0.95}$ can be controlled around 0.05, which implies that the normal approximation of the proposed estimator is valid.

Second, from Figures \ref{relative variance plot 1} and \ref{relative variance plot 2}, the efficiency comparisons among the estimators (O,O), (cP,cP), (N,N) and (S,S) show that the distributions are close to each other. When only the propensity score function is misspecified, variance inflation and shrinkage are both possible. With misspecified outcome regression function, only  variance inflation is possible.

Third, the bias and  standard deviation of $\widehat{\tau}(x_1)$  increase when the covariate dimension grows, see  the comparisons in Figures \ref{relative variance plot 1} and \ref{relative variance plot 2}. Possible explanation to this phenomenon would be that the standard deviations of nuisance models' estimations increase with higher dimension of  covariate.

\section{Conclussion}

In this paper, we investigate the asymptotic behaviours of nine doubly robust estimators (DR), under different combinations of model structures, to provide a relatively complete picture of this methodology.

 When all models are correctly specified, the asymptotic equivalence among all defined estimators does not surprisingly hold. When models are mispecified, we consider local and global misspecifications and some interesting phenomena have been discovered such as asymptotic variance shrinking in some cases due to misspecification. Further, we would recommend semiparametric estimation under dimension reduction structure. This is because  nonparametric estimation severely suffers from the curse of dimensionality whereas parametric estimation may not be sufficiently robust against model structure.
 
 
\section*{Acknowledgements}

The research described herein was supported by a NNSF grant of China and  a grant from the University Grants Council of Hong Kong, Hong Kong, China.
\par

\newpage

\section{Supplementary Material}

The supplementary material contains the detailed proofs of the theorems and propositions, and the additional simulation results.

\subsection{Technical Conditions}

Here we present some conditions to derive the theoretical results. Together with (C1) and (C2) in the main context, the following conditions in the  (C) group are regularity conditions to guarantee the asymptotic properties regardless of the different ways to estimate nuisance models.
\begin{description}
	\item (C3) Density functions involved in this article satisfy the following conditions:
	\begin{description}
		\item (\romannumeral1) For any $x \in \mathcal{X}$, the density function of $X$, $\theta(x)$ is bounded away from 0.
		\item (\romannumeral2) For any $x_1$,  the density function of $X_1$, $f(x_1)$ is bounded away from zero and $s_1$ times continuously differentiable.
		\item (\romannumeral3) Denote the density functions of $A^{\top} X$, $B_1^{\top} X$ and $B_0^{\top} X$ as $\theta_A(\cdot)$, $\theta_{B1}(\cdot)$ and $\theta_{B0}(\cdot)$. For any $x \in \mathcal{X}$, all these density functions are bounded away from 0.
	\end{description}
	 \item (C4) Denote $\mathcal{C}$ as the parameter space of $\beta$. For any $x \in \mathcal{X}$ and $\beta \in \mathcal{C}$, $\widetilde{p}(x; \beta)$ is bounded away from 0 and 1.	
	\item (C5) $\sup_{x_1} E [ {Y(j)}^{2}| X_1 = x_1 ] < \infty$ for $j=0,1$.
	\item (C6) $E |\Psi_1(X, Y, D) - \tau(x_1)|^{2+\kappa_1} \le \infty$, $E |\Psi_2(X, Y, D) - \tau(x_1)|^{2+\kappa_2} \le \infty$, $E |\Psi_3(X, Y, D) - \tau(x_1)|^{2+\kappa_3} \le \infty$, $E |\Psi_4(X, Y, D) - \widetilde{\tau}(x_1)|^{2+\kappa_4} \le \infty$, $\int |K(u)|^{2+\delta} du \le \infty$ for some constants $\kappa_1, \kappa_2, \kappa_3, \kappa_4, \delta \ge 0$.
\end{description}

(C1) and (C2) in the main context are the basic conditions under  Rubin's potential outcome framework, as stated in Section~2. It is obvious that (C4) is an analogue of (C2)($\romannumeral2$). Bounded propensity scores or specified propensity score models, density functions and corresponding conditional moments are required in these conditions, which are common restrictions in the literature, and play important roles in deriving the asymptotic linear expression of the proposed estimators. (C6) ensures the applicability of Lyapunov's Central Limit Theorem here.

Assume some conditions on kernel functions and bandwidths in nonparametric estimation:
\begin{description}
	\item (A1) The kernel functions satisfy the following conditions:
	\begin{description}
		\item (\romannumeral1) $K_1(u)$ is a kernel function of order $s_1$, which is symmetric around zero and $s*$ times continuously differentiable.
		\item (\romannumeral2) $K_2(u)$, $K_3(u)$ and $K_4(u)$ are kernels of order $s_2 \geq p$, $s_3 \geq p$ and $s_4 \geq p$,  which are symmetric around zero and equal to zero outside $\prod_{i=1}^{p}[-1,1]$, with continuous $(s_2+1)$, $(s_3+1)$ and $(s_4+1)$ order derivatives respectively.
		\item (\romannumeral3) $K_5(u)$, $K_6(u)$ and $K_7(u)$ are kernels of order $s_5 \geq p(2)$, $s_6 \geq p(1)$ and $s_7 \geq p(0)$, which are symmetric around zero and equal to zero outside $\prod_{i=1}^{p(2)}[-1,1]$, $\prod_{i=1}^{p(1)}[-1,1]$ and $\prod_{i=1}^{p(0)}[-1,1]$, with continuous $(s_5+1)$, $(s_6+1)$ and $(s_7+1)$ order derivatives respectively.
	\end{description}
	\item (A2) As different scenarios require different bandwidths, we put them together in the following. As $n \to \infty$:
	\begin{description}
		\item (\romannumeral1) $h_1 \to 0$, $n h_1^k \to \infty$, $n h_1^{2s_1 + k} \to 0$.
		\item (\romannumeral2) $h_2 \to 0$, $(\ln n)/(nh_2^{p+s_2}) \to 0$.
		\item (\romannumeral3) $h_3 \to 0$, $(\ln n)/(nh_3^{p+s_3}) \to 0$.
		\item (\romannumeral4) $h_4 \to 0$, $(\ln n)/(nh_4^{p+s_4}) \to 0$.
		\item (\romannumeral5) $h_5 \to 0$, $(\ln n)/(nh_5^{p+s_5}) \to 0$.
		\item (\romannumeral6) $h_6 \to 0$, $(\ln n)/(nh_6^{p+s_6}) \to 0$.
		\item (\romannumeral7) $h_7 \to 0$, $(\ln n)/(nh_7^{p+s_7}) \to 0$.
		\item (\romannumeral8) $h_2^{2s_2}h_1^{ - 2 s_2 - k} \to 0$, $nh_1^k h_2^{2s_2} \to 0$.
		\item (\romannumeral9) $h_3^{2s_3}h_1^{ - 2 s_3 - k} \to 0$, $nh_1^k h_3^{2s_3} \to 0$.
		\item (\romannumeral10) $h_4^{2s_4}h_1^{ - 2 s_4 - k} \to 0$, $nh_1^k h_4^{2s_4} \to 0$.
		\item (\romannumeral11) $h_5^{2s_5}h_1^{ - 2 s_5 - k} \to 0$, $nh_1^k h_5^{2s_5} \to 0$.
		\item (\romannumeral12) $h_6^{2s_6}h_1^{ - 2 s_6 - k} \to 0$, $nh_1^k h_6^{2s_6} \to 0$.
		\item (\romannumeral13) $h_7^{2s_7}h_1^{ - 2 s_7 - k} \to 0$, $nh_1^k h_7^{2s_7} \to 0$.
		\item (\romannumeral14) $n h_1^k h_2^{s_2} h_3^{s_3} \to 0$, $n h_1^k h_2^{s_2} h_4^{s_4} \to 0$, $n h_1^k h_2^{s_2} h_6^{s_6} \to 0$, $n h_1^k h_2^{s_2} h_7^{s_7} \to 0$, $n h_1^k h_5^{s_5} h_3^{s_3} \to 0$, $n h_1^k h_5^{s_5} h_4^{s_4} \to 0$, $n h_1^k h_5^{s_5} h_6^{s_6} \to 0$, $n h_1^k h_5^{s_5} h_7^{s_7} \to 0$.
	\end{description}
\end{description}

Remind that $K_j(u), j=2,3,\dots,7$, and $h_j, j=2,3,\dots,7$ are corresponding kernels and bandwidths in nonparametric and semiparametric estimators of nuisance models. When only parametric methods are applied to estimate nuisance models, no  conditions above, but (A1)(\romannumeral1) and (A2)(\romannumeral1) are required.

The conditions in (A1)(\romannumeral2) and (A1)(\romannumeral3) are required when at least one misspecified model is involved. Epanechnikov kernel of corresponding order can be a candidate of $K_j(u), j=2,\dots,7$. \cite{abrevaya2015estimating} stated that this restriction on the bounded support can be relaxed to exponential tails.

(A2)(\romannumeral2)-(\romannumeral14) place restrictions on the convergence rates of different bandwidths to ensure remainders of the linear expression negligible. (A2)(\romannumeral14), involving more than 2 bandwidths, can be regarded as an interaction term, which makes it handleable to determine those convergence rates. Here we provide a naive idea to accomplish this task based on linear programming. Assume the corresponding bandwidths converge to $0$ in such a manner, $h_j = a_j n^{- \eta_j}, j = 1,\dots,7$, where $\eta_j > 0$. With predetermined $s_j, j=1,\dots,7$ and (A2), the problem goes to a linear programming task to find out the feasible region of $\eta_j$. For a more detailed example, reader can refer to Section 5.

Lastly, we give a condition to ensure  the desired convergence rates of the estimators under semiparametric dimension reduction structure, which will be a favour when pursuing the asymptotic properties of $\widehat{\tau}(x_1)$:
\begin{description}
	\item (B1) $\widehat{A} - A = O_p \left( n^{-1/2} \right)$, $\widehat{B}_1 - B_1 = O_p \left( n^{-1/2} \right)$, $\widehat{B}_0 - B_0 = O_p \left( n^{-1/2} \right)$
\end{description}
These can be achieved by standard estimations in the literature, see the relevant references such as Li (1991), ????

In summary, these conditions are rather standard.

\subsection{Proof of Theorem 1} \label{proof: all correct}
\setcounter{equation}{0}

Recall that
\begin{align*}
	\widehat{\tau}(x_1) = & \frac{\sum_{i=1}^n \left[ \frac{D_i}{\widehat{p}_i} \left( Y_i - \widehat{m}_{1i} \right) - \frac{1 - D_i}{1 - \widehat{p}_i} \left( Y_i - \widehat{m}_{0i} \right) + \widehat{m}_{1i} - \widehat{m}_{0i} \right] K_1 \left( \frac{X_{1i} - x_1}{h_1} \right)}{\sum_{t=1}^n K_1 \left( \frac{X_{1t} - x_1}{h_1} \right)} \\
	= & \frac{\sum_{i=1}^n \left[ \frac{D_i}{\widehat{p}_i} \left( Y_i - \widehat{m}_{1i} \right) + \widehat{m}_{1i} \right]}{\sum_{t=1}^n K_1 \left( \frac{X_{1t} - x_1}{h_1} \right)} - \frac{\sum_{i=1}^n \left[ \frac{1 - D_i}{1 - \widehat{p}_i} \left( Y_i - \widehat{m}_{0i} \right) + \widehat{m}_{0i} \right]}{\sum_{t=1}^n K_1 \left( \frac{X_{1t} - x_1}{h_1} \right)} \\
	=: & \widehat{\tau}_1(x_1) - \widehat{\tau}_0(x_1).
\end{align*}

Let
\begin{align*}
	\tau(x_1) = & E \left[ \left. \frac{D}{p(X)}[Y-m_1(X)] - \frac{1-D}{1-p(X)}[Y-m_0(X)] + m_1(X) - m_0(X) \right| X_1 = x_1 \right] \\
	= & E \left[ \left. \frac{D}{p(X)}[Y-m_1(X)] + m_1(X) \right| X_1 = x_1 \right] \\
	& - E \left[ \left. \frac{1-D}{1-p(X)}[Y-m_0(X)] + m_0(X) \right| X_1 = x_1 \right] \\
	= & \tau_1(x) - \tau_0(x).
\end{align*}

For the very first move, we look for the asymptotic linear expression of $\sqrt{nh_1^k}[\widehat{\tau}(x_1) - \tau(x_1)]$. Note that
\begin{align*}
	& \sqrt{nh_1^k} \left[ \widehat{\tau}(x_1) - \tau(x_1) \right] \\
	= & \sqrt{nh_1^k} \left\{ \left[ \widehat{\tau}_1(x_1) - \tau_1(x_1) \right] - \left[ \widehat{\tau}_0(x_1) - \tau_0(x_1) \right] \right\} \\
	= & \frac{1}{\widehat{f}(x_1)} \frac{1}{\sqrt{nh_1^k}} \sum_{i=1}^n \left[ \frac{D_i}{\widehat{p}_i} \left( Y_i - \widehat{m}_{1i} \right) + \widehat{m}_{1i} - \tau_1(x_1) \right] K_1 \left( \frac{X_{1i} - x_1}{h_1} \right) \\
	& - \frac{1}{\widehat{f}(x_1)} \frac{1}{\sqrt{nh_1^k}} \sum_{i=1}^n \left[ \frac{1 - D_i}{1 - \widehat{p}_i} \left( Y_i - \widehat{m}_{0i} \right) + \widehat{m}_{0i} - \tau_0(x_1) \right] K_1 \left( \frac{X_{1i} - x_1}{h_1} \right).
\end{align*}
$\widehat{f}(x_1)\xrightarrow{ \mathbb{P} } f(x_1)$ ensures that we can use Slutsky's Theorem later. So we can first consider the asymptotic linear expression of
\begin{align}\label{J1}
	J(x_1) = \frac{1}{\sqrt{nh_1^k}} \sum_{i=1}^n \left[ \frac{D_i}{\widehat{p}_i} \left( Y_i - \widehat{m}_{1i} \right) + \widehat{m}_{1i} - \tau_1(x_1) \right] K_1 \left( \frac{X_{1i} - x_1}{h_1} \right).
\end{align}

Consider several combinations of estimation of nuisance functions. Now we list them as below.

{\normalsize
Scenario 1. $p(x)$ parametrically estimated (correctly specified), $m_1(x)$ parametrically estimated (correctly specified)

Scenario 2. $p(x)$ nonparametrically estimated, $m_1(x)$ nonparametrically estimated

Scenario 3. $p(x)$ semiparametrically estimated, $m_1(x)$ semiparametrically estimated

Scenario 4. $p(x)$ parametrically estimated (correctly specified), $m_1(x)$ nonparametrically estimated

Scenario 5. $p(x)$ parametrically estimated (correctly specified), $m_1(x)$ semiparametrically estimated

Scenario 6. $p(x)$ nonparametrically estimated, $m_1(x)$ parametrically estimated (correctly specified)

Scenario 7. $p(x)$ nonparametrically estimated, $m_1(x)$ semiparametrically estimated

Scenario 8. $p(x)$ semiparametrically estimated, $m_1(x)$ parametrically estimated (correctly specified)

Scenario 9. $p(x)$ semiparametrically estimated, $m_1(x)$ nonparametrically estimated}

\textbf{Scenario 1: $p(x)$ and $m_1(x)$ are  parametrically estimated.} From standard parametric estimation argument,
\begin{align*}
	\sup_{x \in \mathcal{X}} \left| \widetilde{p}(x;\widehat{\beta}) - p(x) \right| & = \sup_{x \in \mathcal{X}} \left| \widetilde{p}(x;\widehat{\beta})-\widetilde{p}(x;\beta^*) \right| = O_p \left( \frac{1}{\sqrt{n}} \right), \\
	\sup_{x \in \mathcal{X}} \left| \widetilde{m}_1(x;\widehat{\gamma}_1) - m_1(x) \right| & = \sup_{x \in \mathcal{X}} \left| \widetilde{m}_1(x;\widehat{\gamma}_1) - \widetilde{m}_1(x; \gamma_1^*) \right| = O_p \left( \frac{1}{\sqrt{n}} \right).
\end{align*}

We start from (\ref{J1}):
\begin{align}\label{J2}
	J(x_1) = & \frac{1}{\sqrt{nh_1^k}} \sum_{i=1}^n \left[ \frac{D_i}{\widetilde{p}(X_i; \widehat{\beta})} \left[ Y_i - \widetilde{m}_1(X_1; \widehat{\gamma}_1) \right] + \widetilde{m}_1(X_1; \widehat{\gamma}_1) - \tau_1(x_1) \right] K_1 \left( \frac{X_{1i} - x_1}{h_1} \right) \notag \\
	= & \frac{1}{\sqrt{nh_1^k}} \sum_{i=1}^n \left[ \frac{D_i Y_i}{p(X_i)} - \frac{D_i - p(X_i)}{p(X_i)} m_1(X_i) - \tau_1(x_1) \right] K_1 \left( \frac{X_{1i} - x_1}{h_1} \right) \notag \\
	& + \frac{1}{\sqrt{nh_1^k}} \sum_{i=1}^n \left[ \frac{D_i (m_{1i}^{+} - Y_i)}{{p_i^+}^2} \right] \left[ \widetilde{p}(X_i ; \widehat{\beta}) - p(X_i) \right] K_1 \left( \frac{X_{1i} - x_1}{h_1} \right) \notag \\
	& + \frac{1}{\sqrt{nh_1^k}} \sum_{i=1}^n \left( \frac{p_i^{+} - D_i}{p_i^{+}} \right) \left[ \widetilde{m}_1(X_i;\widehat{\gamma}_1) - m_1(X_i)\right]  K_1 \left( \frac{X_{1i} - x_1}{h_1} \right) \notag \\
	=:  & J_{11}(x_1)+J_{12}(x_1)+J_{13}(x_1)
\end{align}
where $p_i^+$ lies between $p(X_i)$ and $\widetilde{p}(X_i;\widehat{\beta})$, $m_{1i}^+$ lies between $m_1(X_i)$ and $\widetilde{m}_1(X_i;\widehat{\gamma}_1)$.

Bounding $J_{12}(x)$ as
\begin{align*}
	|J_{12}(x)| = & \frac{1}{\sqrt{nh_1^k}} \left| \sum_{i=1}^n \left[ \frac{D_i (m_{1i}^+ - Y_i) }{{p_i^+}^2} \right] \left[ \widetilde{p}(X_i;\widehat{\beta}) - p(X_i) \right] K_1 \left( \frac{X_{1i} - x_1}{h_1} \right) \right| \\
	\le & \sqrt{nh_1^k} \sup_{x \in \mathcal{X}} \left| \widetilde{p}(x;\widehat{\beta}) - p(x) \right| \times \frac{1}{nh_1^k} \sum_{i=1}^n \left| K_1 \left( \frac{X_{1i} - x_1}{h_1} \right) \right| \left| \frac{D_i ( m_{1i}^{+} - Y_i ) }{{p_i^+}^2} \right|
\end{align*}
where $\sup_{x \in \mathcal{X}} \left| \widetilde{p}(x;\widehat{\beta}) - p(x) \right| = O_p \left( \frac{1}{\sqrt{n}} \right)$, $\left| \frac{D_i ( m_{1i}^{+} - Y_i ) }{{p_i^+}^2} \right|$ is bounded due to condition (C2)(\romannumeral2), (C4) and (C5), $\frac{1}{nh_1^k} \sum_{i=1}^n \left| K_1 \left( \frac{X_{1i} - x_1}{h_1} \right) \right| = O_p(1)$ by the standard nonparametric estimation argument. Thus $|J_{12}(x_1)| \le o_p(1)\cdot O_p(1) = o_p(1)$.
With the similar arguments, we can also bound the last term as $|J_{13}(x_1)| = o_p(1)$. So far, we've proved $J_{12}(x_1)$ and $J_{13}(x_1)$ converge to 0 in probability. Hence, according to Slutsky's Theorem, together with (\ref{J2}), we have
\begin{align}\label{J3}
	& \sqrt{nh_1^k} [\widehat{\tau}_1(x_1) - \tau_1(x_1) ] = \frac{1}{\widehat{f}(x_1)} J(x_1) = \frac{1}{\widehat{f}(x_1)} J_{11}(x_1) + o_p(1) \notag \\
	= & \frac{1}{\widehat{f}(x_1)} \frac{1}{\sqrt{nh_1^k}} \sum_{i=1}^n \left[ \frac{D_i}{p(X_i)} [Y_i - m_1(X_i)] + m_1(X_i) - \tau_1(x_1) \right] K_1 \left( \frac{X_{1i} - x_1}{h_1} \right) + o_p(1)
\end{align}

\textbf{Scenario 2: $p(x)$ and $m_1(x)$ are nonparametrically estimated.} From the standard nonparametric estimation argument, under conditions (A1)(\romannumeral1), (\romannumeral2), (\romannumeral3), we have
\begin{align*}
	\sup_{x \in \mathcal{X}} \left| \widehat{p}(x) - p(x) \right| & = O_p \left( h_2^{s_2} + \sqrt{\frac{\ln n}{nh_2^p}} \right)= o_p \left( h_2^{\frac{s_2}{2}} \right) \\
	\sup_{x \in \mathcal{X}} \left| \widehat{m}_1(x) - m_1(x) \right| & = O_p \left( h_3^{s_3} + \sqrt{\frac{\ln n}{nh_3^p}} \right) = o_p \left( h_3^{\frac{s_3}{2}} \right)
\end{align*}

Rewrite (\ref{J1}):
\begin{align}\label{J4}
	J(x_1) = & \frac{1}{\sqrt{nh_1^k}} \sum_{i=1}^n \left[ \frac{D_i}{\widehat{p}(X_i)} \left[ Y_i - \widehat{m}_1(X_i) \right] + \widehat{m}_1(X_i) - \tau_1(x_1) \right] K_1 \left( \frac{X_{1i} - x_1}{h_1} \right) \notag \\
	= & \frac{1}{\sqrt{nh_1^k}} \sum_{i=1}^n \left[ \frac{D_i}{p(X_i)} [Y_i - m_1(X_i)] + m_1(X_i) - \tau_1(x_1) \right] K_1 \left( \frac{X_{1i} - x_1}{h_1} \right) \notag \\
	& + \frac{1}{\sqrt{nh_1^k}} \sum_{i=1}^n \frac{D_i [m_1(X_i)-Y_i]}{p^2(X_i)} \left[ \widehat{p}(X_i) - p(X_i) \right] K \left( \frac{X_{1i} - x_1}{h_1} \right) \notag \\
	& + \frac{1}{\sqrt{nh_1^k}} \sum_{i=1}^n \frac{p(X_i) - D_i}{p(X_i)} \left[ \widehat{m}_1(X_i) - m_1(X_i) \right] K_1 \left( \frac{X_{1i} - x_1}{h_1} \right) + 0 \notag \\
	& + \frac{1}{\sqrt{nh_1^k}} \sum_{i=1}^n \frac{D_i \left[ Y_i - m_{1i}^+ \right]}{{p_i^+}^3} \left[ \widehat{p}(X_i) - p(X_i) \right]^2 K_1 \left( \frac{X_{1i} - x_1}{h_1} \right) \notag \\
	& + \frac{1}{\sqrt{nh_1^k}} \sum_{i=1}^n \frac{2D_i}{{p_i^+}^2} \left[ \widehat{p}(X_i) - p(X_i) \right] \left[ \widehat{m}_1(X_i) - m_1(X_i) \right] K_1 \left( \frac{X_{1i} - x_1}{h_1} \right) \notag \\
	=: \, & \,  J_{21}(x_1) + J_{22}(x_1) + J_{23}(x_1) + 0 + J_{24}(x_1) + J_{25}(x_1)
\end{align}
where $p_i^+$ lies between $p(X_i)$ and $\widehat{p}(X_i)$, $m_{1i}^+$ lies between $m_1(X_i)$ and $\widehat{m}_1(X_i)$.

Rewrite $J_{22}(x_1)$ as
\begin{align*}
	J_{22}(x_1) = & \frac{1}{\sqrt{nh_1^k}} \sum_{i=1}^n \frac{D_i [m_1(X_i)-Y_i]}{p^2(X_i)} \left[ \widehat{p}(X_i) - p(X_i) \right] K_1 \left( \frac{X_{1i} - x_1}{h_1} \right) \\
	= &  \frac{1}{\sqrt{n}} \sum_{i=1}^n \frac{D_i [m_1(X_i)-Y_i]}{p^2(X_i)} \frac{1}{\sqrt{h_1^k}} \left[ \widehat{p}(X_i) - p(X_i) \right] K_1 \left( \frac{X_{1i} - x_1}{h_1} \right)
\end{align*}
In which $E \left[ \left. \frac{D_i [m_1(X_i)-Y_i]}{p^2(X_i)} \right| X_i \right] = 0$ and thus $\frac{D_i [m_1(X_i)-Y_i]}{p^2(X_i)}$ is independent of $X_i$ for every $i$; $\sup_{x \in \mathcal{X}} \left| \widehat{p}(x) - p(x) \right| = O_p \left( h_2^{s_2} + \sqrt{\frac{\ln n}{nh_2^p}} \right)= o_p \left( h_2^{\frac{s_2}{2}} \right) $, By  condition (A1)(\romannumeral8) and CLT, $\frac{1}{\sqrt{h_1^p}} \left[ \widehat{p}(X_i) - p(X_i) \right] K_1 \left( \frac{X_{1i} - x_1}{h_1} \right) = o_p(1)$.  and then $J_{22}(x_1)=o_p(1)$.

Similarly, $|J_{23}(x_1)| = o_p(1)$. Deal with  $J_{24}(x_1)$ by using the decomposition as
\begin{align*}
	|J_{24}(x_1)| = & \left| \frac{1}{\sqrt{nh_1^k}} \sum_{i=1}^n \frac{D_i \left[ Y_i - m_{1i}^+ \right]}{{p_i^+}^3} \left[ \widehat{p}(X_i) - p(X_i) \right]^2 K_1 \left( \frac{X_{1i} - x_1}{h_1} \right) \right| \\
	\le & \sqrt{nh_1^k} \sup_{x \in \mathcal{X}} \left[ \widehat{p}(X_i) - p(X_i) \right]^2 \times \frac{1}{nh_1^k} \sum_{i=1}^n \left| K_1 \left( \frac{X_{1i} - x_1}{h_1} \right) \right| \left| \frac{D_i \left[ Y_i - m_{1i}^+ \right]}{{p_i^+}^3} \right|
\end{align*}
in which $\sup_{x \in \mathcal{X}} \left[ \widehat{p}(x) - p(x) \right] =  o_p \left( h_2^{s_2} \right)$. Then under condition (A1)(\romannumeral8), $\sqrt{nh_1^k} \sup_{x \in \mathcal{X}} \left[ \widehat{p}(x) - p(x) \right]^2 = o_p(1)$. Under conditions (C2)(\romannumeral2) and (C5), $\left| \frac{D_i \left[ Y_i - m_{1i}^+ \right]}{{p_i^+}^3} \right|$ is bounded. Again by the standard argument for handling  nonparametric estimation, $ \frac{1}{nh_1^k} \sum_{i=1}^n \left| K \left( \frac{X_{1i} - x_1}{h_1} \right) \right| = O_p(1)$. Thus, we can obtain that $|J_{24}(x_1)|=o_p(1) \cdot O_p(1) = o_p(1)$. In a similar way, $|J_{25}(x_1)|=o_p(1)$ can also be proved. Here we have derived that $J_{22}(x_1)$, $J_{23}(x_1)$, $J_{24}(x_1)$ and $J_{25}(x_1)$ can be bounded as $o_p(1)$. Together with (\ref{J4}), we can obtain that
\begin{align}\label{J5}
	& \sqrt{nh_1^k} \left[ \widehat{\tau}_1(x_1) - \tau_1(x_1) \right] = \frac{1}{\widehat{f}(x_1)} J(x_1) = \frac{1}{\widehat{f}(x_1)} J_{21}(x_1) + o_p(1) \notag \\
	= & \frac{1}{\widehat{f}(x_1)} \frac{1}{\sqrt{nh_1^k}} \sum_{i=1}^n \left[ \frac{D_i}{p(Z_i)} [Y_i - m_1(X_i)] + m_1(X_i) - \tau_1(x_1) \right] K_1 \left( \frac{X_{1i} - x_1}{h_1} \right) \notag \\
&+ o_p(1).
\end{align}

\textbf{Scenario 3: $p(x)$ and $m_1(x)$ are semiparametrically estimated.} Under conditions (A2)(\romannumeral5) and (\romannumeral6), (\romannumeral7),
\begin{align*}
	\sup_{x\in \mathcal{X}} \left| \widehat{g}(A^{\top} x) - g(A^{\top} z) \right| & = O_p \left( h_5^{s_5}+\sqrt{\frac{ln(n)}{n h_5^{p(2)}}} \right) = o_p \left( h_5^\frac{s_5}{2} \right), \\
	\sup_{x \in \mathcal{X}} \left| \widehat{r}_1(B_1^{\top} z) - r_1(B_1^{\top} z) \right| & = O_p \left(h_6^{s_6} + \sqrt{\frac{\ln n}{n h_6^{p(1)}}} \right) = o_p \left( h_6^\frac{s_6}{2} \right).
\end{align*}
Note that under condition (B1), we can first discuss the asymptotic distribution by assuming that the projection matrices $A$, $B_0$ and $B_1$ are given. Then
\begin{align}\label{J6}
	J(x_1) = & \frac{1}{\sqrt{nh_1^k}} \sum_{i=1}^n \left[ \frac{D_i}{\widehat{g}(A^{\top} X_i)} \left[ Y_i - \widehat{r}_1(B_1^{\top} X_i) \right] + \widehat{r}_1(B_1^{\top} X_i) - \tau_1(x_1) \right] K_1 \left( \frac{X_{1i} - x_1}{h_1} \right) \notag \\
	= &  \frac{1}{\sqrt{nh_1^k}} \sum_{i=1}^n \left[ \frac{D_i}{p(X_i)} [Y_i - m_1(X_i)] + m_1(X_i) - \tau_1(x_1) \right] K_1 \left( \frac{X_{1i} - x_1}{h_1}\right) \notag \\
	& + \frac{1}{\sqrt{nh_1^k}} \sum_{i=1}^n \frac{D_i[m_1(X_i) - Y_i]}{p^2(X_i)} \left[ \widehat{g}(A^{\top} X_i) - g(A^{\top} X_i) \right] K_1 \left( \frac{X_{1i} - x_1}{h_1} \right) \notag \\
	& + \frac{1}{\sqrt{nh_1^k}} \sum_{i=1}^n \frac{p(X_i) - D_i}{p(X_i)} \left[ \widehat{r}_1(B_1^{\top} X_i) - r_1(B_1^{\top} X_i) \right] K_1 \left( \frac{X_{1i} - x_1}{h_1}\right) + 0 \notag \\
	& + \frac{1}{\sqrt{nh_1^k}} \sum_{i=1}^n \frac{D_i[Y_i - r_{1i}^+]}{{g_i^+}^3} \left[ \widehat{g}(A^{\top} X_i) - g(A^{\top} X_i) \right]^2 K_1 \left( \frac{X_{1i} - x_1}{h_1}\right) \notag \\
	& + \frac{1}{\sqrt{nh_1^k}} \sum_{i=1}^n \frac{2D_i}{{g_i^+}^2}\left[ \widehat{r}_1(B_1^{\top} X_i) - r_1(B_1^{\top} X_i) \right] \left[ \widehat{g}(A^{\top} X_i) - g(A^{\top} X_i) \right] K_1 \left(\frac{X_{1i} - x_1}{h_1}\right) \notag \\
	=: \,  & J_{31}(x_1)+J_{32}(x_1)+J_{33}(x_1)+0+J_{34}(x_1)+J_{35}(x_1)
\end{align}
where $g_i^+$ lies between $g(A^{\top} X_i)$ and $\widehat{g}(A^{\top} X_i)$, $m_{1i}^+$ lies between $r_1(B_1^{\top} X_i)$ and $\widehat{r}_1(B_1^{\top} X_i)$. Then we deal with all terms one by one.

Consider $J_{32}(x_1)$ and $J_{33}(x_1)$.  We have
\begin{align*}
	J_{32}(x_1) = & \frac{1}{\sqrt{nh_1^k}} \sum_{i=1}^n \frac{D_i[m_1(X_i) - Y_i]}{p^2(X_i)} \left[ \widehat{g}(A^{\top} X_i) - g(A^{\top} X_i) \right] K_1 \left( \frac{X_{1i} - x_1}{h_1} \right) \\
	= & \frac{1}{\sqrt{n}} \sum_{i=1}^n \frac{D_i[r_1(B_1^{\top} X_i) - Y_i]}{p^2(A^{\top} X_i)} \frac{1}{\sqrt{h_1^k}} \left[ \widehat{g}(A^{\top} X_i) - g(A^{\top} X_i) \right] K_1 \left( \frac{X_{1i} - x_1}{h_1} \right)
\end{align*}
Again $E \left[ \left. \frac{D_i[m_1(X_i) - Y_i]}{p^2(X_i)} \right| X_i \right] = 0$. Then, $\frac{D_i[m_1(X_i) - Y_i]}{p^2(X_i)}$ is independent of $X_i$ for every $i$; $\sup_{x \in \mathcal{X}} \left| \widehat{g}(A^{\top} X_i) - g(A^{\top} X_i) \right| = O_p \left( h_5^{s_5} + \sqrt{\frac{\ln n}{nh_5^p}} \right)= o_p \left( h_5^{\frac{s_5}{2}} \right) $. Condition (A2)(\romannumeral11) yields that  $\frac{1}{\sqrt{h_1^k}} \left[ \widehat{g}(A^{\top} X_i) - g(A^{\top} X_i) \right] K_1 \left( \frac{X_{1i} - x_1}{h_1} \right) = o_p(1)$. The application of CLT  yields that  $J_{32}(x_1) = o_p(1)$. Also, we can prove $J_{33}(x_1) = o_p(1)$ similarly.

Deal with $J_{34}(x_1)$ and $J_{35}(x_1)$. We have
\begin{align*}
	|J_{34}(x_1)| = & \left| \frac{1}{\sqrt{nh_1^k}} \sum_{i=1}^n \frac{D_i[Y_i - r_{1i}^+]}{{g_i^+}^3} \left[ \widehat{g}(A^{\top} X_i) - g(A^{\top} X_i) \right]^2 K_1 \left(\frac{X_{1i} - x_1}{h_1} \right) \right| \\
	\le & \sqrt{nh_1^k} \sup_{x \in \mathcal{X}} \left[ \widehat{g}(A^{\top} X_i)-g(A^{\top} X_i) \right]^2 \times \frac{1}{nh_1^k} \sum_{i=1}^n \left| K_1 \left( \frac{X_{1i} - x_1}{h_1} \right) \right| \left| \frac{D_i[Y_i - r_{1i}^+]}{{g_i^+}^3} \right|.
\end{align*}
Also $\sup_{x \in \mathcal{X}} \left[ \widehat{g}(A^{\top} X_i) - g(A^{\top} X_i) \right] =  o_p \left( h_6^{s_6} \right)$. Condition  (A2)(\romannumeral11) implies that
\begin{align*}
	\sqrt{nh_1^k} \sup_{x \in \mathcal{X}} \left[ \widehat{g}(A^{\top} X_i) - g(A^{\top} X_i) \right]^2 = o_p(1)
\end{align*}
Under conditions (C2)(\romannumeral2) and (C5), $\left| \frac{D_i \left[ Y_i - r_{1i}^+ \right]}{{g_i^+}^3} \right|$ is bounded. Again, $ \frac{1}{nh_1^k} \sum_{i=1}^n \left| K_1 \left( \frac{X_{1i} - x_1}{h_1} \right) \right| = O_p(1)$. We can then achieve $|J_{34}(x_1)| = o_p(1) \cdot O_p(1) = o_p(1)$. This is also the way to prove $|J_{35}(x_1)| = o_p(1)$. In this way, the asymptotic negligibility of $J_{32}(x_1)$, $J_{33}(x_1)$, $J_{34}(x_1)$ and $J_{35}(x_1)$ has been proved. Together with (\ref{J6}), it can be derived that
\begin{align}\label{J7}
	& \sqrt{nh_1^k} \left[ \widehat{\tau}_1(x_1) - \tau_1(x_1) \right] = \frac{1}{\widehat{f}(x_1)} J(x_1) = \frac{1}{\widehat{f}(x_1)} J_{31}(x_1) + o_p(1) \notag \\
	= & \frac{1}{\widehat{f}(x_1)} \frac{1}{\sqrt{nh_1^k}} \sum_{i=1}^n \left[ \frac{D_i}{p(X_i)} [Y_i - m_1(X_i)] + m_1(X_i) - \tau_1(x) \right] K_1 \left( \frac{X_{1i} - x_1}{h_1} \right) \notag \\
& + o_p(1).
\end{align}
Consider equations (\ref{J3}), (\ref{J5}) and (\ref{J7}), which imply that the asymptotic linear expressions of $\sqrt{nh_1^k}[\widehat{\tau}_1(x_1) - \tau_1(x_1)]$ are identical among scenarios 1, 2 and 3. It is obvious that under the conditions of Theorem 1, in any scenario mentioned above, the asymptotic linear expression remains the same, which leads to the same asymptotic distribution. With the asymptotic linear expression, we  can further derive the asymptotic distribution. First, we have the decomposition as
\begin{align}\label{J8}
	& \sqrt{nh_1^k}[\widehat{\tau}(x_1) - \tau(x_1)] \notag \\
	= &\frac{1}{\widehat{f}(x_1)} \frac{1}{\sqrt{nh_1^k}} \sum_{i=1}^n [\Psi_1(X_i, Y_i, D_i) - \tau(x_1)] K_1 \left( \frac{X_{1i} - x_1}{h_1} \right) + o_p(1) \notag \\
	= & \frac{1}{\widehat{f}(x_1)} \frac{1}{\sqrt{nh_1^k}} \sum_{i=1}^n [\Psi_1(X_i, Y_i, D_i) - \tau(X_{1i})] K_1 \left( \frac{X_{1i} - x_1}{h_1} \right) \notag \\
	& + \frac{1}{\widehat{f}(x_1)} \frac{1}{\sqrt{nh_1^k}} \sum_{i=1}^n [\tau(X_{1i}) - \tau(x_1)] K_1 \left( \frac{X_{1i} - x_1}{h_1} \right) + o_p(1) \notag \\
	=: \, & \frac{1}{\widehat{f}(x_1)}I_1(x_1) + \frac{1}{\widehat{f}(x_1)}I_2(x_1) + o_p(1).
\end{align}
Consider $I_1(x_1)$ first. Note that $\tau(X_{1i}) = E [ \Psi_1(X_i, Y_i, D_i) | X_1 = X_{1i} ]$. Then $\Psi_1(X_i, Y_i, D_i) - \tau(X_{1i})$ is independent of $X_{1i}$, $K_1 \left( \frac{X_{1i} - x_1}{h_1} \right)$ only depends on $n$ and $X_{1i}$. Thus
\begin{align*}
	E \left\{ [\Psi_1(X, Y, D) - \tau(X_1)] K_1 \left( \frac{X_1 - x_1}{h_1} \right) \right\} = E[\Psi_1(X, Y, D) - \tau(X_1)] E K_1 \left( \frac{X_1 - x_1}{h_1} \right) = 0.
\end{align*}
Also, $\left\{ [\Psi(X_i, Y_i, D_i) - \tau(X_{1i})] K_1 \left( \frac{X_{1i} - x_1}{h_1} \right) \right\}_{i=1}^n $ independently and identically distributed.
We now check  the condition of Lyapunov's CLT: $\exists \kappa >0$, s.t.
\begin{align*}
	\sum_{i=1}^n E \left| \frac{1}{\sqrt{n}} K_1 \left(\frac{X_{1i} - x_1}{h_1}\right)[\Psi_1(X_i, Y_i, D_i) - \tau(X_{1i})] \frac{1}{\sqrt{h_1^k}}\right|^{2+\kappa} \to 0 \qquad (n \to \infty)
\end{align*}
 Under condition (C6), letting $C = E|\Psi_1(X_i, Y_i, D_i) - \tau(X_{1i})|< \infty$, we have
\begin{align*}
	& \sum_{i=1}^n E \left| \frac{1}{\sqrt{n}}K_1(\frac{X_{1i} - x_1}{h_1})[\Psi_1(X_i, Y_i, D_i) - \tau(X_{1i})]\frac{1}{\sqrt{h_1^k}} \right|^{2+\kappa}  \\
	= & \left(\frac{1}{\sqrt{nh_1^k}}\right)^\kappa E \left|\Psi_1(X, Y, D) - \tau(X_1)|^{2+\kappa} E |K_1 \left(\frac{X_1 - x_1}{h_1} \right) \right|^{2+\kappa} \frac{1}{h_1^k} \\
	\le &  \left(\frac{1}{\sqrt{nh_1^k}}\right)^\kappa C E \left|K_1 \left( \frac{X_1 - x_1}{h_1}\right)\right|^{2+\kappa} \frac{1}{h_1^k}
\end{align*}
where
\begin{align*}
	\frac{1}{h_1^k} E \left| K_1 \left( \frac{X_1 - x_1}{h_1} \right) \right|^{2+\kappa} = \int K_1^{2+\kappa}(u)f(x_1 + h_1 u)du \to f(x_1)\int K_1^{2+\kappa}(u)du < \infty .
\end{align*}
Thus
\begin{align*}
	\left(\frac{1}{\sqrt{nh_1^k}}\right)^\kappa E \left| \Psi_1(X, Y, D) - \tau(X) \right|^{2+\kappa} E \left| K_1 \left( \frac{X_1 - x_1}{h_1} \right)\right|^{2+\kappa} \frac{1}{h_1^k} \to 0  (n \to \infty) .
\end{align*}
The Lyapunov's condition is satisfied and then
\begin{align}\label{J9}
	I_1(x) \xrightarrow{d} \mathcal{N}(0, V)
\end{align}
where $V=\lim_{n \to \infty} Var \left\{ \frac{1}{\sqrt{n h_1^k}} \sum_{i=1}^n [\Psi_1(X_i, Y_i, D_i) - \tau(x_1)] K_1 \left( \frac{X_{1i} - x_1}{h_1} \right) \right\}$. To compute the variance $V$, we can see that
\begin{align*}
	& Var \left\{ \frac{1}{\sqrt{n h_1^k}} \sum_{i=1}^n [\Psi_1(X_i, Y_i, D_i) - \tau(x_1)] K_1 \left( \frac{X_{1i} - x_1}{h_1} \right) \right\} \\
	= & E \left[ \frac{1}{\sqrt{n h_1^k}} \sum_{i=1}^n [\Psi_1(X_i, Y_i, D_i) - \tau(x_1)] K_1 \left( \frac{X_{1i} - x_1}{h_1} \right) \right]^2  \\
	= & h_1^k E\left\{ E \left[ \left. \left[ [\Psi_1(X, Y, D) - \tau(x_1)]\frac{1}{h_1^k} K_1 \left( \frac{X_1 - x_1}{h_1} \right) \right]^2 \right| X_1 \right] \right\} \\
	= & h_1^k \int \left( \frac{1}{h_1^k} K_1\left( \frac{X_1 - x_1}{h_1} \right) \right)^2 E \left[ \left.  [\Psi(X_i, Y_i, D_i) - \tau(x_1)]^2 \right| X \right] d F_{X_1} \\
	= & h_1^k \int \left( \frac{1}{h_1^k} K_1 \left( \frac{t - x_1}{h_1} \right) \right)^2 E \left[ \left.  [\Psi(X_i, Y_i, D_i) - \tau(x_1)]^2 \right| X_1 = t \right] f(t) dt \\
	= & h_1^k \frac{1}{h_1^k} \int K_1^2(u) E \left[ \left.  [\Psi_1(X, Y, D) - \tau(x_1)]^2 \right| X=x_1 + h_1 u \right] f(x_1 + h_1 u) du \\
	= & \sigma_1^2 (x_1) f(x_1) \int K_1^2(u) du + O(h_1^k)
\end{align*}
where
\begin{align*}
	\sigma_1^2(x_1) = E \left\{ \left. [\Psi_1(X,Y,D) - \tau(x_1)]^2 \right| X_1 = x_1 \right\}.
\end{align*}
Consider $I_2(x_1)$. We have
\begin{align*}
	I_2(x_1)= \frac{1}{\sqrt{nh_1^k}} \sum_{i=1}^n [\tau(X_{1i})-\tau(x_1)] K_1 \left(\frac{X_{1i} - x_1}{h_1} \right)
\end{align*}
where
\begin{align*}
	& E \left\{ \frac{1}{\sqrt{nh_1^k}} \sum_{i=1}^n [\tau(X_{1i})-\tau(x_1)] K_1 \left( \frac{X_{1i} - x_1}{h_1} \right) \right\} \\
	= & \sqrt{nh_1^k} \int [\tau(x_1 + h_1 u) - \tau(x_1)] K_1(u) f(x_1 + h_1 u)du
	=  \sqrt{nh_1^k}O_p(h_1^{s_1})
	=  o_p(1).
\end{align*}
Note that its variance is as
\begin{align*}
	& Var \left\{ \frac{1}{\sqrt{nh_1^k}} \sum_{i=1}^n [\tau(X_{1i}) - \tau(x_1)] K_1 \left(\frac{X_{1i} - x_1}{h_1}\right) \right\} \\
	= & E \left\{ \frac{1}{\sqrt{nh_1^k}} \sum_{i=1}^n [\tau(X_{1i}) -\tau(x_1)] K_1 \left(\frac{X_{1i} - x_1}{h_1}\right) \right\}^2 \\
	& - \left[ E \{ \frac{1}{\sqrt{nh_1^k}} \sum_{i=1}^n [\tau(X_{1i})-\tau(x_1)] K_1 \left(\frac{X_{1i} - x_1}{h_1}\right) \} \right]^2 \\
	= & \int [\tau(x_1 + h_1 u) - \tau(x)]^2 K_1^2(u) f(x_1 + h_1 u)du \\
	& - nh_1^k \left[ \tau(x_1 + h_1 u) - \tau(x_1)] K_1 (u) f(x_1 + h_1 u)du \right]^2
	=  o_p(1).
\end{align*}
so $I_2(x_1)=o_p(1)$.

Combining with (\ref{J8}), (\ref{J9}) and $I_2(x_1)=o_p(1)$, we can obtain that
\begin{align*}
	I_1(x_1) + I_2(x_1) \xrightarrow{d} \mathcal{N}\left( 0, \sigma_1^2 (x_1) f(x_1) \int K_1^2(u)du \right)
\end{align*}
and
\begin{align*}
	\sqrt{nh_1^k}[\widehat{\tau}(x_1) - \tau(x_1)] = \frac{1}{\widehat{f}(x_1)} [I_1(x_1) + I_2(x_1)] \xrightarrow{d} \mathcal{N}\left( 0, \frac{\sigma_1^2 (x_1) \int K_1^2(u)du }{f(x_1)}  \right).
\end{align*}

Now we consider the cases with unknown ${A}$, ${B}_1$ and and ${B}_0$. Note that under condition (B1), $\widehat{A}$, $\widehat{B}_1$ and $\widehat{B}_0$ converge in probability to $A$, $B_1$ and $B_0$ respectively at the rate of $O_p\left( \frac{1}{\sqrt{n}} \right)$. Following the similar arguments in \cite{hu2014estimation}, we can see easily that the asymptotic distribution retains. We then do not give the details for space saving.
Now we can conclude that, under the conditions of Theorem 1, regardless of which estimation method (parametric, nonparametric, semiparametric dimension reduction) used to estimate nuisance models,
\begin{align*}
	\sqrt{nh_1^l}[\widehat{\tau}(x_1) - \tau(x_1)] \xrightarrow{d} \mathcal{N}\left( 0, \frac{\sigma_1^2 (x_1) \int K_1^2(u)du }{f(x_1)}  \right).
\end{align*}
The proof is done.  \hfill$\Box$
\subsection{Proofs of Theorems 2 and 3}
We now consider global misspecification cases. Similar to the proof of Theorem 1, we first consider the asymptotic linear expression of $J(x_1)$.

\textbf{Scenario 1:  $m_1(x)$ is nonparametrically estimateed}. In this case, we have
\begin{align*}
	\sup_{x \in \mathcal{X}} \left| \widetilde{p}(x;\widehat{\beta})-\widetilde{p}(x;\beta^*) \right| & = O_p \left( \frac{1}{\sqrt{n}} \right), \\
	\sup_{x \in \mathcal{X}} \left| \widehat{m}_1(x) - m_1(x) \right| & = O_p \left( h_3^{s_3} + \sqrt{\frac{\ln n}{nh_3^p}} \right) = o_p \left( h_3^{\frac{s_3}{2}} \right).
\end{align*}

We can further rewrite (\ref{J1}) as
\begin{align}\label{J10}
	J(x_1) = & \frac{1}{\sqrt{nh_1^k}} \sum_{i=1}^n \left[ \frac{D_i}{\widetilde{p}(X_i; \widehat{\beta})} \left[ Y_i - \widehat{m}_1(X_i) \right] + \widehat{m}_1(X_i) - \tau_1(x_1) \right] K_1 \left( \frac{X_{1i} - x_1}{h_1} \right) \notag \\
	= & \frac{1}{\sqrt{nh_1^k}} \sum_{i=1}^n \left[ \frac{D_i}{\widetilde{p}(X_i;\beta^*)} [Y_i - m_1(X_i)] + m_1(X_i) - \tau_1(x_1) \right] K_1 \left( \frac{X_{1i} - x_1}{h_1} \right) \notag \\
	& + \frac{1}{\sqrt{nh_1^k}} \sum_{i=1}^n \frac{D_i [m_1(X_i)-Y_i]}{\widetilde{p}^2(X_i;\beta^*)} \left[ \widetilde{p}(X_i;\widehat{\beta}) - \widehat{p}(X_i;\beta^*) \right] K_1 \left( \frac{X_{1i} - x_1}{h_1} \right) \notag \\
	& + \frac{1}{\sqrt{nh_1^k}} \sum_{i=1}^n \frac{\widetilde{p}(X_i;\beta^*) - D_i}{\widetilde{p}(X_i;\beta^*)} \left[ \widehat{m}_1(X_i) - m_1(X_i) \right] K_1 \left( \frac{X_{1i} - x_1}{h_1} \right) + 0 \notag \\
	& + \frac{1}{\sqrt{nh_1^k}} \sum_{i=1}^n \frac{D_i \left[ Y_i - m_{1i}^+ \right]}{{p_i^+}^3} \left[ \widetilde{p}(X_i;\widehat{\beta}) - \widetilde{p}(X_i;\beta^*) \right]^2 K_1 \left( \frac{X_{1i} - x_1}{h_1} \right) \notag \\
	& + \frac{1}{\sqrt{nh_1^k}} \sum_{i=1}^n \frac{2D_i}{{p_i^+}^2} \left[ \widetilde{p}(X_i;\widehat{\beta}) - \widetilde{p}(X_i;\beta^*) \right] \left[ \widehat{m}_1(X_i) - m_1(X_i) \right] K_1 \left( \frac{X_{1i} - x_1}{h_1} \right) \notag \\
	=: \,  & J_{41}(x_1) + J_{42}(x_1) + J_{43}(x_1) + 0 + J_{44}(x_1) + J_{45}(x_1)
\end{align}
where $p_i^+$ lies between $\widetilde{p}(X_i;\beta^*)$ and $\widetilde{p}(X_i; \widehat{\beta})$, $m_{1i}^+$ lies between $m_1(X_i)$ and $\widehat{m}_1(X_i)$.

As we can prove $J_{42}(x_1)$ and $J_{44}(x_1)$ are $o_p(1)$ in the same way as the proof for  $J_{12}(x_1) = o_p(1)$ in scenario 1 of Subsection~\ref{proof: all correct}, the details are then omitted. For $J_{45}(x_1)$, obviously
\begin{align*}
	|J_{45}(x_1)| & = \left| \frac{1}{\sqrt{nh_1^k}} \sum_{i=1}^n \frac{2D_i}{{p_i^+}^2} \left[ \widetilde{p}(X_i;\widehat{\beta}) - \widetilde{p}(X_i;\beta^*) \right] \left[ \widehat{m}_1(X_i) - m_1(X_i) \right] K_1 \left( \frac{X_{1i} - x_1}{h_1} \right) \right| \\
	& = \sqrt{nh_1^k} \sup_{x \in \mathcal{X}} \left| \widetilde{p}(x;\widehat{\beta}) - \widetilde{p}(x;\beta^*) \right| \sup_{x \in \mathcal{X}} \left| \widehat{m}_1(x) - m_1(x) \right| \frac{1}{nh_1^k} \sum_{i=1}^n \left| K_1 \left( \frac{X_{1i} - x_1}{h_1} \right) \right| \left| \frac{2D_i}{{p_i^+}^2} \right| \\
	& = o_p(1).
\end{align*}

Consider $J_{43}(x_1)$. Denote that
\begin{align*}
	& \lambda_1(X_i) = E \left[ \left. \frac{\widetilde{p}(X_i; \beta^*) - D_i}{\widetilde{p}(X_i; \beta^*)} \right| X_i \right],  & \mu_{1i} = \frac{\widetilde{p}(X_i; \beta^*) - D_i}{\widetilde{p}(X_i; \beta^*)} - \lambda_1(X_i), \\
	& \epsilon_{1i} = Y_i - m_1(X_i), & \rho_{ij} = \frac{K_3 \left( \frac{X_i - X_j}{h_3} \right)}{\sum_{t=1}^n D_t K_3 \left( \frac{X_i - X_t}{h_3} \right)}.
\end{align*}

We first give a lemma to show their asymptotics below, which is useful for the proof of the theorem.
\begin{lemma}
\label{nonparametric outcome regression} Under condition (A1)(\romannumeral2), the outcome regression estimator satisfies
\begin{align*}
	| \rho_{ij} - \rho_{ji} | \leq \frac{C_n}{nh_3^p} K_3 \left( \frac{X_i - X_j}{h_3} \right)
\end{align*}
where $C_n = O_p(h_3)$ and does not depend on $i$, $j$.
\end{lemma}

\noindent \textbf{Proof.} Note that $\rho_{ij} = \rho_{ji} = 0$, if $|| X_i - X_j ||_{\infty} > h_3$.  We now consider the event that  $|| X_i - X_j ||_{\infty} \leq h_3$. For all $i$,
\begin{align*}
	\frac{1}{nh_3^p} \sum_{t=1}^n D_t K_3 \left( \frac{X_i - X_t}{h_3} \right) = \frac{\sum_{t=1}^n D_t K_3 \left( \frac{X_i - X_t}{h_3} \right)}{\sum_{t=1}^n K_3 \left( \frac{X_i - X_t}{h_3} \right)} \frac{1}{nh_3^p} \sum_{t=1}^n K_3 \left( \frac{X_i - X_t}{h_3} \right) = \widehat{p}(X_i) \widehat{\theta}(X_i).
\end{align*}
Then,
\begin{align*}
	|\rho_{ij} - \rho_{ji}| = & \frac{1}{nh_3^p} \left| K_3 \left( \frac{X_i - X_j}{h_3} \right) \right| \left| \widehat{p}^{-1}(X_i) \widehat{\theta}^{-1}(X_i) - \widehat{p}^{-1}(X_j) \widehat{\theta}^{-1}(X_j) \right| \\
	\leq & \frac{1}{nh_3^p} \left| K_3 \left( \frac{X_i - X_j}{h_3} \right) \right|\left\{ \left| \frac{1}{\widehat{p}(X_i) \widehat{\theta}(X_i)} - \frac{1}{p(X_i) \theta(X_i)} \right| + \left| \frac{1}{p(X_i) \theta(X_i)} - \frac{1}{p(X_j) \theta(X_j)} \right| \right.  \\
	& + \left. \left| \frac{1}{p(X_j) \theta(X_j)} - \frac{1}{\widehat{p}(X_j) \widehat{\theta}(X_j)}  \right| \right\}.
\end{align*}
Again by the standard arguments for dealing with nonparametric estimation as we have used  before and $s_3 \ge p$,
\begin{align*}
	\sup_{x \in \mathcal{X}} | \widehat{p}(x) - p(x) | & = O_p \left( h_3^{s_3} + \sqrt{\frac{\ln n}{n h_3^p}} \right) = o_p(h_3), \\
	\sup_{x \in \mathcal{X}} | \widehat{\theta}(x) - \theta(x) | & = O_p \left( h_3^{s_3} + \sqrt{\frac{\ln n}{n h_3^p}} \right) = o_p(h_3).
\end{align*}
Recall the conditions (C2)(\romannumeral2) and (C3)(\romannumeral1) that $p(x)$ and $\theta(x)$ are bounded away from $0$ and $1$. The two equations above implies that $\widehat{p}(x)$ and $\widehat{\theta}(x)$ uniformly converge to $p(x)$ and $\theta(x)$ respectively. We can also obtain that $\widehat{p}(x)$ and $\widehat{\theta}(x)$ are bounded away from $0$ in probability for $n$ large enough. Then
\begin{align*}	 & \sup_{x \in \mathcal{X}} \left| \frac{1}{p(x) \theta(x)} - \frac{1}{\widehat{p}(x) \widehat{\theta}(x)} \right| = \sup_{x \in \mathcal{X}} \frac{ \left| \widehat{p}(x) \widehat{\theta}(x) - p(x) \theta(x) \right|}{p(x) \theta(x) \widehat{p}(x) \widehat{\theta}(x)} \\
	\le & \sup_{x \in \mathcal{X}} \frac{ \widehat{p}(x) | \widehat{\theta}(x) - \theta(x) | + | \widehat{p}(x) - p(x) | \theta(x) }{p(x) \theta(x) \widehat{p}(x) \widehat{\theta}(x)} = o_p(1).
\end{align*}
This leads to that   $\left| \frac{1}{\widehat{p}(X_i) \widehat{\theta}(X_i)} - \frac{1}{p(X_i) \theta(X_i)} \right| = o_p(1)$, $\left| \frac{1}{\widehat{p}(X_j) \widehat{\theta}(X_j)} - \frac{1}{p(X_j) \theta(X_j)} \right| = o_p(1)$ uniformly over all $X_i, \, X_j$. By the Lipchitz continuity, we have  $\left| \frac{1}{p(X_i) \theta(X_i)} - \frac{1}{p(X_j) \theta(X_j)} \right| = O_p(h_3)$ uniformly over all $X_i, \, X_j$. Altogether,  we have that the summation in curly brace is $O_p(h_3)$. Therefore, there exists a $C_n = O_p(h_3)$ such that
\begin{align*}
	| \rho_{ij} - \rho_{ji} | \leq \frac{C_n}{nh_3^p} K_3 \left( \frac{X_i - X_j}{h_3} \right).	
\end{align*}
The proof  is completed. \hfill$\Box$

Now we come back to handle the term $J_{43}(x_1)$ that can be decomposed as
\begin{align}
	J_{43}(x_1) = & \frac{1}{\sqrt{nh_1^k}} \sum_{i=1}^n \frac{\widetilde{p}(X_i;\beta^*) - D_i}{\widetilde{p}(X_i;\beta^*)} \left[ \widehat{m}_1(X_i) - m_1(X_i) \right] K_1 \left( \frac{X_{1i} - x_1}{h_1} \right) \notag \\
	= & \frac{1}{\sqrt{nh_1^k}} \sum_{i=1}^n \lambda_1(X_i) \left[ \widehat{m}_1(X_i) - m_1(X_i) \right] K_1 \left( \frac{X_{1i} - x_1}{h_1} \right) \notag \\
	& + \frac{1}{\sqrt{nh_1^k}} \sum_{i=1}^n \mu_{1i} \left[ \widehat{m}_1(X_i) - m_1(X_i) \right] K_1 \left( \frac{X_{1i} - x_1}{h_1} \right) \notag \\
	= & \frac{1}{\sqrt{nh_1^k}} \sum_{i=1}^n K_1 \left( \frac{X_{1i} - x_1}{h_1} \right) \lambda_1(X_i) \left[ \sum_{j=1}^n \rho_{ij} D_j [\epsilon_{1j} + m_1(X_j)] - m_1(X_i) \right] \notag \\
	& + \frac{1}{\sqrt{nh_1^k}} \sum_{i=1}^n \mu_{1i} \left[ \widehat{m}_1(X_i) - m_1(X_i) \right] K_1 \left( \frac{X_{1i} - x_1}{h_1} \right) \notag \\
	= & \frac{1}{\sqrt{nh_1^k}} \sum_{i=1}^n K_1 \left( \frac{X_{1i} - x_1}{h_1} \right) \lambda_1(X_i) \frac{D_i \epsilon_{1i}}{p(X_i)} \notag \\
	& + \frac{1}{\sqrt{nh_1^k}} \sum_{i=1} \frac{D_i \epsilon_{1i}}{p(X_i)} \left[ p(X_i) \sum_{j=1}^n \rho_{ji} K_1 \left( \frac{X_{1j} - x_1}{h_1} \right) \lambda_1(X_j) - K_1 \left( \frac{X_{1i} - x_1}{h_1} \right) \lambda_1(X_i) \right] \notag \\
	& + \frac{1}{\sqrt{nh_1^k}} K_1 \left( X_{1i} - x_1 \right) \lambda_1(X_i) \left[ \sum_{j=1}^n \rho_{ij} D_j m_1(X_j) - m_1(X_i) \right] \notag \\
	& + \frac{1}{\sqrt{nh_1^k}} \sum_{i=1}^n \mu_{1i} \left[ \widehat{m}_1(X_i) - m_1(X_i) \right] K_1 \left( \frac{X_{1i} - x_1}{h_1} \right) \notag \\
	:= \, \,  & J_{431}(x_1) + J_{432}(x_1) + J_{433}(x_1) + J_{434}(x_1). \label{J43}
\end{align}

We first  prove that $J_{43k}(x_1)=o_p(1)$ for $k=2,3,4$. Consider  $J_{432}(x_1)$  by using the following decomposition:
\begin{align*}
	& \frac{1}{\sqrt{h_1^k}} \left[ p(X_i) \sum_{j=1}^n \rho_{ji} K_1 \left( \frac{X_{1j} - x_1}{h_1} \right) \lambda_1(X_j) - K_1 \left( \frac{X_{1i} - x_1}{h_1} \right) \lambda_1(X_i) \right]  \\
	= & \frac{1}{\sqrt{h_1^k}} p(X_i) \sum_{j=1}^n (\rho_{ij} - \rho_{ji}) K_1 \left( \frac{X_{1j} - x_1}{h_1} \right) \lambda_1(X_j) \\
	& + \frac{1}{\sqrt{h_1^k}} \left[ p(X_i) \sum_{j=1}^n \rho_{ij} K_1 \left( \frac{X_{1j} - x_1}{h_1} \right) \lambda_1(X_j) - K_1 \left( \frac{X_{1i} - x_1}{h_1} \right) \lambda_1(X_i) \right] \\
	:= \, & L_1 + L_2.
\end{align*}
 $L_1$ can be bounded as
\begin{align*}
	L_1 \leq & \frac{1}{\sqrt{h_1^k}} \sup_{i} \left| p(X_i) \sum_{j=1}^n (\rho_{ij} - \rho_{ji}) K_1 \left( \frac{X_{1j} - x_1}{h_1} \right) \lambda_1(X_j) \right| \\
	\leq & \frac{1}{\sqrt{h_1^k}} \sup_i \sum_{j=1}^n |\rho_{ij}-\rho_{ji}| \left| K_1 \left( \frac{X_{1j} - x_1}{h_1} \right) \lambda_1(X_j) \right| \\
	\leq & \frac{1}{\sqrt{h_1^k}} \frac{M C_n}{n h_3^p} \sup_i \sum_{j=1}^n |K_3 \left( \frac{X_i - X_j}{h_3} \right)| \left| K_1 \left( \frac{X_{1j} - x_1}{h_1} \right) \lambda_1(X_j) \right| \\
	\leq & \frac{M C_n}{h_3} \frac{h_3}{\sqrt{h_1^k}} \sup_i \frac{1}{nh_3^p} \sum_{j=1}^n |K_3 \left( \frac{X_i - X_j}{h_3} \right)| \left| K_1 \left( \frac{X_{1j} - x_1}{h_1} \right) \lambda_1(X_j) \right| \\
	= & O_p(1) \cdot o_p(1) \cdot O_p(1) \\
	= & o_p(1).
\end{align*}
Then  $L_2$ can be handled by noting
\begin{align*}
	& \frac{1}{\sqrt{h_1^k}} \left[ p(X_i) \sum_{j=1}^n \rho_{ij} K_1 \left( \frac{X_{1j} - x_1}{h_1} \right) \lambda_1(X_j) - K_1 \left( \frac{X_{1i} - x_1}{h_1} \right) \lambda_1(X_i) \right] \\
	= & \frac{1}{\sqrt{h_1^k}} \left[ p(X_i) \frac{\sum_{s=1}^n K_3 \left( \frac{X_i - X_s}{h_3} \right)}{\sum_{t=1}^n D_t K_3 \left( \frac{X_i - X_t}{h_3} \right)} \sum_{j=1}^n \frac{D_j K_3 \left( \frac{X_j - X_i}{h_3} \right)}{\sum_{t=1}^n D_t K_3 \left( \frac{X_i - X_t}{h_3} \right)} K_1 \left( \frac{X_{1j} - x_1}{h_1} \right) \lambda_1(X_j) \right. \\
	 & - \left. K_1 \left( \frac{X_{1i} - x_1}{h_1} \right) \lambda_1(X_i) \right] \\
	 = & \frac{1}{\sqrt{h_1^k}} \left\{ \left[ p(X_i) \frac{\sum_{s=1}^n K_3 \left( \frac{X_i - X_s}{h_3} \right)}{\sum_{t=1}^n D_t K_3 \left( \frac{X_i - X_t}{h_3} \right)} \sum_{j=1}^n \frac{D_j K_3 \left( \frac{X_j - X_i}{h_3} \right)}{\sum_{t=1}^n D_t K_3 \left( \frac{X_i - X_t}{h_3} \right)} K_1 \left( \frac{X_{1j} - x_1}{h_1} \right) \lambda_1(X_j) \right. \right. \\
	 & - \left. p(X_i) \frac{\sum_{s=1}^n K_3 \left( \frac{X_i - X_s}{h_3} \right)}{\sum_{t=1}^n D_t K_3 \left( \frac{X_i - X_t}{h_3} \right)} K_1 \left( \frac{X_{1i} - x_1}{h_1} \right) \lambda_1(X_i) \right] \\
	 & + \left[ p(X_i) \frac{\sum_{s=1}^n K_3 \left( \frac{X_i - X_s}{h_3} \right)}{\sum_{t=1}^n D_t K_3 \left( \frac{X_i - X_t}{h_3} \right)} K_1 \left( \frac{X_{1i} - x_1}{h_1} \right) \lambda_1(X_i) \right. \\
	 & - \left. \left. p(X_i) \frac{1}{p(X_i)} K_1 \left( \frac{X_{1i} - x_1}{h_1} \right) \lambda_1(X_i) \right] \right\} \\
	 = & \frac{1}{\sqrt{h_1^k}} O_p \left( \frac{h_3^{s_3}}{h_1^{s_3}} + h_3^{s_3} \right)
	 =   O_p \left( \frac{h_3^{s_3}}{h_1^{s_3 + k/2}} \right).
\end{align*}
Then under Condition (A2)(\romannumeral9), $L_2 = o_p(1)$ and,  together with the bound for $L_1$, we have
\begin{align*}
	\frac{1}{\sqrt{h_1^k}} \left[ p(X_i) \sum_{j=1}^n \rho_{ji} K_1 \left( \frac{X_{1j} - x_1}{h_1} \right) \lambda_1(X_j) - K_1 \left( \frac{X_{1i} - x_1}{h_1} \right) \lambda_1(X_i) \right] = o_p(1).
\end{align*}
Since $\{ \epsilon_{1i} \}_{i=1}^n$ are mutually independent given $\{ Z_i \}_{i=1}^n$, it follows that $J_{432}(x_1) = o_p(1)$.

Second, bound  $J_{433}(x_1)$ by noting that
\begin{align*}
	|J_{433}(x_1)| = & \left| \frac{1}{\sqrt{nh_1^k}} K_1 \left( X_{1i} - x_1 \right) \lambda_1(X_i) \left[ \sum_{j=1}^n \rho_{ij} D_j m_1(X_j) - m_1(X_i) \right] \right| \\
	\leq & \sqrt{nh_1^k} \sup_{x \in \mathcal{X}} \left| \sum_{j=1}^n \rho_{ij} D_j m_1(X_j) - m_1(X_i) \right| \frac{1}{nh_1^k} \sum_{i=1}^n \left| K_1 \left( X_{1i} - x_1 \right) \lambda_1(X_i) \right| \left| \lambda_1(X_i) \right| \\
	= & \sqrt{nh_1^k} O_p \left( h_3^{s_3} \right) O_p(1)
	=  o_p(1).
\end{align*}
A similar argument to bound $J_{23}(x_1)$ can lead to  $J_{434}(x_1) = o_p(1)$. Altogether and combining  (\ref{J43}), we have $J_{43}(x_1) = J_{431}(x_1) + o_p(1)$. Recalling that  $J_{42}(x_1)$, $J_{44}(x_1)$ and $J_{45}(x_1)$ have all been proved to be $o_p(1)$, together with (\ref{J10}), we can conclude that
\begin{align}
	J(x_1) = & J_{41}(x_1) + J_{431}(x_1) + o_p(1) \notag  \\
	= & \frac{1}{\sqrt{nh_1^k}} \sum_{i=1}^n \left[ \frac{D_i}{\widetilde{p}(X_i;\beta^*)} [Y_i - m_1(X_i)] + m_1(X_i) - \tau(x_1) \right] K_1 \left( \frac{X_{1i} - x_1}{h_1} \right) \notag \\
	& + \frac{1}{\sqrt{nh_1^k}} \sum_{i=1}^n K_1 \left( \frac{X_{1i} - x_1}{h_1} \right) \lambda_1(X_i) \frac{D_i \epsilon_{1i}}{p(X_i)} + o_p(1) \notag \\
	= & \frac{1}{\sqrt{nh_1^k}} \sum_{i=1}^n \left[ \frac{D_i}{p(X_i)} [Y_i - m_1(X_i)] + m_1(X_i) - \tau_1(x) \right] K_1 \left( \frac{X_{1i} - x_1}{h_1} \right) + o_p(1). \label{Jcon}
\end{align}
The proof is finished. \hfill$\Box$

\textbf{Scenario 2: $m_1(x)$ is semiparametrically estimated. }
First, we have
\begin{align*}
	\sup_{x \in \mathcal{X}} \left| \widetilde{p}(x;\widehat{\beta})-\widetilde{p}(x;\beta^*) \right| & = O_p \left( \frac{1}{\sqrt{n}} \right), \\
	\sup_{x \in \mathcal{X}} \left| \widehat{r}_1(B_1^{\top} x) - r_1(B_1^{\top} x) \right| & = O_p \left( h_6^{s_6} + \sqrt{\frac{\ln n}{nh_6^{p(1)}}} \right) = o_p \left( h_6^{\frac{s_6}{2}} \right),
\end{align*}
We can further decompose the term in  (\ref{J1}) as:
\begin{align}
	J(x_1) = & \frac{1}{\sqrt{nh_1^k}} \sum_{i=1}^n \left[ \frac{D_i}{\widetilde{p}(X_i; \widehat{\beta})} \left[ Y_i - \widehat{r}_1(B_1^{\top} X_i) \right] + \widehat{r}_1(B_1^{\top} X_i) - \tau_1(x_1) \right] K_1 \left( \frac{X_{1i} - x_1}{h_1} \right) \notag \\
	= & \frac{1}{\sqrt{nh_1^k}} \sum_{i=1}^n \left[ \frac{D_i}{\widetilde{p}(X_i;\beta^*)} [Y_i - m_1(X_i)] + m_1(X_i) - \tau(x_1) \right] K_1 \left( \frac{X_{1i} - x_1}{h_1} \right) \notag \\
	& + \frac{1}{\sqrt{nh_1^k}} \sum_{i=1}^n \frac{D_i [r_1(B_1^{\top} X_i)-Y_i]}{\widetilde{p}^2(X_i;\beta^*)} \left[ \widetilde{p}(X_i;\widehat{\beta}) - \widehat{p}(X_i;\beta^*) \right] K_1 \left( \frac{X_{1i} - x_1}{h_1} \right) \notag \\
	& + \frac{1}{\sqrt{nh_1^k}} \sum_{i=1}^n \frac{\widetilde{p}(X_i;\beta^*) - D_i}{\widetilde{p}(X_i;\beta^*)} \left[ \widehat{r}_1(B_1^{\top} X_i) - r_1(B_1^{\top} X_i) \right] K_1 \left( \frac{X_{1i} - x_1}{h_1} \right) + 0 \notag \\
	& + \frac{1}{\sqrt{nh_1^k}} \sum_{i=1}^n \frac{D_i \left[ Y_i - r_{1i}^+ \right]}{{p_i^+}^3} \left[ \widetilde{p}(X_i;\widehat{\beta}) - \widetilde{p}(X_i;\beta^*) \right]^2 K_1 \left( \frac{X_{1i} - x_1}{h_1} \right) \notag \\
	& + \frac{1}{\sqrt{nh_1^k}} \sum_{i=1}^n \frac{2D_i}{{p_i^+}^2} \left[ \widetilde{p}(X_i;\widehat{\beta}) - \widetilde{p}(X_i;\beta^*) \right] \left[ \widehat{r}_1(B_1^{\top} X_i) - r_1(B_1^{\top} X_i) \right] K_1 \left( \frac{X_{1i} - x_1}{h_1} \right) \notag \\
	:= \, \, & J_{51}(x_1) + J_{52}(x_1) + J_{53}(x_1) + 0 + J_{54}(x_1) + J_{55}(x_1) \label{J13}
\end{align}
where $p_i^+$ lies between $\widetilde{p}(X_i;\beta^*)$ and $\widetilde{p}(X_i; \widehat{\beta})$, $r_{1i}^+$ lies between $r_1(B_1^{\top} X_i)$ and $\widehat{r}_1(B_1^{\top} X_i)$.

Due to the similarity in the proof as the above,  we omit the details for proving $J_{52}(x_1)$, $J_{54}(x_1)$ and $J_{55}(x_1)$ to be $o_p(1)$. Now  consider $J_{53}(x_1)$. Denote that
\begin{align*}
	& \lambda_2(X_i) = E \left[ \left. \frac{\widetilde{p}(X_i; \beta^*) - D_i}{\widetilde{p}(X_i; \beta^*)} \right| X_i \right]  & \mu_{2i} = \frac{\widetilde{p}(X_i; \beta^*) - D_i}{\widetilde{p}(X_i; \beta^*)} - \lambda_2(X_i) \\
	& \epsilon_{2i} = Y_i - r_1(B_1^{\top} X_i) & \nu_{ij} = \frac{K_6 \left( \frac{B_1^{\top} X_i - B_1^{\top} X_j}{h_6} \right)}{\sum_{t=1}^n D_t K_6 \left( \frac{B_1^{\top} X_i - B_1^{\top} X_t}{h_6} \right)}
\end{align*}
$J_{53}(x_1)$ can be rewritten as
\begin{align}
	J_{53}(x_1) = & \frac{1}{\sqrt{nh_1^k}} \sum_{i=1}^n \frac{\widetilde{p}(X_i;\beta^*) - D_i}{\widetilde{p}(X_i;\beta^*)} \left[ \widehat{r}_1(B_1^{\top} X_i) - r_1(B_1^{\top} X_i) \right] K_1 \left( \frac{X_{1i} - x_1}{h_1} \right) \notag \\
	= & \frac{1}{\sqrt{nh_1^k}} \sum_{i=1}^n \lambda_2(X_i) \left[ \widehat{r}_1(B_1^{\top} X_i) - r_1(B_1^{\top} X_i) \right] K_1 \left( \frac{X_{1i} - x_1}{h_1} \right) \notag \\
	& + \frac{1}{\sqrt{nh_1^k}} \sum_{i=1}^n \mu_{2i} \left[ \widehat{r}_1(B_1^{\top} X_i) - r_1(B_1^{\top} X_i) \right] K_1 \left( \frac{X_{1i} - x_1}{h_1} \right) \notag \\
	= & \frac{1}{\sqrt{nh_1^k}} \sum_{i=1}^n K_1 \left( \frac{X_{1i} - x_1}{h_1} \right) \lambda_2(X_i) \left[ \sum_{j=1}^n \rho_{ij} D_j [\epsilon_{2j} + r_1(B_1^{\top} X_j)] - r_1(B_1^{\top} X_i) \right] \notag \\
	& + \frac{1}{\sqrt{nh_1^k}} \sum_{i=1}^n \mu_{2i} \left[ \widehat{r}_1(B_1^{\top} X_i) - r_1(B_1^{\top} X_i) \right] K_1 \left( \frac{X_{1i} - x_1}{h_1} \right) \notag \\
	= & \frac{1}{\sqrt{nh_1^k}} \sum_{i=1} \frac{D_i \epsilon_{2i}}{p(X_i)} \left[ p(X_i) \sum_{j=1}^n \nu_{ji} K_1 \left( \frac{X_{1j} - x_1}{h_1} \right) \lambda_2(X_j) \right] \notag \\
	& + \frac{1}{\sqrt{nh_1^k}} K_1 \left( \frac{X_{1j} - x_1}{h_1} \right) \lambda_2(X_i) \left[ \sum_{j=1}^n \nu_{ij} D_j m_1(X_j) - m_1(X_i) \right] \notag \\
	& + \frac{1}{\sqrt{nh_1^k}} \sum_{i=1}^n \mu_{2i} \left[ \widehat{m}_1(X_i) - m_1(X_i) \right] K_1 \left( \frac{X_{1i} - x_1}{h_1} \right) \notag \\
	:= & J_{531}(x_1) + J_{532}(x_1) + J_{533}(x_1). \label{J53}
\end{align}

It is obvious that $J_{532}(x_1) = o_p(1)$ and $J_{533}(x_1) = o_p(1)$. To derive that $J_{531}(x_1) = o_p(1)$, we start by writing
\begin{align*}
	& \frac{1}{\sqrt{h_1^k}} \left[ p(X_i) \sum_{j=1}^n \nu_{ji} K_1 \left( \frac{X_{1j} - x_1}{h_1} \right) \lambda_2(X_j) \right]  \\
	= & \frac{1}{\sqrt{h_1^k}} p(X_i) \sum_{j=1}^n (\nu_{ij} - \nu_{ji}) K_1 \left( \frac{X_{1j} - x_1}{h_1} \right) \lambda_2(X_j) \\
	& + \frac{1}{\sqrt{h_1^k}} \left[ p(X_i) \sum_{j=1}^n \nu_{ij} K_1 \left( \frac{X_{1j} - x_1}{h_1} \right) \lambda_2(X_j) \right] \\
	:= \, \, & L_3 + L_4.
\end{align*}
Similarly for proving $L_1 = o_p(1)$ above, we can show that $L_3 = o_p(1)$, and thus omit the details. To handle $L_4 = o_p(1)$, we denote $q_{B1}(z)$ as the density of $B_1^{\top} X$, and
\begin{align*}
	\theta_{B1}(B_1^{\top} X) & = E[Y(1) | B_1^{\top} X], \\
	\widehat{\theta}_{B1}(B_1^{\top} x) & = \frac{\sum_{j=1}^n D_j Y_j K_6 \left( \frac{B_1^{\top} X_j - B_1^{\top} x}{h_6} \right)}{\sum_{t=1}^n D_t K_6 \left( \frac{B_1^{\top} X_t - B_1^{\top} x}{h_6} \right)}, \\
	\widehat{q}_{B1}(B_1^{\top} x) & = \frac{\sum_{j=1}^n D_j K_6 \left( \frac{B_1^{\top} X_j - B_1^{\top} x}{h_6} \right)}{\sum_{t=1}^n K_6 \left( \frac{B_1^{\top} X_t - B_1^{\top} x}{h_6} \right)}.
\end{align*}
Let $T_1 = \frac{B_1^{\top} X - B_1^{\top} X_i}{h_6}$, $T_1 = \frac{X_1 - X_{1i}}{h_6}$, $T_3 = \frac{X - X_i}{h_6}$. To deal with  $L_4$, consider the conditional expectation that can be derived as:
\begin{align*}
	& E \left\{ \left. p(X_i) \sum_{j=1}^n \nu_{ij} K_1 \left( \frac{X_{1j} - x_1}{h_1} \right) \lambda_2(X_j) \right| X_i \right\} \\
	= & E \left\{ \left. p(X_i) \frac{\frac{1}{nh_6^{k(1)}} \sum_{j=1}^n K_6 \left( \frac{B_1^{\top} X_j - B_1^{\top} X_i}{h_6} \right) K_1 \left( \frac{X_{1j} - x_1}{h_1} \right) \lambda_2(X_j)}{\widehat{\theta}_{B1}(B_1^{\top} X_i) \widehat{q}_{B1}(B_1^{\top} X_i) } \right| X_i \right\} \\
	= & \frac{[1 + o_p(1)] p(X_i)}{h_6^{k(1)} \theta_{B1}(B_1^{\top} X_i) q_{B1}(B_1^{\top} X_i)} \int K_6 \left( \frac{B_1^{\top} u - B_1^{\top} X_i}{h_6} \right) K_1 \left(  \frac{u - x_1}{h_1} \right) \lambda_2(u) \theta(u) du \\
	= & \frac{h_6^p [1 + o_p(1)] p(X_i)}{h_6^{k(1)} \theta_{B1}(B_1^{\top} X_i) q_{B1}(B_1^{\top} X_i)} \int K_6 \left( t_1 \right) K_1 \left( \frac{X_{1i} - x_1}{h_1} + t2 \frac{h_6}{h_1} \right) \lambda_2(X_i + t_3 \frac{h_6}{h_1}) \theta(X_i + t_3 h_6) d t_3 \\
	= & h_6^{p - p(1)} \frac{p(X_i)}{q_{B1}(B_1^{\top} X_i)} \frac{\theta(X_i)}{\theta_{B1}(B_1^{\top} X_i)} K_1 \left( \frac{X_{1i} - x_1}{h_1} \right) \lambda_2(X_i) + O_p \left( \frac{h_6^{p - p(1) + s_6}}{h_1^{s_6}} \right) \\
	= & O_p \left(  h_6^{p - p(1)} + \frac{h_6^{p - p(1) + s_6}}{h_1^{s_6}} \right).
\end{align*}
Then, when $s_6 < (2 s_6 + k) (p - p(1))$, $L_4 = O_p \left( \frac{h_6^{p - p(1)}}{h_1^{l/2}} +  \frac{h_6^{p - p(1) + s_6}}{h_1^{s_6 + l/2}} \right) = o_p(1) $, $J_{531}(x_1) = o_p(1)$. Together with (\ref{J53}), $J_{53}(x_1) = o_p(1)$. Recall  that we have proved that $J_{52}(x_1)$, $J_{54}(x_1)$ and $J_{55}(x_1)$ can be bounded by $o_p(1)$. With (\ref{J5}), we can eventually derive the asymptotically linear representation as
\begin{align}
	J(x_1) = & J_{51}(x_1) + o_p(1) \notag \\
	= & \frac{1}{\sqrt{nh_1^k}} \sum_{i=1}^n \left[ \frac{D_i}{\widetilde{p}(X_i;\beta^*)} [Y_i - m_1(X_i)] + m_1(X_i) - \tau(x_1) \right] K_1 \left( \frac{X_{1i} - x_1}{h_1} \right) \nonumber\\
 &+ o_p(1). \label{Jlin2}
\end{align}
The proof is completed. \hfill$\Box$

\textbf{Scenario 3: $m_1(x)$ is parametrically estimated (correctly specified). }
With the similar argument for proving scenario 1 in Theorem 1, we can easily derive that
\begin{align}
	J(x_1) = & \frac{1}{\sqrt{nh_1^k}} \sum_{i=1}^n \left[ \frac{D_i}{\widetilde{p}(X_i;\beta^*)} [Y_i - m_1(X_i)] + m_1(X_i) - \tau(x_1) \right] K_1 \left( \frac{X_{1i} - x_1}{h_1} \right) + o_p(1). \label{Jlin3}
\end{align}
With (\ref{Jcon}), it will be easy to further deduct the asymptotically linear expression of the proposed estimator. When the outcome regression functions are nonparametrically estimated, recalling  the relation between the following and $\sqrt{nh_1^k} \left[ \widehat{\tau}(x_1) - \tau(x_1) \right]$ and $J(x_1)$ defined in (\ref{J1}), we can derive that
\begin{align*}
	\sqrt{nh_1^k}[\widehat{\tau}(x_1) - \tau(x_1)]
	= \frac{1}{\widehat{f}(x_1)} \frac{1}{\sqrt{nh_1^k}} \sum_{i=1}^n [\Psi_1(X_i, Y_i, D_i) - \tau(x_1)] K_1 \left( \frac{X_{1i} - x_1}{h_1} \right).
\end{align*}
According to (\ref{Jlin2}) and (\ref{Jlin3}), when the outcome regression functions are semiparametrically or parametrically estimated, we have a similar representation as
\begin{align*}
	 \sqrt{nh_1^k}[\widehat{\tau}(x_1) - \tau(x_1)]
	= &\frac{1}{\widehat{f}(x_1)} \frac{1}{\sqrt{nh_1^k}} \sum_{i=1}^n [\Psi_2(X_i, Y_i, D_i) - \tau(x_1)] K_1 \left( \frac{X_{1i} - x_1}{h_1} \right).
\end{align*}
Similarly as the proof for Theorem 1, we can derive that under the conditions of Theorem 2, when the outcome regression functions are estimated nonparametrically, we have
\begin{align*}
	\sqrt{nh_1^k} \left[ \widehat{\tau}(x_1) - \tau(x_1) \right] \xrightarrow{d} \mathcal{N}\left( 0, \frac{\sigma_1^2 (x_1) \int K_1^2(u)du }{f(x_1)} \right),
\end{align*}
and when the outcome regression functions are estimated semiparametrically or parametrically, we have
\begin{align*}
	\sqrt{nh_1^k} \left[ \widehat{\tau}(x_1) - \tau(x_1) \right] \xrightarrow{d} \mathcal{N}\left( 0, \frac{\sigma_2^2 (x_1) \int K_1^2(u)du }{f(x_1)} \right).
\end{align*}
The proof of Theorem 2 is concluded. \hfill$\Box$

As for Theorem 3. the proof can be very similar to the proof for Theorem 2. Here we only give a crucial lemma in this proof and omit the details of the proof.

\begin{lemma}
\label{nonparametric propensity score} Under condition (A2)(\romannumeral8), the propensity score estimator satisfies
\begin{align*}
	| \omega_{ij} - \omega_{ji} | \leq \frac{E_n}{nh_2^p} K_2 \left( \frac{X_i - X_j}{h_2} \right)
\end{align*}
where $E_n = O_p(h_2)$  free of  $i$ and $j$.
\end{lemma}

\subsection{Proof of Theorem 4}

This is the case with local misspecification. To check the asymptotic efficiency through the variance comparison, we now compute the difference between $\sigma_1^2(x_1)$ and $\sigma_2^2(x_1)$:
\begin{align}
	& \sigma_2^2 (x_1) - \sigma_1^2 (x_1) \notag \\
	= & E \left\{ \left. \frac{p(X)-\widetilde{p}(X;\beta^*)}{\left[ \widetilde{p}(X;\beta^*) \right]^2} Var(Y|D=1,X) + \frac{\widetilde{p}(X;\beta^*)-p(X)}{\left[ 1-\widetilde{p}(X;\beta^*) \right]^2} Var(Y|D=0,X) \right| X_1 = x_1 \right\} \notag \\
	= & E \left\{  \frac{p(X)-\widetilde{p}(X;\beta_0)+\widetilde{p}(X;\beta_0)-\widetilde{p}(X;\beta^*)}{\left[ \widetilde{p}(X;\beta^*) \right]^2} Var(Y|D=1,X) \right. \notag \\
	& + \left. \left. \frac{\widetilde{p}(X;\beta^*)-\widetilde{p}(X;\beta_0)+\widetilde{p}(X;\beta_0)-p(X)}{\left[ 1-\widetilde{p}(X;\beta^*) \right]^2} Var(Y|D=0,X) \right| X_1 = x_1 \right\} \label{sd2-1}
\end{align}
and the difference between $\sigma_1^2(x_1)$ and $\sigma_3^2(x_1)$:
\begin{align}
	& \sigma_3^2 (x_1) - \sigma_1^2 (x_1) \notag \\
	= & E \left\{ \left. \left[ \left( 1-\frac{D}{p(X)} \right) \left[ \widetilde{m}_1(X;\gamma_1^*)-m_1(X) \right] - \left( 1-\frac{1-D}{1-p(X)}\right) \left[ \widetilde{m}_0(X;\gamma_0^*) - m_0(X) \right]  \right]^2 \right| X_1 = x_1 \right\} \notag \\
	= & E \left\{ \frac{1-p(X)}{p(X)} \left[ \widetilde{m}_1(X;\gamma_1^*)-\widetilde{m}_1(X;\gamma_{10})+\widetilde{m}_1(X;\gamma_{10})-m_1(X) \right]^2 \right. \notag \\
	& + \frac{p(X)}{1-p(X)} \left[ \widetilde{m}_0(X;\gamma_0^*)-\widetilde{m}_0(X;\gamma_{00})+\widetilde{m}_0(X;\gamma_{00})-m_0(X) \right]^2 \notag  \\
	& + \left[ \widetilde{m}_1(X;\gamma_1^*)-\widetilde{m}_1(X;\gamma_{10})+\widetilde{m}_1(X;\gamma_{10})-m_1(X) \right] \notag \\
	& \times \left. \left. \left[ \widetilde{m}_0(X;\gamma_0^*) -\widetilde{m}_0(X;\gamma_{00})+\widetilde{m}_0(X;\gamma_{00}) - m_0(X) \right] \right| X_1 = x_1 \right\}. \label{sd3-1}
\end{align}
Recall  that as the definitions, for all $x \in \mathcal{X}$, there exists $\beta_0,\gamma_{10},\gamma_{00}$, such that
\begin{align*}
	p(x) & = \widetilde{p}(x;\beta_0)[1+c_n a(x)], \\
	m_1(x) & = \widetilde{m}_1(x;\gamma_{10}) + d_{1n} b_1(x), \\
	m_0(x) & = \widetilde{m}_0(x;\gamma_{00}) + d_{0n} b_0(x).
\end{align*}
That is, $p(x) - \widetilde{p}(x;\beta_0) = O(c_n)$, $m_1(x) - \widetilde{m}_1(x;\gamma_{10}) = O(d_{1n})$, and $m_0(x) - \widetilde{m}_0(x;\gamma_{00}) = O(d_{0n})$. So now we only need to consider $\widetilde{p}(x;\beta_0) - \widetilde{p}(x;\beta^*)$, $\widetilde{m}_1(x;\gamma_{10}) - \widetilde{m}_1(x;\gamma_1^*)$ and $\widetilde{m}_0(x;\gamma_{00}) - \widetilde{m}_0(x;\gamma_0^*)$.
Note that $\beta^*, \gamma_1^*, \gamma_0^*$ are the limits of  the maximum likelihood estimators $\widehat{\beta},\widehat{\gamma_1},\widehat{\gamma_0}$ respectively. Discuss $\beta^*$ first. Given the propensity score function, $D$ is bernoulli distributed. We can respectively obtain, as  the propensity score function would be misspecified, the quasi-likelihood function and the quasi-log likelihood function of the unknown parameter $\beta$:
\begin{align*}
	\widetilde{L}(\beta) = \prod_{i=1}^n \widetilde{p}(X_i;\beta)^{D_i} [ 1 - \widetilde{p}(X_i;\beta) ]^{1-D_i} f(X_i),
\end{align*}
and
\begin{align*}
	\widetilde{l}(\beta) = \sum_{i=1}^n D_i \ln \widetilde{p}(X_i;\beta) + (1-D_i) \ln [ 1 - \widetilde{p}(X_i;\beta) ] + \ln f(X_i).
\end{align*}
Then, $\widehat{\beta}$ and $\beta^*$ satisfy that
\begin{align*}
	\widehat{\beta} = \mathop{argmax}_\beta \frac{1}{n} \widetilde{l}(\beta),
\qquad
	\beta^* = \mathop{argmax}_\beta E \left[ g(W;\beta) \right].
\end{align*}
where $g(W;\beta) = D \ln \widetilde{p}(X;\beta) + (1-D) \ln [ 1 - \widetilde{p}(X;\beta) ] + \ln f(X).$
By the mean value theorem,
\begin{align*}
	E \left[ \left. \frac{\partial g(W,\beta)}{\partial \beta} \right|_{\beta = \beta_0} \right] - E \left[ \left. \frac{\partial g(W,\beta)}{\partial \beta} \right|_{\beta = \beta^*} \right] = E \left[ \left. \frac{\partial^2 g(W,\beta)}{\partial \beta \partial \beta^{\top}} \right|_{\beta = \widetilde{\beta}} \right] (\beta^* - \beta_0),
\end{align*}
and
\begin{align*}
	E \left[ \left. \frac{\partial^2 g(W,\beta)}{\partial \beta \partial \beta^{\top}} \right|_{\beta = \widetilde{\beta}} \right] (\beta^* - \beta_0) = E \left[ \left. \frac{\partial g(W,\beta)}{\partial \beta} \right|_{\beta = \beta_0} \right]
\end{align*}
where $\widetilde{\beta}$ takes the value between $\beta_0$ and $\beta^*$. Note that
\begin{align*}
	& E \left[ \left. \frac{\partial g(W,\beta)}{\partial \beta} \right|_{\beta = \beta_0} \right] \\
	= & E \left[ \frac{D[1+c_n a(X)]}{p(X)} \left. \frac{\partial \widetilde{p}(X;\beta)}{\partial \beta} \right|_{\beta = \beta_0} - \frac{1-D}{1-\widetilde{p}(X;\beta)} \left. \frac{\partial \widetilde{p}(X;\beta)}{\partial \beta} \right|_{\beta = \beta_0} \right] \\
	= & E \left[ (1+c_n a(X)) \left. \frac{\partial \widetilde{p}(X;\beta)}{\partial \beta} \right|_{\beta = \beta_0} - \frac{1-p(X)}{1-p(X)/[1+c_n a(X)]} \left. \frac{\partial \widetilde{p}(X;\beta)}{\partial \beta} \right|_{\beta = \beta_0} \right] \\
	= & E \left[ \frac{c_n a(X) + c_n^2 a^2(X)}{[1+c_n a(X) - p(X)} \left. \frac{\partial \widetilde{p}(X;\beta)}{\partial \beta} \right|_{\beta = \beta_0} \right] \\
	= & O(c_n).
\end{align*}
Assume that $E \left[ \frac{\partial^2 \ln g(U,\beta)}{\partial \beta \partial \beta^{\top}}  \right]$ is non-singular for any $\beta$. We have
\begin{align*}
	& \beta^* - \beta_0 \\
	= & \left\{ E \left[ \left. \frac{\partial^2 g(W,\beta)}{\partial \beta \partial \beta^{\top}} \right|_{\beta = \widetilde{\beta}} \right] \right\}^{-1} E \left[ \left. \frac{\partial g(W,\beta)}{\partial \beta} \right|_{\beta = \beta_0} \right]\\
	= & \left\{ E \left[ \left. \frac{\partial^2 g(W,\beta)}{\partial \beta \partial \beta^{\top}} \right|_{\beta = \widetilde{\beta}} \right] \right\}^{-1} O(c_n)
	=  O(c_n).
\end{align*}
The application of Taylor expansion yields that $\widetilde{p}(x;\beta_0) - \widetilde{p}(x;\beta^*) = O(c_n)$. Similar argument is devoted to deriving that $\widetilde{m}_1(x;\gamma_{10}) - \widetilde{m}_1(x;\gamma_1^*) = O(d_{1n})$ and $\widetilde{m}_0(x;\gamma_{00}) - \widetilde{m}_0(x;\gamma_0^*) = O(d_{0n})$.
Together with these results, we continue to calculate the quantities in (\ref{sd2-1}) and (\ref{sd3-1}) to derive that
\begin{align*}
	\sigma_2^2 (x_1) - \sigma_1^2 (x_1) &= O(c_n),\\
	\sigma_3^2 (x_1) - \sigma_1^2 (x_1) &= O(d_{1n}^2) + O(d_{0n}^2) + O(d_{1n} d_{0n}).
\end{align*}
These differences show that when only the propensity score function is  or only the outcome regression functions are locally misspecified, the asymptotic distribution remains the same as that without misspecification. \hfill$\Box$

\subsection{Proofs of Theorems 5 and 6}

Consider the cases with all models misspecified. The proof of Theorem 5 will be very similar to the proof of scenario 1 in Theorem 6 except that the asymptotic linear expression can be as
\begin{align*}
	\sqrt{nh_1^k} \left[ \widehat{\tau}(x_1)-\tau(x_1) \right] =& \frac{1}{\sqrt{nh_1^k}}\frac{1}{f(x_1)} \sum_{i=1}^n [\Psi_4(X_i,Y_i,D_i)-\tau(x_1)] K_1 \left( \frac{X_{1i}-x_1}{h_1} \right)\\
& + o_p(1).
\end{align*}
As the unbiasedness no longer holds, we then compute the bias term. A decomposition is as follows:
\begin{align*}
	& E \left\{ \sqrt{n h_1^k} \left[ \widehat{\tau}(x_1) - \tau(x_1) \right] \right\} \\
	= & \sqrt{n h_1^k} E \left\{ \left( \frac{D}{\widetilde{p}(X;\beta^*)} - \frac{D}{p(X)} \right) [Y - m_1(X)] + \left( 1 - \frac{D}{\widetilde{p}(X;\beta^*)} \right) [\widetilde{m}_1(X;\gamma_1^*) - m_1(X)] \right. \\
	& \left. \left. - \left( \frac{D}{\widetilde{p}(X;\beta^*)} - \frac{D}{p(X)} \right) [Y - m_0(X)] + \left( 1 - \frac{D}{\widetilde{p}(X;\beta^*)} \right) [\widetilde{m}_0(X;\gamma_0^*) - m_0(X)] \right| X_1 = x_1 \right\}\\
	= &  \sqrt{n h_1^k} E \left\{ \frac{\left[ m_1(X) - \widetilde{m}_1(X;\gamma_1^*) \right] \left[ p(X) - \widetilde{p}(X;\beta^*) \right]}{\widetilde{p}(X;\beta^*)} \right. \\
	 & \left. \left. - \frac{\left[ m_0(X) - \widetilde{m}_0(X;\gamma_0^*) \right] \left[ \widetilde{p}(X;\beta^*) - p(X) \right]}{1 - \widetilde{p}(X;\beta^*)} \right| X_1 = x_1 \right\} \\
	 := \, \, & \sqrt{n h_1^k} \, bias(x_1).
\end{align*}
Let
\begin{align*}
	\widetilde{\tau}(x_1) = & E \left\{ \left[ \frac{D}{\widetilde{p}(X;\beta^*)}[Y - \widetilde{m}_1(X;\gamma_1^*)]  \right. \right.  \\
	& -\left. \left. \left. \frac{1 - D}{1 - \widetilde{p}(X;\beta^*)}[Y - \widetilde{m}_0(X;\gamma_0^*)] + \widetilde{m}_1(X;\gamma_1^*) - \widetilde{m}_0(X;\gamma_0^*) \right] \right| X_1 = x_1 \right\}.
\end{align*}
The variance term of $\sqrt{nh_1^k} \left[ \widehat{\tau}(x_1)-\tau(x_1) - bias(x_1) \right]$ can be derived as:
\begin{align*}
	& Var \left\{ \frac{1}{\sqrt{n h_1^k}} \frac{1}{f(x_1)} \sum_{i=1}^n [\Psi_4(X_i, Y_i, D_i) - \tau(x_1) - bias(x_1)] K_1 \left( \frac{X_{1i} - x_1}{h_1} \right) \right\} \\
	= & E \left\{ \frac{1}{\sqrt{n h_1^k}} \frac{1}{f(x_1)} \sum_{i=1}^n [\Psi_4(X_i, Y_i, D_i) - \widetilde{\tau}(x_1)] K_1 \left( \frac{X_{1i} - x_1}{h_1} \right) \right\}^2 \\
	= & \frac{h_1^k}{f^2(x_1)} E \left\{ E \left[ \left. \left[ [\Psi_4(X, Y, D)- \widetilde{\tau}(x_1)] \frac{1}{h_1^k} K_1 \left( \frac{X_1-x_1}{h_1} \right) \right]^2 \right| X_1 \right] \right\} \\
	= & \frac{h_1^k}{f^2(x_1)} \frac{1}{h_1^k} \int K_1^2(u) E \left[ \left.  [\Psi_4(X, Y, D)- \widetilde{\tau}(x_1)]^2 \right| X_1=x_1 + h_1 u \right] f(x_1 + h_1 u) du \\
	= & \frac{\sigma_4^2 (x_1) \int K_1^2(u) du}{f(x_1)} + O(h_1^k).
\end{align*}
where
\begin{align*}
	\sigma_4^2 (x_1) = E \left[ \left.  [\Psi_4(X, Y, D)- \widetilde{\tau}(x_1)]^2 \right| X_1=x_1 \right].
\end{align*}
With the same argument to derive the asymptotic distribution in Theorem 1, we can obtain that
\begin{align}
	\sqrt{nh_1^k} \left[ \widehat{\tau}(x_{1})-\tau(x_{1}) - bias(x_1) \right] \xrightarrow{d} N \left( 0,\frac{\sigma_4^2 (x_1) \int K_1^2(u) du}{f(x_1)}\right). \label{allgm1}
\end{align}
The proof of Theorem 5 is completed.   \hfill$\Box$

Note that Theorem 6 is a variant of Theorem 5. To derive the asymptotic distribution, we only need to consider the bias term and the variance term based on (\ref{allgm1}) when all nuisance models are locally misspecified. From the definitions of misspecified models before, the bias term can be bounded by
\begin{align*}
	bias(x_1) = O_p(c_n d_{1n}) + O_p(c_n d_{0n}).
\end{align*}
This result implies that if  the convergence rates of $c_n d_{1n}$ and $c_n d_{1n}$ are faster than $O \left( \frac{1}{\sqrt{nh_1^k}} \right)$, the bias term vanishes asymptotically. By the central limit theorem, we can also derive the asymptotic normality with the variance term. By (\ref{allgm1}) and (\ref{allgm1}) when $c_n$, $d_{1n}$ and $d_{0n}$ all converge to $0$, we have
\begin{align*}
	\frac{\sigma_4^2 (x_1) \int K_1^2(u) du}{f(x_1)} = \frac{\sigma_1^2 (x_1) \int K_1^2(u) du}{f(x_1)} + o(1).
\end{align*}
With Slutsky's Theorem, we can conclude that, when all nuisance models are locally misspecified,  $c_n d_{1n} = o\left( \frac{1}{\sqrt{n h_1^k}} \right)$, and $c_n d_{0n} = o\left( \frac{1}{\sqrt{n h_1^k}} \right)$, we then have
\begin{align*}
	\sqrt{nh_1^k}[\widehat{\tau}(x_1) - \tau(x_1)] \xrightarrow{d} \mathcal{N}\left( 0, \frac{\sigma_1^2 (x_1) \int K_1^2(u) du}{f(x_1)} \right).
\end{align*}
Then the proof is completed. \hfill$\Box$

\subsection{A simple justification for  Remark~5}
\setcounter{equation}{0}

As we showed in the proof of Theorem~4, 
\begin{align*}
	\sigma_2^2 (x_1) - \sigma_1^2 (x_1) = & E \left\{ \left. \frac{p(X)-\widetilde{p}(X;\beta^*)}{\left[ \widetilde{p}(X;\beta^*) \right]^2} Var(Y|D=1,X) \right| X_1 = x_1 \right\} \\
	& + E \left\{ \left. \frac{\widetilde{p}(X;\beta^*)-p(X)}{\left[ 1-\widetilde{p}(X;\beta^*) \right]^2} Var(Y|D=0,X) \right| X_1 = x_1 \right\}.
\end{align*}
This difference cannot be showed either positive or negative for all $x_1$. The example in Remark~2 confirms this.  For $\sigma_3^2 (x_1)$ we have
\begin{align*}
	& \sigma_3^2 (x_1) - \sigma_1^2 (x_1) \\
	= & E \left\{ \left. \left[ \left( 1-\frac{D}{p(X)} \right) \left[ \widetilde{m}_1(X;\gamma_1^*)-m_1(X) \right] - \left( 1-\frac{1-D}{1-p(X)}\right) \left[ \widetilde{m}_0(X;\gamma_0^*) - m_0(X) \right]  \right]^2 \right| X_1 = x_1 \right\} \\
	&\ge  \, 0.
\end{align*}
In other words, the variance with $\sigma_3^2 (x_1)$ can be larger than that of the estimators with all models correctly specified. Further,
\begin{align*}\label{allgm2}
	& \sigma_4^2(x_1) - \sigma_1^2(x_1) \\
	= & Var(\Psi_4(X,Y,D)|X_1 = x_1) - Var(\Psi_1(X,Y,D)|X_1 = x_1) \\
	= & E \left\{ \left. \left( \frac{p(X)}{\widetilde{p}^2(X;\beta^*)} - \frac{1}{p(X)} \right) Var(Y|X,D=1) \right| X_1 = x_1 \right\} \\
	& - E \left\{ \left. \left( \frac{1 - p(X)}{[1 - \widetilde{p}(X;\beta^*)]^2} - \frac{1}{1-p(X)} \right) Var(Y|X,D=0) \right| X_1 = x_1 \right\} \\
	& + E \left\{ \left. \frac{p(X)}{[\widetilde{p}(X;\beta^*)]^2} [m_1(X) - \widetilde{m}_1(X;\gamma_1^*)]^2 + \frac{1-p(X)}{[1-\widetilde{p}(X;\beta^*)]^2} [m_0(X) - \widetilde{m}_0(X;\gamma_0^*)]^2 \right| X_1 = x_1 \right\} \\
	& + 2 E \left\{ \left. [\widetilde{m}_1(X;\gamma_1^*) - \widetilde{m}_0(X;\gamma_0^*)] [m_1(X) - m_0(X) - \widetilde{m}_1(X;\gamma_1^*) + \widetilde{m}_0(X;\gamma_0^*)] \right| X_1 = x_1 \right\} \notag \\
	& + \tau^2(x_1) - \widetilde{\tau}^2(x_1).
\end{align*}
Again $\sigma_4^2(x_1)$ cannot be easily judged whether it is larger than $\sigma_1^2(x_1)$ or not. \hfill$\Box$

\newpage

\subsection{Additional Simulation Results}

\begin{table}[h!]
\caption{The simulation results under model 2 (part 1)}
\label{result 1}
\resizebox{\linewidth}{!}{\begin{tabular}{|c|r|rrrrr|rrrrr|}
\hline
 & \multicolumn{1}{c}{}   & \multicolumn{5}{c}{n=500} & \multicolumn{5}{c}{n=5000} \\
\cline{3-12}
DRCATE & \multicolumn{1}{c}{$x_1$} & \multicolumn{1}{c}{bias} & \multicolumn{1}{c}{sam-SD} & \multicolumn{1}{c}{MSE} & \multicolumn{1}{c}{$P_{0.05}$} & \multicolumn{1}{c}{$P_{0.95}$} & \multicolumn{1}{c}{bias} & \multicolumn{1}{c}{sam-SD} & \multicolumn{1}{c}{MSE} & \multicolumn{1}{c}{$P_{0.05}$} & \multicolumn{1}{c}{$P_{0.95}$} \\
\hline
\multirow{5}{*}{\begin{tabular}[c]{@{}c@{}}DRCATE\\ (O,O)\end{tabular}} & -0.4 & 0.0012 & 0.2077 & 0.0431 & 0.058 & 0.045 & 0.0006 & 0.2034 & 0.0414                  & 0.045 & 0.047 \\
 & -0.2 & 0.0010 & 0.2139 & 0.0457 & 0.040 & 0.052 & 0.0001 & 0.1988 & 0.0395 & 0.050 & 0.050 \\
 & 0 & -0.0005 & 0.2046 & 0.0418 & 0.048 & 0.059 & 0.0007 & 0.1846 & 0.0341 & 0.050 & 0.050 \\
 & 0.2 & -0.0001 & 0.2226 & 0.0495 & 0.045 & 0.050 & 0.0008 & 0.2034 & 0.0415 & 0.038 & 0.057 \\
 & 0.4 & 0.0024 & 0.3312 & 0.1097 & 0.048 & 0.049 & 0.0010 & 0.3114 & 0.0971 & 0.044 & 0.053 \\
\hline
\multirow{5}{*}{\begin{tabular}[c]{@{}c@{}}DRCATE\\ (cP,cP)\end{tabular}} & -0.4 & 0.0011 & 0.2077 & 0.0431 & 0.056 & 0.048 & 0.0006 & 0.2035 & 0.0415                  & 0.046 & 0.046 \\
 & -0.2 & 0.0009 & 0.2137 & 0.0456 & 0.045 & 0.053 & 0.0001 & 0.1988 & 0.0395 & 0.049 & 0.050 \\
 & 0 & -0.0006 & 0.2044 & 0.0417 & 0.047 & 0.055 & 0.0007 & 0.1846 & 0.0341 & 0.048 & 0.052 \\
 & 0.2 & -0.0001 & 0.2228 & 0.0496 & 0.046 & 0.049 & 0.0007 & 0.2035 & 0.0415 & 0.038 & 0.057 \\
 & 0.4 & 0.0024 & 0.3316 & 0.1100 & 0.047 & 0.047 & 0.0010 & 0.3114 & 0.0971 & 0.044 & 0.053 \\
\hline
\multirow{5}{*}{\begin{tabular}[c]{@{}c@{}}DRCATE\\ (N,N)\end{tabular}} & -0.4 & 0.0002 & 0.2653 & 0.0703 & 0.017 & 0.029 & 0.0004 & 0.2136 & 0.0456                  & 0.057 & 0.052 \\
 & -0.2 & 0.0011 & 0.2300 & 0.0529 & 0.042 & 0.048 & 0.0004 & 0.1990 & 0.0396 & 0.041 & 0.045 \\
 & 0 & 0.0007 & 0.1962 & 0.0385 & 0.048 & 0.051 & 0.0003 & 0.1917 & 0.0367 & 0.041 & 0.052 \\
 & 0.2 & 0.0011 & 0.2299 & 0.0528 & 0.043 & 0.058 & 0.0006 & 0.2122 & 0.0451 & 0.046 & 0.052 \\
 & 0.4 & 0.0041 & 0.3373 & 0.1141 & 0.054 & 0.057 & 0.0003 & 0.3125 & 0.0976 & 0.050 & 0.052 \\
\hline
\multirow{5}{*}{\begin{tabular}[c]{@{}c@{}}DRCATE\\ (S,S)\end{tabular}} & -0.4 & -0.0018 & 0.2058 & 0.0424 & 0.051 & 0.046 & 0.0002 & 0.2501 & 0.0625                  & 0.028 & 0.040 \\
 & -0.2 & -0.0021 & 0.2093 & 0.0439 & 0.056 & 0.039 & -0.0008 & 0.2087 & 0.0436 & 0.046 & 0.047 \\
 & 0 & 0.0000 & 0.2040 & 0.0416 & 0.055 & 0.051 & 0.0011 & 0.1868 & 0.0351 & 0.044 & 0.056 \\
 & 0.2 & 0.0060 & 0.2257 & 0.0518 & 0.031 & 0.068 & 0.0014 & 0.2093 & 0.0441 & 0.047 & 0.059 \\
 & 0.4 & 0.0089 & 0.3409 & 0.1181 & 0.039 & 0.064 & 0.0010 & 0.3298 & 0.1089 & 0.043 & 0.062 \\
\hline
\end{tabular}}
\end{table}

\newpage

\begin{table}[h!]
\caption{The simulation results under model 2 (part 2)}
\label{result 2}
\resizebox{\linewidth}{!}{\begin{tabular}{|c|r|rrrrr|rrrrr|}
\hline
 & \multicolumn{1}{c}{}   & \multicolumn{5}{c}{n=500} & \multicolumn{5}{c}{n=5000} \\
\cline{3-12}
DRCATE & \multicolumn{1}{c}{$x_1$} & \multicolumn{1}{c}{bias} & \multicolumn{1}{c}{sam-SD} & \multicolumn{1}{c}{MSE} & \multicolumn{1}{c}{$P_{0.05}$} & \multicolumn{1}{c}{$P_{0.95}$} & \multicolumn{1}{c}{bias} & \multicolumn{1}{c}{sam-SD} & \multicolumn{1}{c}{MSE} & \multicolumn{1}{c}{$P_{0.05}$} & \multicolumn{1}{c}{$P_{0.95}$} \\
\hline
\multirow{5}{*}{\begin{tabular}[c]{@{}c@{}}DRCATE\\ (O,O)\end{tabular}} & -0.4 & 0.0012 & 0.2077 & 0.0431 & 0.058 & 0.045 & 0.0006 & 0.2034 & 0.0414                  & 0.045 & 0.047 \\
 & -0.2 & 0.0010 & 0.2139 & 0.0457 & 0.040 & 0.052 & 0.0001 & 0.1988 & 0.0395 & 0.050 & 0.050 \\
 & 0 & -0.0005 & 0.2046 & 0.0418 & 0.048 & 0.059 & 0.0007 & 0.1846 & 0.0341 & 0.050 & 0.050 \\
 & 0.2 & -0.0001 & 0.2226 & 0.0495 & 0.045 & 0.050 & 0.0008 & 0.2034 & 0.0415 & 0.038 & 0.057 \\
 & 0.4 & 0.0024 & 0.3312 & 0.1097 & 0.048 & 0.049 & 0.0010 & 0.3114 & 0.0971 & 0.044 & 0.053 \\
\hline
\multirow{5}{*}{\begin{tabular}[c]{@{}c@{}}DRCATE\\ (mP,cP)\end{tabular}} & -0.4 & 0.0011 & 0.2082 & 0.0433 & 0.058 & 0.041 & 0.0006 & 0.2042 & 0.0417                  & 0.050 & 0.045 \\
 & -0.2 & 0.0009 & 0.2123 & 0.0451 & 0.044 & 0.054 & 0.0001 & 0.1974 & 0.0389 & 0.051 & 0.053 \\
 & 0 & -0.0005 & 0.2025 & 0.0410 & 0.048 & 0.058 & 0.0006 & 0.1834 & 0.0337 & 0.050 & 0.052 \\
 & 0.2 & -0.0002 & 0.2222 & 0.0493 & 0.045 & 0.052 & 0.0007 & 0.2030 & 0.0413 & 0.037 & 0.056 \\
 & 0.4 & 0.0025 & 0.3315 & 0.1099 & 0.047 & 0.051 & 0.0011 & 0.3116 & 0.0972 & 0.043 & 0.053 \\
\hline
\multirow{5}{*}{\begin{tabular}[c]{@{}c@{}}DRCATE\\ (mP,N)\end{tabular}} & -0.4 & -0.0011 & 0.2156 & 0.0464 & 0.056 & 0.042 & -0.0005 & 0.2082 & 0.0434 & 0.048 & 0.043 \\
 & -0.2 & -0.0019 & 0.2086 & 0.0436 & 0.061 & 0.036 & -0.0013 & 0.2062 & 0.0428 & 0.058 & 0.036 \\
 & 0 & -0.0028 & 0.2003 & 0.0403 & 0.057 & 0.034 & -0.0011 & 0.1888 & 0.0358 & 0.058 & 0.040 \\
 & 0.2 & 0.0021 & 0.2108 & 0.0445 & 0.052 & 0.058 & -0.0004 & 0.2099 & 0.0440 & 0.054 & 0.044 \\
 & 0.4 & 0.0060 & 0.3258 & 0.1069 & 0.045 & 0.059 & 0.0033 & 0.3276 & 0.1093 & 0.033 & 0.069 \\
\hline
\multirow{5}{*}{\begin{tabular}[c]{@{}c@{}}DRCATE\\ (mP,S)\end{tabular}} & -0.4 & -0.0034 & 0.2215 & 0.0493 & 0.054 & 0.050 & -0.0010 & 0.2119 & 0.0451 & 0.053 & 0.043 \\
 & -0.2 & -0.0055 & 0.2235 & 0.0507 & 0.060 & 0.041 & -0.0034 & 0.2115 & 0.0469 & 0.073 & 0.029 \\
 & 0 & -0.0023 & 0.2049 & 0.0421 & 0.051 & 0.043 & -0.0025 & 0.1895 & 0.0371 & 0.061 & 0.032 \\
 & 0.2 & -0.0003 & 0.2149 & 0.0462 & 0.043 & 0.045 & -0.0003 & 0.1982 & 0.0393 & 0.052 & 0.043 \\
 & 0.4 & 0.0102 & 0.3351 & 0.1148 & 0.032 & 0.068 & 0.0034 & 0.3122 & 0.0997 & 0.039 & 0.058 \\
\hline
\end{tabular}}
\end{table}

\newpage

\begin{table}[h!]
\caption{The simulation results under model 2 (part 3)}
\label{result 3}
\resizebox{\linewidth}{!}{\begin{tabular}{|c|r|rrrrr|rrrrr|}
\hline
 & \multicolumn{1}{c}{}   & \multicolumn{5}{c}{n=500} & \multicolumn{5}{c}{n=5000} \\
\cline{3-12}
DRCATE & \multicolumn{1}{c}{$x_1$} & \multicolumn{1}{c}{bias} & \multicolumn{1}{c}{sam-SD} & \multicolumn{1}{c}{MSE} & \multicolumn{1}{c}{$P_{0.05}$} & \multicolumn{1}{c}{$P_{0.95}$} & \multicolumn{1}{c}{bias} & \multicolumn{1}{c}{sam-SD} & \multicolumn{1}{c}{MSE} & \multicolumn{1}{c}{$P_{0.05}$} & \multicolumn{1}{c}{$P_{0.95}$} \\
\hline
\multirow{5}{*}{\begin{tabular}[c]{@{}c@{}}DRCATE\\ (O,O)\end{tabular}} & -0.4 & 0.0012 & 0.2077 & 0.0431 & 0.058 & 0.045 & 0.0006 & 0.2034 & 0.0414                  & 0.045 & 0.047 \\
 & -0.2 & 0.0010 & 0.2139 & 0.0457 & 0.040 & 0.052 & 0.0001 & 0.1988 & 0.0395 & 0.050 & 0.050 \\
 & 0 & -0.0005 & 0.2046 & 0.0418 & 0.048 & 0.059 & 0.0007 & 0.1846 & 0.0341 & 0.050 & 0.050 \\
 & 0.2 & -0.0001 & 0.2226 & 0.0495 & 0.045 & 0.050 & 0.0008 & 0.2034 & 0.0415 & 0.038 & 0.057 \\
 & 0.4 & 0.0024 & 0.3312 & 0.1097 & 0.048 & 0.049 & 0.0010 & 0.3114 & 0.0971 & 0.044 & 0.053 \\
\hline
\multirow{5}{*}{\begin{tabular}[c]{@{}c@{}}DRCATE\\ (cP,mP)\end{tabular}} & -0.4 & 0.0008 & 0.2179 & 0.0474 & 0.051 & 0.043 & 0.0004 & 0.2204 & 0.0485                  & 0.048 & 0.045 \\
 & -0.2 & 0.0012 & 0.2233 & 0.0498 & 0.049 & 0.052 & 0.0003 & 0.2069 & 0.0428 & 0.049 & 0.054 \\
 & 0 & -0.0004 & 0.2104 & 0.0442 & 0.051 & 0.060 & 0.0008 & 0.1890 & 0.0358 & 0.050 & 0.053 \\
 & 0.2 & -0.0002 & 0.2226 & 0.0495 & 0.048 & 0.049 & 0.0007 & 0.2039 & 0.0416 & 0.036 & 0.054 \\
 & 0.4 & 0.0028 & 0.3407 & 0.1162 & 0.047 & 0.050 & 0.0010 & 0.3170 & 0.1006 & 0.043 & 0.053 \\
\hline
\multirow{5}{*}{\begin{tabular}[c]{@{}c@{}}DRCATE\\ (N,mP)\end{tabular}} & -0.4 & -0.0050 & 0.2225 & 0.0501 & 0.060 & 0.036 & -0.0006 & 0.2227 & 0.0496                  & 0.051 & 0.050 \\
 & -0.2 & -0.0015 & 0.2185 & 0.0477 & 0.056 & 0.043 & -0.0011 & 0.1931 & 0.0375 & 0.054 & 0.039 \\
 & 0 & -0.0020 & 0.2039 & 0.0416 & 0.072 & 0.032 & -0.0013 & 0.1857 & 0.0348 & 0.056 & 0.038 \\
 & 0.2 & 0.0024 & 0.2178 & 0.0475 & 0.042 & 0.051 & 0.0005 & 0.2064 & 0.0426 & 0.039 & 0.056 \\
 & 0.4 & 0.0046 & 0.3259 & 0.1066 & 0.035 & 0.064 & 0.0023 & 0.3324 & 0.1115 & 0.044 & 0.050 \\
\hline
\multirow{5}{*}{\begin{tabular}[c]{@{}c@{}}DRCATE\\ (S,mP)\end{tabular}} & -0.4 & -0.0115 & 0.2117 & 0.0481 & 0.075 & 0.027 & -0.0024 & 0.3260 & 0.1073                  & 0.020 & 0.017 \\
 & -0.2 & -0.0021 & 0.2083 & 0.0434 & 0.051 & 0.051 & -0.0021 & 0.2010 & 0.0412 & 0.065 & 0.033 \\
 & 0 & -0.0018 & 0.2002 & 0.0401 & 0.044 & 0.053 & -0.0005 & 0.2045 & 0.0418 & 0.045 & 0.038 \\
 & 0.2 & 0.0035 & 0.2290 & 0.0527 & 0.044 & 0.071 & 0.0001 & 0.2155 & 0.0464 & 0.054 & 0.054 \\
 & 0.4 & 0.0017 & 0.3460 & 0.1196 & 0.040 & 0.064 & 0.0015 & 0.3519 & 0.1241 & 0.031 & 0.054 \\
\hline
\end{tabular}}
\end{table}

\newpage


\bibhang=1.7pc
\bibsep=2pt
\fontsize{9}{14pt plus.8pt minus .6pt}\selectfont
\renewcommand\bibname{\large \bf References}
\expandafter\ifx\csname
\natexlab\endcsname\relax\def\natexlab#1{#1}\fi
\expandafter\ifx\csname url\endcsname\relax
  \def\url#1{\texttt{#1}}\fi
\expandafter\ifx\csname urlprefix\endcsname\relax\def\urlprefix{URL}\fi

\bibliographystyle{chicago}
\bibliography{ref20200419}

\end{document}